\documentclass[preprint,3p,times,numbers]{elsarticle}
\biboptions{sort&compress}

\usepackage{amsmath}
\usepackage{amssymb}
\usepackage{bm}   
\usepackage{float}
\usepackage{graphicx}
\usepackage{caption}
\usepackage{subcaption}
\usepackage{placeins}
\usepackage{comment}
\usepackage{multirow}
\usepackage{lipsum}
\usepackage{booktabs}

\usepackage{xcolor}

\usepackage{lineno} 

\usepackage{setspace}

\journal{International Journal of Mechanical Sciences}

\begin{document}

\begin{frontmatter}

\title{ {Phase-field digital image correlation for integrated displacement and damage measurements}}

\author{Dingxiang Zhu}

\author{Ye Lu\corref{cor1}}
\ead{yelu@umbc.edu}
\cortext[cor1]{Corresponding author}

\affiliation{
  organization = {Department of Mechanical Engineering, University of Maryland Baltimore County},
  city         = {Baltimore},
  country      = {USA}
}

\onehalfspacing 

\begin{abstract}
This work presents a novel digital image correlation (DIC) framework for full-field measurements of displacement, strain, and damage, based on a phase field (PF) approach. The idea is to take advantage of the ability of the PF method to track complex crack morphologies and to provide a natural way in DIC to perform damage and crack measurements from experimental speckle images, in addition to displacement and strain fields. Moreover, incorporating the damage variable into DIC can improve the displacement accuracy near the crack tip, and can avoid the need of user-defined masks when dealing with cracked samples, which is advantageous when cracks become complex and the manual application of masks becomes challenging. The theoretical formulation of the proposed framework, namely PF-DIC, was presented in detail in the paper, along with a finite element implementation.  Numerical examples have demonstrated the capability of the proposed PF-DIC in terms of capturing different types of cracks while  {providing similar measurement accuracy to that of masked DIC}. Additionally, it is shown that the PF-DIC can be easily adapted to selectively identify critical cracks under specific loading conditions or mechanisms for damage assessment and diagnostic purposes. The proposed DIC framework  {can be used to characterize material defects, support structural health monitoring, and enable a potential unification of PF simulations and experimental fracture measurements}.

\end{abstract}




\end{frontmatter}



\onehalfspacing 
\section{Introduction}

Digital image correlation (DIC) has become an important and widely used optical technique for displacement and strain measurements in mechanical testing due to its non-contact nature,  high accuracy, and ability to capture complex full-field displacement and strain data  {\cite{Peters1982,Hild2006,Pan2018}}. Among its broad range of applications in science and engineering, DIC has enabled quantitative tracking of deformation evolution \cite{Adam2013}, high-temperature strain determination \cite{Pan2010},  dynamic testing \cite{Sasso2018}, fracture analysis  {\cite{Tian2024,huo2022characterization}}, and comprehensive material characterization \cite{bouclier2026pinn}. In particular, for fracture analysis, DIC provides essential experimental data to characterize the fracture behavior of materials \cite{Huo2021,Wu2022,kunecky2025mode}, identify the underlying mechanisms \cite{heinzmann2024investigation}, and validate theoretical and numerical predictions \cite{mishra2026fracture}. Furthermore, a recent review on outstanding issues and emerging frontiers in fracture mechanics \cite{yang2026outstanding} has highlighted the growing need to integrate experimental and simulation data for predictive fracture modeling in real materials. In this context, DIC emerges as a unique and highly valuable tool for bridging experiments and simulations.

Various DIC algorithms exist and can be broadly categorized into two groups: local DIC, commonly referred to as subset-based DIC, and global DIC, also known as finite element-based DIC (FE-DIC). The  details and comparisons of the various methods can be found in \cite{Bruck1989,Hild2012,Wang2016}. Generally speaking, in local DIC, the region of interest is partitioned into multiple subsets, and the displacement in each subset is determined separately, without explicitly considering the continuity between adjacent subsets  {\cite{Yang2021}}. By contrast, global DIC  {\cite{chapelier2021free,elguedj2011isogeometric,Lu2024}} describes the displacement field in a continuous manner over the entire region of interest and incorporates FE approximations  into the correlation process. Nevertheless, the standard DIC frameworks (both local and global versions) fundamentally rely on the assumption that the deformation is continuous in the region of interest or at least in the subsets. When discontinuity  occurs, such as in the presence of cracks, this assumption can result in substantial matching errors or even correlation failure. Consequently, the ability of conventional DIC to deal with complex cracks remains limited \cite{BECKER2012,VascoOlmo2019,Hassan2021}. 

To address the challenge of handling cracks in DIC, a common practice is to either exclude image subsets/pixels  \cite{Helm2008,Zhu2023,Chen2023} near discontinuities from the correlation analysis or remove the erroneous values after the correlation in a post-processing stage \cite{Gehri2020,Gehri2022,Seemab2023}. These approaches usually require  extensive manual intervention to determine suitable masks for complex crack patterns, and no unified automatic strategy is available to accurately recover crack paths and morphologies directly from the displacement measurements  {\cite{Vincens2024,Zhao2025}}. An alternative approach is the subset-splitting method \cite{Poissant2009}, which divides a crack-affected subset into two parts along the crack trajectory and performs separate correlation analyses on each part to avoid the erroneous correlation. However, the effectiveness of the method still strongly depends on the accurate identification of the split line location, which is often difficult to obtain in practice  {\cite{Han2018,Dong2025}}, although some improvement strategies have been proposed  {\cite{Pan2015,Zhang2024,Becker2023}}. Another important work that could resolve the discontinuity issue with cracks is the extended DIC (X-DIC) method \cite{Rthor2007} that incorporates the concept of extended finite element method into global DIC. Unlike standard FE-DIC, the X-DIC method can accommodate displacement discontinuity with enriched FE approximations \cite{Rthor2007(2),Rthor2008,Nguyen2011}, although the enrichment functions may require prior knowledge of the crack geometry.  More recently, machine learning-based approaches have also been explored to deal with crack detection and fracture characterization in DIC  {\cite{strohmann2021automatic,holzmond2024enhancing,vsvcerba2025inflection}}. For example,  {Sadeghian et al. \cite{Sadeghian2025} integrated machine learning and DIC to assess damage in engineering materials.} Wang et al. \cite{wang2022digital} proposed a deep learning led semantic segmentation algorithm to
automatically track changes in DIC-measured strain fields and extract damage information in carbon fiber reinforced plastic  composites.  Leu et al. \cite{Leu2026} used machine learning to automate mask generation for DIC analysis of specimens with complex geometries. Despite these advances, machine learning-based methods still face significant limitations, such as strong dependence on training data and potentially high computational costs due to data preparation and model training.

In this work, we proposed an alternative approach to deal with this challenge. Inspired by the phase field simulations of fracture \cite{Francfort1998,Miehe2010,Lu2020}, we proposed the so-called phase field DIC (PF-DIC) framework that integrates an implicit representation of cracks into the DIC analysis. As a result, the  PF-DIC is expected to benefit from the implicit crack   tracking, natural treatment of crack initiation and propagation, and the ability to easily handle complex crack topologies. Similarly to the fracture modeling, the specific formulation of the PF-DIC will be presented as a variational problem that minimizes an energy functional unifying the grayscale correlation and the damage analysis. Importantly, the framework allows a flexible definition of damage driving force and offers great advantages in terms of damage diagnosis and potential coupling with various multiphysics damage mechanisms. Furthermore, due to the involvement of a damage variable, the PF-DIC eliminates the need for a pre-defined mask and provides a damage field as an additional measurement closely integrated with displacement and strain fields from the image correlation process. It was found that the PF-DIC can effectively detect various types of cracks under different loading conditions, such as tensile and bending tests,  track potential crack propagation, and selectively identify critical cracks for diagnostic purposes. Compared to the standard FE-DIC, the PF-DIC can improve the measurement accuracy near the crack tip without requiring a pre-defined mask. These were confirmed by numerical examples using both synthetic and experimental speckle images of fractured samples.

The proposed PF-DIC opens many opportunities for both experimental characterization and simulations of fracture of materials. For example, its automatic damage detection capability enables in-situ monitoring of crack propagation and crack-tip displacement without the need for manual mask application. Furthermore, the resulting experimental data can be directly compared with and potentially integrated into PF fracture simulations due to the unified representation of cracks.  {This would be beneficial, as PF simulations have been widely adopted for fracture and fatigue modeling \cite{wu2020phase,li2023review,yang2026fatigue} for a broad range of materials, such as composites \cite{bui2021review}, concrete \cite{feng2018phase,fang2026multiscale,meng2026multi}, rocks \cite{zhou2018phase,zhang2017modification}, ceramics \cite{carollo2018modeling}, biomaterials \cite{wu2020fracture}, quasicrystals \cite{li2023phase}, polymers \cite{miehe2014phase,konica2021reaction}, and alloys \cite{sharma2026phase,zhang2026phase,castro2026coupled}. Leveraging the experimental data could potentially improve the predictive capability of simulation models and help extract important fracture properties of materials.}

The remainder of the paper is organized as follows. Section 2 introduces the proposed PF-DIC framework, including the conservation of optical flow, the diffusive representation of cracks, the PF-DIC formulation, as well as the corresponding solution procedure and implementation details. Section 3 presents numerical examples, including both synthetic and experimental speckle images, to demonstrate the performance of PF-DIC. Finally, the paper closes with some concluding remarks.

\label{introduction}

\section{The proposed PF-DIC framework}

\subsection{Conservation of optical flow}
We start from introducing the basic idea of DIC. Generally speaking, the conventional DIC analysis is to extract the displacement fields from images, using the concept of conservation of the optical flow  {\cite{Hild2006}}. Considering two images of the same specimen with a speckle pattern printed on surfaces, we can designate one as the reference image, with $f(\boldsymbol{X})\in [0,1]$ denoting its grayscale value of the position $\boldsymbol{X}$, and another as the deformed image, with $g(\boldsymbol{x})\in [0,1]$ denoting the grayscale of the position  $\boldsymbol{x}$. It is clear that if  $\boldsymbol{X}$ and $\boldsymbol{x}$ are taken as the initial and final positions of the same material point on the images,  we should have: $g(\boldsymbol{x})=f(\boldsymbol{X})$, i.e., the conservation of the optical flow. Therefore, given that $\boldsymbol{X}=\boldsymbol{x}-\boldsymbol{u}$, with $\boldsymbol{u}$ denoting the displacement of the point, we have  
\begin{equation}
g({\boldsymbol{x}}) = f({\boldsymbol{x}} - \boldsymbol{u}),
\label{eq:optical_flow}
\end{equation}
which provides a fundamental equation to solve for the unknown displacement $\boldsymbol{u}$. In practice, the equity between $f$ and $g$ can  be approximated by minimizing the the following residual 
\begin{equation}
r(\boldsymbol{u}) = g(\boldsymbol{x}) - f(\boldsymbol{x} - \boldsymbol{u}).
\label{eq:residual}
\end{equation}

Considering that the residual should be minimized over the entire domain of the images, the best approximation of the displacement field $\boldsymbol{u}$ can be found by solving the following minimization problem
\begin{equation}
\min \, \int_{\mathit{\Omega}} r^2(\boldsymbol{u}) \,\mathrm{d}V
= \int_{\mathit{\Omega}} \left( g(\boldsymbol{x}) - f(\boldsymbol{x} - \boldsymbol{u}) \right)^2 \, \mathrm{d}V,
\label{eq:minimization}
\end{equation}
where $\mathit{\Omega}$ denotes the region of interest (ROI) in the images. With an appropriate spatial discretization, this nonlinear optimization problem can be solved either by  an iterative scheme, e.g., the steepest descent \cite{Long2012}, Newton–Raphson \cite{Bruck1989}, Gauss–Newton \cite{Black1998},  Levenberg–Marquardt \cite{Cheng2002}, and inverse compositional Guass-Newton methods \cite{Passieux2019}, or by a direct method using linearization techniques \cite{Besnard2006,Lu2024}. To deal with problems involving fractured samples, a common approach in DIC is to apply a mask to exclude the cracking areas in the computational domain \cite{Chen2023,Tian2024}, i.e.,
\begin{equation}
\min \, \int_{\mathit{\Omega}\setminus \mathit{\Omega}_{\text{mask}}} r^2(\boldsymbol{u}) \,\mathrm{d}V
= \int_{\mathit{\Omega}\setminus \mathit{\Omega}_{\text{mask}}} \left( g(\boldsymbol{x}) - f(\boldsymbol{x} - \boldsymbol{u}) \right)^2 \, \mathrm{d}V,
\label{eq:minimization_crack}
\end{equation}
where $\mathit{\Omega}_{\text{mask}}$ denotes a pre-defined mask containing the pixels to be removed from the DIC images. This approach is effective but needs an explicit description of the mask, which can be challenging when the cracks become complex and requires extensive human intervention to identify the cracks from images.  In the proposed PF-DIC, we adopt a PF approach to describe the cracks in an implicit way.

\subsection{Diffusive representation of cracks}
\label{sec:diff_crack}
\begin{figure}[htbp]
  \centering
  \begin{subfigure}{0.25\textwidth}
    \centering
    \includegraphics[height=3.5cm]{   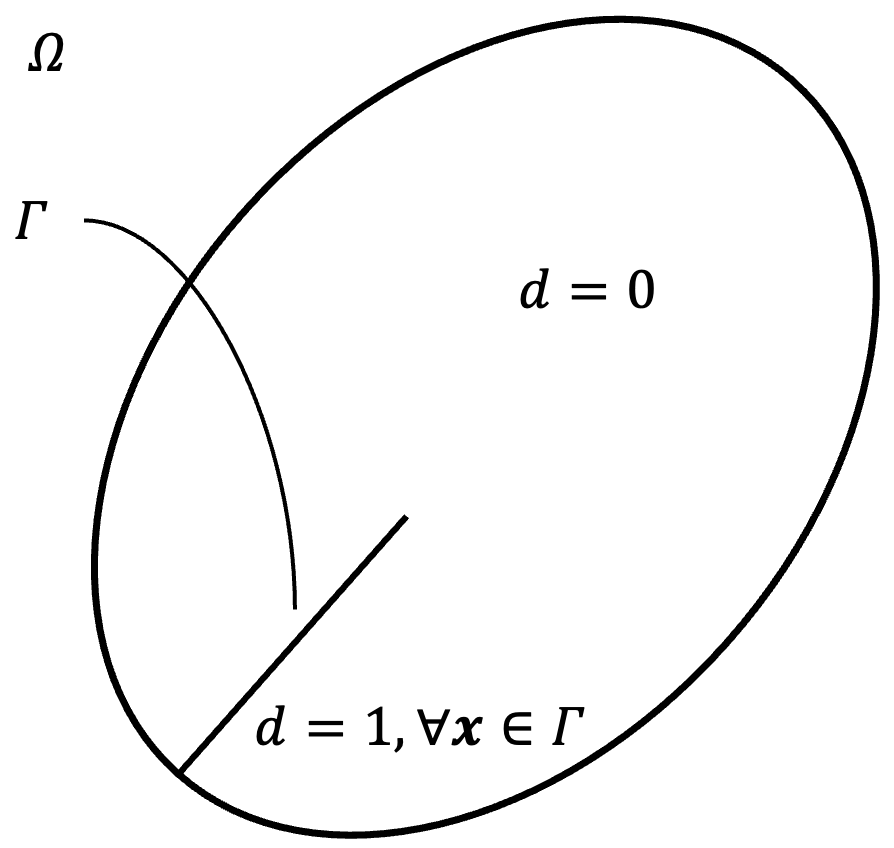}
    \caption{Discrete crack}
    \label{fig:chap1_fig1}
  \end{subfigure}
  \hspace{0.02\textwidth}
  \begin{subfigure}{0.25\textwidth}
    \centering
    \includegraphics[height=3.5cm]{   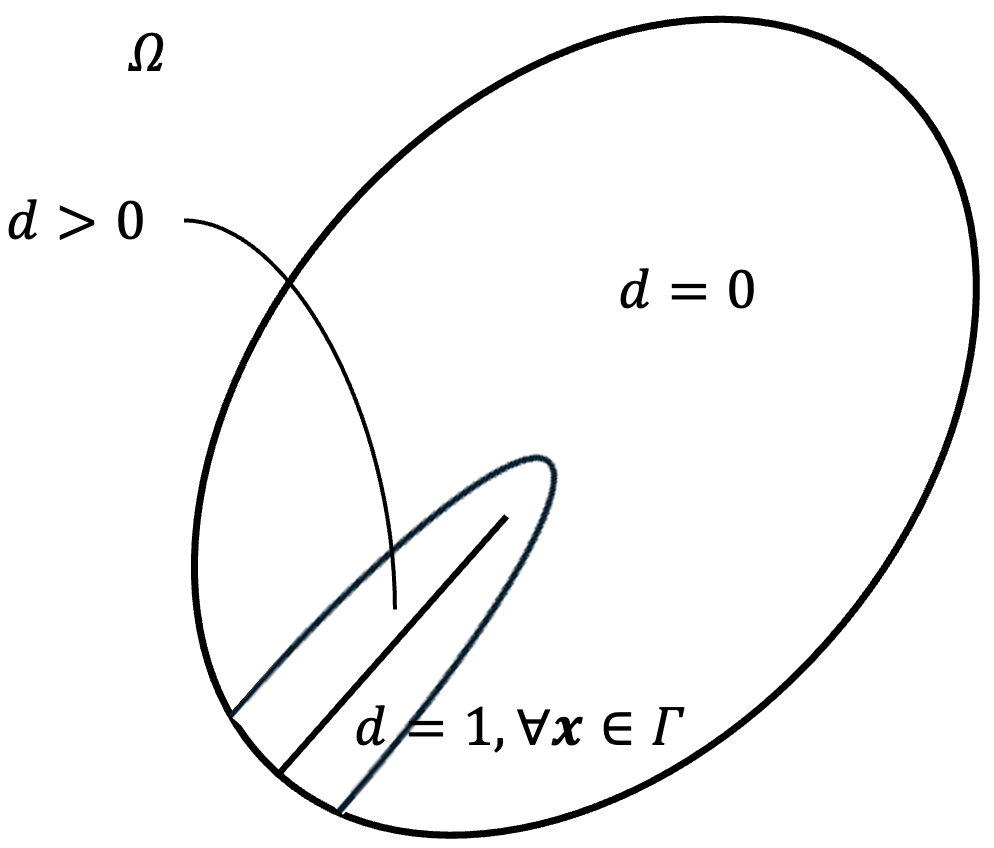}
    \caption{Diffuse crack}
    \label{fig:chap1_fig2}
  \end{subfigure}
  \caption{ {Approximation of crack topology using (a) discrete representation and (b) diffusive representation.}}
  \label{fig:chap1_compare}
\end{figure}

In the PF approach, the crack topology is represented by a diffusive approximation (e.g., \cite{Bourdin2000,Miehe2010}). Let us consider a $\mathit{n}$-dimensional domain $\mathit{\Omega} \subset \mathbb{R}^n$ with a fully open crack $\mathit{\Gamma} \subset \mathbb{R}^{n-1}$. The crack field can then be
 described by a damage function $d(x)$, which is equal to 0 everywhere, except for $\mathit{\Gamma}$. For a fully damaged point on $\mathit{\Gamma}$, we consider that $d=1$. This leads to a sharp crack topology, as described in Fig.~\ref{fig:chap1_compare}\textup{(\subref*{fig:chap1_fig1})}. The diffusive approximation consists in introducing a continuous transition between the crack and the undamaged zone, in which $d \in (0,1)$, representing a partially damaged zone around the crack. A diffusive representation of $\mathit{\Gamma}$ is illustrated in Fig.~\ref{fig:chap1_compare}\textup{(\subref*{fig:chap1_fig2})}. The width of the
transition zone is controlled by a characteristic length $l_c$, and a damage field $d(x)$ of this type can be considered as a solution to the following problem
\begin{equation}
d - l_c^2 \Delta  d=0,
\end{equation}
where $\Delta$ denotes the Laplacian. The crack surface density in $\mathit{\Omega}$ can then be defined as
\begin{equation}
\gamma(d,\nabla d)
= \frac{1}{2l_c}\left(d^2 + l_c^2 \left| \nabla d \right|^2\right),
\end{equation}
and the total crack surface is approximated by its integral over the domain
\begin{equation}
\mathit{\Gamma} \approx \mathit{\Gamma}(d)
= \int_{\mathit{\Omega}} \gamma(d,\nabla d)\,\mathrm{d}V
= \int_{\mathit{\Omega}} \frac{1}{2l_c}\left(d^2 + l_c^2 \left| \nabla d \right|^2\right) \,\mathrm{d}V.
\end{equation}


In this way, the crack topology is described implicitly by a partial differential equation. This can facilitate the tracking of complex crack topologies, such as branched, curved, or multi-crack networks. In our work, we developed the PF approach for DIC analysis.

\subsection{Phase field DIC}
Based on the diffusive approximation of cracks, the proposed PF-DIC can be defined as a variational problem that minimizes the following type of functional 
\begin{equation}
\mathit{\Pi}(\boldsymbol{u},d)= \int_{\mathit{\Omega}}
\mathcal{L}(\boldsymbol{u},d)\ \mathrm{d}V
=
\underbrace{\frac{1}{2}
\int_{\mathit{\Omega}} M(d)\, r^2 (\boldsymbol{u}) \,\mathrm{d}V}_{\text{Grayscale residual}}
+
\underbrace{\frac{1}{2}\int_{\mathit{\Omega}}
\left(
\frac{\mu}{l_c}\, d^2
+
l_c\, \lvert \nabla d \rvert^2
\right)
\,\mathrm{d}V}_{\text{Diffusive cracks}}
-
\underbrace{\frac{1}{2}\int_{\mathit{\Omega}} (1-M(d))\, \mathit{\Phi}(\boldsymbol{u})\,\mathrm{d}V}_{\text{Damage driving force energy}}
,
\label{eq:objective}
\end{equation}
where $M(d)$ is a standard degradation function: $M(d)=(1-d)^2$, $r$ is the grayscale residual defined in Eq.~\eqref{eq:residual}, $\mu$ is the weighting parameter that balances the contributions between the damage field $d$ and its gradient $\nabla d$ to the total crack surface, $\mathit{\Phi}(\boldsymbol{u})$ is an energetic driving force function that should be defined according to specific needs or hypotheses, which will be discussed later. We remark that the weighting parameter $\mu$ does not change the physical meaning of $l_c$ in the diffusive approximation of cracks, i.e., a characteristic length controlling the width of the transition zone.  This $\mu$ can be considered as a constant once an appropriate value is identified for specimens. Throughout this work, we chose $\mu=0.1$, $l_c=5$ .  Based on the above definition, the solutions of PF-DIC can be found by 
\begin{equation}
(\boldsymbol{u}, d)
=
\mathrm{Arg}\left\{
\inf_{\boldsymbol{u}^{\ast},\,\mathit{d}^{\ast}}
\mathit{\Pi}(\boldsymbol{u}^{\ast},\,\mathit{d}^{\ast})
\right\},
\end{equation}
with an appropriate approximation for $\boldsymbol{u}$ and $d$, e.g., FE approximation, and a suitable optimization algorithm. A practical approach is to solve the problem in a staggered way by alternating the solutions for $\boldsymbol{u}$ and $d$.  

Now, considering that $\delta \mathit{\Pi}=0$, we have 
\begin{equation}
\delta \mathit{\Pi} 
=
\int_{\mathit{\Omega}}
\left(
\frac{\partial \mathcal{L}}{\partial \boldsymbol{u}} \delta \boldsymbol{u}
+
\frac{\partial \mathcal{L}}{\partial d}\, \delta d
\right)
\,\mathrm{d}V
= 0,
\end{equation}
which leads to two subproblems for the displacement and damage fields
\begin{equation}
\begin{cases}
\dfrac{\partial \mathcal{L}}{\partial \boldsymbol{u}} = 0,
& \forall\, \delta \boldsymbol{u}, \\[6pt]
\dfrac{\partial \mathcal{L}}{\partial d} = 0,
& \forall\, \delta d .
\end{cases}
\end{equation}
Therefore, the solutions $\boldsymbol{u}$ and $d$ should satisfy the following equations
\begin{itemize}
    \item For displacement
    \begin{equation}
        M(d)\,
r \frac{\partial r}{\partial \boldsymbol{u}} \
=
\frac{1}{2} (1-M(d))\,\frac{\partial \mathit{\Phi}}{\partial \boldsymbol{u}},
\label{eq:disp}
    \end{equation}
    \item For damage 
    \begin{equation}
\frac{\mu}{l_c}\,d
-l_c\,\Delta d  = (1-d)(r^2
+\mathit{\Phi}
),
\label{eq:damage}
    \end{equation}
\end{itemize}
where $\frac{\partial r}{\partial \boldsymbol{u}}$ and $\frac{\partial \mathit{\Phi}}{\partial \boldsymbol{u}}$ can be derived from specific approximations of $r$ and $\mathit{\Phi}$. In the following section, we will discuss some specific choices for $\mathit{\Phi}$.

\subsection{Choices of damage driving force energy} 
In the PF-DIC, the damage driving force energy plays an important role in detecting cracks, as shown in Eq.~\eqref{eq:damage}. Hence, the energetic driving force function $\mathit{\Phi}$ needs to be carefully chosen. In general, we expect the value of this function to be relatively large around the actual crack $\mathit{\Gamma}$, to be consistent with the diffusive approximation. A simple yet practical choice of $\mathit{\Phi}$ is
\begin{equation}
\mathit{\Phi}=\mathcal{A}
\!\left(
\boldsymbol{\varepsilon}, r,f
\right)=w\tanh (\frac{f\,r^2\,
\boldsymbol{\varepsilon}:\boldsymbol{\varepsilon}-a}{c}),
\label{eq:fracture_energy}
\end{equation}
where  {the strain} $\boldsymbol{\varepsilon}
=
\nabla_{\mathrm{s}} \boldsymbol{u}$,  {$w, a, c$    are  the parameters to be calibrated. Specifically, $w$ is a magnitude coefficient, $a$ is a threshold for the local damage energy, and $c$ is a parameter controlling the window size of tanh}. In this case, the function  $\mathit{\Phi}$  is designed to identify regions where the displacement gradient $\nabla \boldsymbol{u}$, the grayscale residual $r$, and the grayscale of image $f$ are relatively large.   This is based on some common observations listed below.
\begin{itemize}
    \item The gradient of displacement, i.e., strain $\boldsymbol{\varepsilon}$, tends to be large across the crack tip, due to the displacement discontinuity. 
    \item The grayscale residual, i.e., $r=g-f$, tends to be large around cracks. This is due to the loss of materials and the potential measurement uncertainty of the displacement field around cracks.
    \item The grayscale value of image, i.e., $f$ or $g$, tends to be large or close to the background on cracks. Here, we assume that the grayscale of the background is uniformly set to 1.
\end{itemize}
Therefore, Eq.~\eqref{eq:fracture_energy} is expected to capture cracks on deformed samples, as the product of the three quantifies would distinguish the cracks from undamaged regions.  {For example, a region with a large strain $\boldsymbol{\varepsilon}$ but with a small grayscale residual $r$ would not be considered as a crack,  as the value of this function would not be high. Eq.~\eqref{eq:fracture_energy} can be considered as a damage indicator to some extent. The influence of the parameters $w,a,c$ will be discussed later.}

Other choices of $\mathit{\Phi}$ can also be considered in the PF-DIC, by modifying the inputs in Eq.~\eqref{eq:fracture_energy}. For example, we can consider that only part of the total strain energy would contribute to the damage 
\begin{equation}
\mathit{\Phi}=\mathcal{A}
\!\left(
\boldsymbol{\varepsilon}^+, r,f
\right)=w\tanh (\frac{f\,r^2\,
\boldsymbol{\varepsilon}^+:\boldsymbol{\varepsilon}^+-a}{c}),
\label{eq:fracture_energy_select}
\end{equation}
where $\boldsymbol{\varepsilon}^+$ denotes the positive part of the strain, which can be obtained by various decomposition methods in the literature on PF modeling of fracture \cite{amor2009regularized,Miehe2010,zhang2022assessment} or by  a user-defined decomposition. In this case,  {Eq.~\eqref{eq:fracture_energy_select} only highlights the region where the tensile strain energy is significant,  enabling the PF-DIC to identify cracks that are only sensitive to this specific driving force. In a general sense, the modification of $\mathit{\Phi}$  can enable the selective detection of critical cracks under certain fracture modes or mechanisms}. This is particularly important for the characterization of defects in materials and structures. This capability of PF-DIC will be demonstrated in numerical examples.

Another potential choice of $\mathit{\Phi}$ is to include the grayscale value of images in a separate function
\begin{equation}
\mathit{\Phi}=\mathcal{A}_1
\!\left(
\boldsymbol{\varepsilon}, r,f
\right)+\mathcal{A}_2
\!\left(f
\right)=w_1\tanh (\frac{f\,r^2\,
\boldsymbol{\varepsilon}:\boldsymbol{\varepsilon}-a_1}{c_1})+w_2\tanh (\frac{f\,-a_2}{c_2}),
\label{eq:fracture_energy_w_background}
\end{equation}
where  $w_{1,2}$, $a_{1,2}$ and $c_{1,2}$ are parameters to be calibrated. In this case, the first term $\mathcal{A}_1$ plays a similar role as previously, but with the help of the second term $\mathcal{A}_2$ to detect the background pixels in the images. This is advantageous when cracks are wide open, as the grayscale value can be used as a direct indicator for damaged regions in such scenarios. Additionally, with appropriate parameters, $\mathcal{A}_2$  can help automatically create a damage field on the background of images, excluding  unnecessary information for the displacement analysis without additional trimmings for the images.  

The flexibility of choosing  $\mathit{\Phi}$ should be considered as an advantage of PF-DIC. It can be defined according to user-specifications, and the parameters of functions can be calibrated based on some baseline experimental samples. In this work, we will demonstrate the influence of this function using the several choices mentioned earlier. The specific parameters will be provided later, along with numerical results.

\subsection{A special case of strong form} 
The general strong form problem of PF-DIC is given by Eqs.~\eqref{eq:disp} and \eqref{eq:damage}, which is intrinsically a nonlinear two-way coupled problem. To simplify the numerical implementation, we can adopt the following approximations. 

First, let us consider an incremental solution procedure for $\boldsymbol{u}$ with $\boldsymbol{u}_{i+1}=\boldsymbol{u}_i+\delta \boldsymbol{u}$, where $\boldsymbol{u}_i$ is assumed known and $\delta \boldsymbol{u}$ is the unknown incremental displacement at the $i$-th iteration. The grayscale residual can then be written as
\begin{equation}
    r(\boldsymbol{u}_{i+1}) = g(\boldsymbol{x}) - f(\boldsymbol{x} - \boldsymbol{u}_i-\delta \boldsymbol{u})\approx r(\boldsymbol{u}_{i}) +\nabla f (\boldsymbol{x} - \boldsymbol{u}_i)\cdot \delta \boldsymbol{u}.
\label{eq:residual_linearization}
\end{equation}
where $r(\boldsymbol{u}_{i}) =g(\boldsymbol{x}) - f(\boldsymbol{x} - \boldsymbol{u}_i)$, and the first-order Taylor expansion of $f$ was used with 
\begin{equation}
f(\boldsymbol{x} - \boldsymbol{u}_i-\delta \boldsymbol{u})
\approx
f(\boldsymbol{x} - \boldsymbol{u}_i)-\nabla f (\boldsymbol{x} - \boldsymbol{u}_i)\cdot \delta \boldsymbol{u}.
\end{equation}
Therefore, we have
\begin{equation}
    \frac{\partial r}{\partial \boldsymbol{u}_{i+1}} =\nabla f (\boldsymbol{x} - \boldsymbol{u}_i).
    \label{eq:approx_partial_r}
\end{equation}
In addition, we can notice that given the definition of $M(d)$, the rhs of Eq.\eqref{eq:disp}  diminishes as $d\rightarrow0$. This leads to 
\begin{equation}
    \frac{1}{2} (1-M(d))\,\frac{\partial \mathit{\Phi}}{\partial \boldsymbol{u}}\approx0, \quad \forall\boldsymbol{x}\in \mathit{\Omega}/\mathit{\Gamma}.
    \label{eq:approx_rhs}
\end{equation}

Now, plugging the above approximations into Eq.\eqref{eq:disp} yields the following strong form problem for the displacement field
\begin{equation}
        M(d)\,
(g - f(\boldsymbol{x} - \boldsymbol{u}_i)+\nabla f (\boldsymbol{x} - \boldsymbol{u}_i)\cdot \delta \boldsymbol{u}) \nabla f (\boldsymbol{x} - \boldsymbol{u}_i)\
=
 0,
    \end{equation}
which is equivalent to a Gauss-Newton algorithm for solving $\boldsymbol{u}$. In practice, we might consider a modified version of the method with $\nabla f (\boldsymbol{x} - \boldsymbol{u}_i)\approx\nabla f (\boldsymbol{x} )$. The resulting strong form problem then becomes
\begin{itemize}
    \item For displacement
    \begin{equation}
        M(d)\,
(g - f(\boldsymbol{x} - \boldsymbol{u}_i)+\nabla f \cdot \delta \boldsymbol{u}) \nabla f \
=
 0,
\label{eq:disp_mod}
    \end{equation}
    \item For damage 
    \begin{equation}
\frac{\mu}{l_c}\,d
-l_c\,\Delta d  = (1-d)(r^2(\boldsymbol{u}_{i+1})
+\mathit{\Phi}(\boldsymbol{u}_{i+1})
).
\label{eq:damage_mod}
    \end{equation}
\end{itemize}

We remark that other approximations can also be used for the displacement equation in the PF-DIC. Here, the modified Gauss-Newton method is considered for illustration purposes.

\subsection{Finite element formulation}
The PF-DIC problem can be discretized with the FE method by considering that
\begin{equation}
\delta\boldsymbol{\mathit{u}}({\boldsymbol{x}}) = \boldsymbol{\mathit{N}}_{\mathit{u}}(\boldsymbol{{x}})\,
\boldsymbol{\mathit{U}},
\qquad
\mathit{d}(\boldsymbol{{x}})
=
\boldsymbol{{N}}_{\mathit{d}}(\boldsymbol{{x}})\,
\boldsymbol{\mathit{d}},
\label{eq:FE_interpolation}
\end{equation}
where $\boldsymbol{\mathit{N}}_{\mathit{u}}$ and $\boldsymbol{\mathit{N}}_{\mathit{d}}$ are the FE shape functions, $\boldsymbol{\mathit{U}}$ and $\boldsymbol{\mathit{d}}$ are the nodal displacement increment and damage fields. The FE discretized equations can generally be written as 
\begin{itemize}
    \item For displacement
    \begin{equation}
       \frac{\partial \mathit{\Pi}(\boldsymbol{u},d)}{\partial \boldsymbol{\mathit{u}}}=
 0 \quad\rightarrow \quad \boldsymbol{R}_u(\boldsymbol{{U}},\boldsymbol{\mathit{d}})=0,
 \label{eq:disp_discrete}
    \end{equation}
    \item For damage 
    \begin{equation}
\frac{\partial \mathit{\Pi}(\boldsymbol{u},d)}{\partial d}=
 0 \quad\rightarrow \quad \boldsymbol{R}_d(\boldsymbol{{U}},\boldsymbol{\mathit{d}})=0,
 \label{eq:damage_discrete}
    \end{equation}
\end{itemize}
where $\boldsymbol{R}_u$ and $\boldsymbol{R}_d$ are the residuals to be minimized during the solution procedure. Taking  Eqs.~\eqref{eq:disp_mod} and \eqref{eq:damage_mod} as an example, we derive a specific form of the discretized equations in the following. We start from the weak formulation corresponding to Eqs.~\eqref{eq:disp_mod} and \eqref{eq:damage_mod}, i.e., 
    \begin{equation}
        \int_{\mathit{\Omega}} (\delta \boldsymbol{u}^{*})^T M(d)\,
(g - f(\boldsymbol{x} - \boldsymbol{u}_i)+\nabla f \cdot \delta \boldsymbol{u}) \nabla f \ \mathrm{d}V 
=
 0,
    \end{equation}
    and
    \begin{equation}
\int_{\mathit{\Omega}} (d^{*} \frac{\mu}{l_c}\,d
+l_c (\nabla d^{*})^T \,\nabla d ) \ \mathrm{d}V = \int_{\mathit{\Omega}} d^{*} (1-d)(r^2(\boldsymbol{u}_{i+1})
+\mathit{\Phi}(\boldsymbol{u}_{i+1}) 
)\ \mathrm{d}V .
    \end{equation}
    where $\delta \boldsymbol{u}^{*}$ and $d^{*}$ are the test functions. By rearrangement and considering that $\delta\boldsymbol{\mathit{u}}^*= \boldsymbol{\mathit{N}}_{\mathit{u}}\,
\boldsymbol{\mathit{U}}^*
$,  $
\mathit{d}^*
=
\boldsymbol{{N}}_{\mathit{d}}\,
\boldsymbol{\mathit{d}}^*,$ and $\boldsymbol{B}=\partial \boldsymbol{N}/\partial \boldsymbol{x}$, we can obtain the following discretized form of the problem
    \begin{equation}
        \boldsymbol{U}^{*T} \boldsymbol{Q}_u+\boldsymbol{U}^{*T} \boldsymbol{K}_u \, \boldsymbol{U}
=
 0,
    \end{equation}
    and
    \begin{equation}
\boldsymbol{d}^{*T} \boldsymbol{M}_d \, \boldsymbol{d}
+\boldsymbol{d}^{*T} \boldsymbol{K}_d\, \boldsymbol{d}  = \boldsymbol{d}^{*T} \boldsymbol{Q}_d ,
    \end{equation}
    with 
    \begin{equation}
\begin{cases}
\displaystyle
         \boldsymbol{Q}_u=\int_{\mathit{\Omega}}  \boldsymbol{N}_u^T (1-d)^2\,
(g - f(\boldsymbol{x} - \boldsymbol{u}_i))\ \mathrm{d}V,\\[2ex]
\boldsymbol{K}_u=\int_{\mathit{\Omega}}  \boldsymbol{N}_u^T (1-d)^2\, \nabla f\,\nabla^T f \, \boldsymbol{N}_u  \ \mathrm{d}V,\\[2ex]
\boldsymbol{M}_d=\int_{\mathit{\Omega}} \boldsymbol{N}_d^{T} (\frac{\mu}{l_c}+r^2(\boldsymbol{u}_{i+1})
+\mathit{\Phi}(\boldsymbol{u}_{i+1}) )\,\boldsymbol{N}_d \ \mathrm{d}V,\\[2ex]
\boldsymbol{K}_d=\int_{\mathit{\Omega}} l_c \boldsymbol{B}_d^T \boldsymbol{B}_d  \ \mathrm{d}V,\\[2ex]
\boldsymbol{Q}_d=\int_{\mathit{\Omega}} \boldsymbol{N}_d^{T} (r^2(\boldsymbol{u}_{i+1})
+\mathit{\Phi}(\boldsymbol{u}_{i+1}) 
)\ \mathrm{d}V. 
\label{eq:discrete_operator}
\end{cases}
    \end{equation}
Therefore, the final discrete problem using a modified Gauss-Newton method can be written as 
\begin{itemize}
    \item For displacement
\begin{equation}
        \boldsymbol{R}_u(\boldsymbol{{U}},\boldsymbol{\mathit{d}})= \boldsymbol{Q}_u+ \boldsymbol{K}_u \, \boldsymbol{U}=0,
        \label{eq:disp_discrete_gauss}
    \end{equation}
    \item For damage 
    \begin{equation}
 \boldsymbol{R}_d(\boldsymbol{{U}},\boldsymbol{\mathit{d}})=(\boldsymbol{M}_d 
+ \boldsymbol{K}_d)\, \boldsymbol{d}  -  \boldsymbol{Q}_d=0.
\label{eq:damage_discrete_gauss}
    \end{equation}
\end{itemize}

In this derivation, we assumed that $\nabla f$ has been obtained from an approximation of $f(x)$. This can be done in many different ways, as discussed in \cite{Lu2024}. Solving the above equations will give the solution $\delta\boldsymbol{\mathit{u}}$ and $d$ at the $i$-th iteration. The final solution is obtained by updating $\boldsymbol{u}_{i+1}=\boldsymbol{u}_i+\delta \boldsymbol{u}$ until convergence. The detailed solution procedure and numerical implementation will be described in the next.

\subsection{Solution procedure and implementation details}
We present here a staggered algorithm to solve the discrete equations for the displacement and damage solutions, although a monolithic approach can also be considered. The general procedure is summarized below. 
\begin{itemize}
  \item \textbf{Initialization}: $\boldsymbol{u}_0=0$, $d=0$, $\forall \boldsymbol{x}\in \mathit{\Omega}$, and $r(\boldsymbol{u}_0)=g-f$.
  \item \textbf{Loop} $i=0,1,2, \dots$

  \item \textbf{Displacement update}:
  \begin{enumerate}
    \item Compute a displacement increment $\delta \boldsymbol{u}$ by solving $\boldsymbol{R}_u= 0$ with  $d$ fixed.
  \item Update the displacement field: $\boldsymbol{u}_{i+1}\leftarrow\boldsymbol{u}_i+\delta \boldsymbol{u}$.   
  \end{enumerate}
    \item \textbf{Damage update}:
    \begin{enumerate}
        \item Update $r(\boldsymbol{u}_{i+1})$ and $\mathit{\Phi}(\boldsymbol{u}_{i+1})$.
    \item Update the damage field $d$ by solving $\boldsymbol{R}_d= 0$  with  $\boldsymbol{u}_{i+1}$ fixed.
    \end{enumerate}
    \item \textbf{Convergence check}:
     If 
    $\lVert M(d)\, \delta \boldsymbol{u}\rVert_2\leq \epsilon \,\lVert M(d)\, \boldsymbol{u}_1\rVert_2 $,
    where $\epsilon$ is a small threshold, then, END; Otherwise,
     $(\cdot)_i \leftarrow (\cdot)_{i+1}$ and continue.     
\end{itemize}

\begin{figure}[!htbp]
  \centering
  \includegraphics[height=3.5cm]{   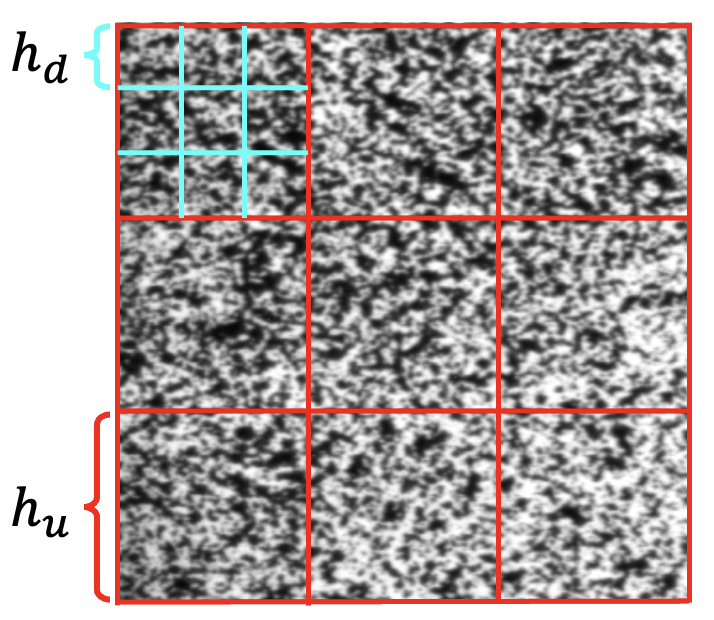}
  \caption{Meshes on a speckle pattern. $h_u$: element (subset) size for displacement, $h_d$: element size for damage.}
  \label{fig:hd_hu_mesh}
\end{figure}

In the above procedure, the degradation function $M(d)$ is involved in the convergence criterion to exclude the displacement on pixels that represent cracks. With an appropriate $\epsilon$, this iterative procedure can converge in a few (two or three) iterations for a good accuracy, as shown in the examples. Furthermore, this procedure can be applied to different optimization schemes and can generally work with both global and local versions of DIC. In our implementation, we used a FE approach, and considered that the displacement and damage fields, i.e., $\boldsymbol{u}$ and $d$, are discretized with two meshes of different resolutions. Fig.~\ref{fig:hd_hu_mesh} illustrates an example of the two meshes in the ROI, with $h_u$ denoting the element size for displacement and $h_d$ denoting the one for damage. Particularly, we chose a fine mesh with $h_d=1$ (pixel)  for the damage field, as it provides a precise description for cracks, whereas a coarser mesh was used for the displacement, as it provides relatively smoother displacement solutions. As a result, a mapping between the two meshes is needed, particularly when updating the damage driving force energy $\mathit{\Phi}(\boldsymbol{u}_{i+1})$. Taking  Eq.~\eqref{eq:fracture_energy} as an example, we can see that updating $\mathit{\Phi}(\boldsymbol{u}_{i+1})=\mathcal{A}
\!\left(
\boldsymbol{\varepsilon}, r,f
\right)$ requires the values of $\boldsymbol{\varepsilon}(\boldsymbol{u}_{i+1})$ and $r(\boldsymbol{u}_{i+1})$ on a fine mesh of $d$, while $\boldsymbol{u}_{i+1}$ is represented on a coarse mesh. A straightforward approach to handle this mapping is to use the FE shape function and the nodal solution of $\boldsymbol{u}_{i+1}$ to evaluate the strain $\boldsymbol{\varepsilon}$ and the grayscale residual $r$ directly at the center of each pixel, so as to link the coarse mesh solution to the fine mesh. This can be done in a natural way when the integration in Eq.~\eqref{eq:discrete_operator} uses a  scheme in which the quadrature points coincide with the centers of pixels.  Additionally, we found that if an averaging operator is further applied to $r$ and $f$ for calculating $\mathcal{A}
\!\left(
\boldsymbol{\varepsilon}, r,f
\right)$, a better damage solution can be obtained. Denoting this averaging operator by $\langle\cdot\rangle=\frac{1}{N} \sum_{i=1}^N (\cdot)$ with $N$ as the number of pixels in a predefined subregion of ROI,  we then have  
\begin{equation}
\mathit{\Phi}(\boldsymbol{u}_{i+1})=\mathcal{A}
\!\left(
\boldsymbol{\varepsilon}, \langle r \rangle,\langle f \rangle
\right).   
\end{equation}
Here we assumed the ROI has been split into some non-overlapping  subregions with $N$ pixels each for the averaging operator, and they can be defined separately for $\langle r \rangle$ and $\langle f \rangle$. In our work, the specific values of $N$ for the averaging operator are 16 (based on a $4\times4$ subregion) and 64 (based on a $8\times8$ subregion) for $\langle r \rangle$ and $\langle f \rangle$, respectively. This stays the same for all the numerical examples presented in the next section.

\FloatBarrier

\section{Numerical examples}
This section presents the numerical results obtained by the proposed PF-DIC. In particular, we focus on the accuracy of the displacement and damage fields. The strain can be calculated in a standard way with the given displacement field.  We used both synthetic and experimental images of fractured samples to demonstrate the capability of the PF-DIC. 

\subsection{Generation of synthetic images}  
To generate relevant speckle images that mimic fractured samples, we first performed FE simulations of  specimens with some prescribed cracks to determine the reference displacement fields under different loading conditions. These simulations include standard single edge notched tensile tests and three-point bending experiments. The prescribed cracks are represented by a diffusive approximation using the PF approach, as explained in Section \ref{sec:diff_crack}. The modeling and implementation details for  PF modeling of fracture can be found in the relevant literature (e.g., \cite{Miehe2010,Lu2020}). The only difference is that the damage field was assumed already known and we only needed to compute the corresponding displacement.

Once the simulated displacement field is known, we can then use it to deform a reference speckle image $f$, thereby generating its deformed counterpart $g$. The reference speckle image can be obtained in many different ways. In our work, we used a reference image from the DIC challenge 2.0 \cite{reu2022dic}, given the high quality speckle patterns that mimic realistic experimental ones, which has been widely adopted for testing DIC algorithms. Furthermore, to represent the cracks in the speckle image, we set $f=1$ for the pixels where $d\geq0.5$. The deformation of the "cracked" image $f$ was done by inputting a pixel-wise displacement field to a backward mapping, using the OpenCV library in Python. 
In this way, we generated the synthetic reference and deformed images for the examples in Section \ref{sec:single_edge}-\ref{sec:bending}.



\subsection{Single edge notched tensile test}
\label{sec:single_edge}
\begin{figure}[t]
\centering
\hspace*{-25mm}
\begin{subfigure}[b]{0.2\textwidth}
  \centering
  \includegraphics[height=5cm,keepaspectratio]{ 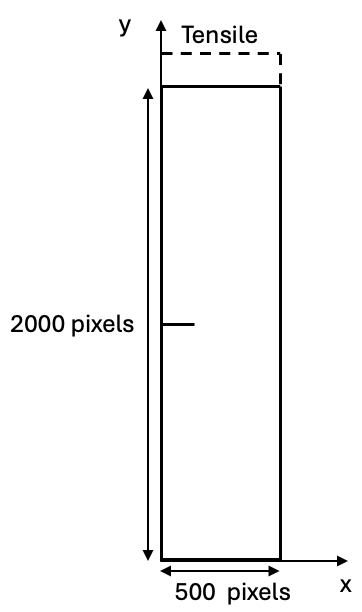}
  \caption{Geometry}
  \label{fig:chap3_geometry}
\end{subfigure}
\hspace{0.05\textwidth}
\begin{subfigure}[b]{0.2\textwidth}
  \centering
  \includegraphics[height=4.5cm,keepaspectratio]{ 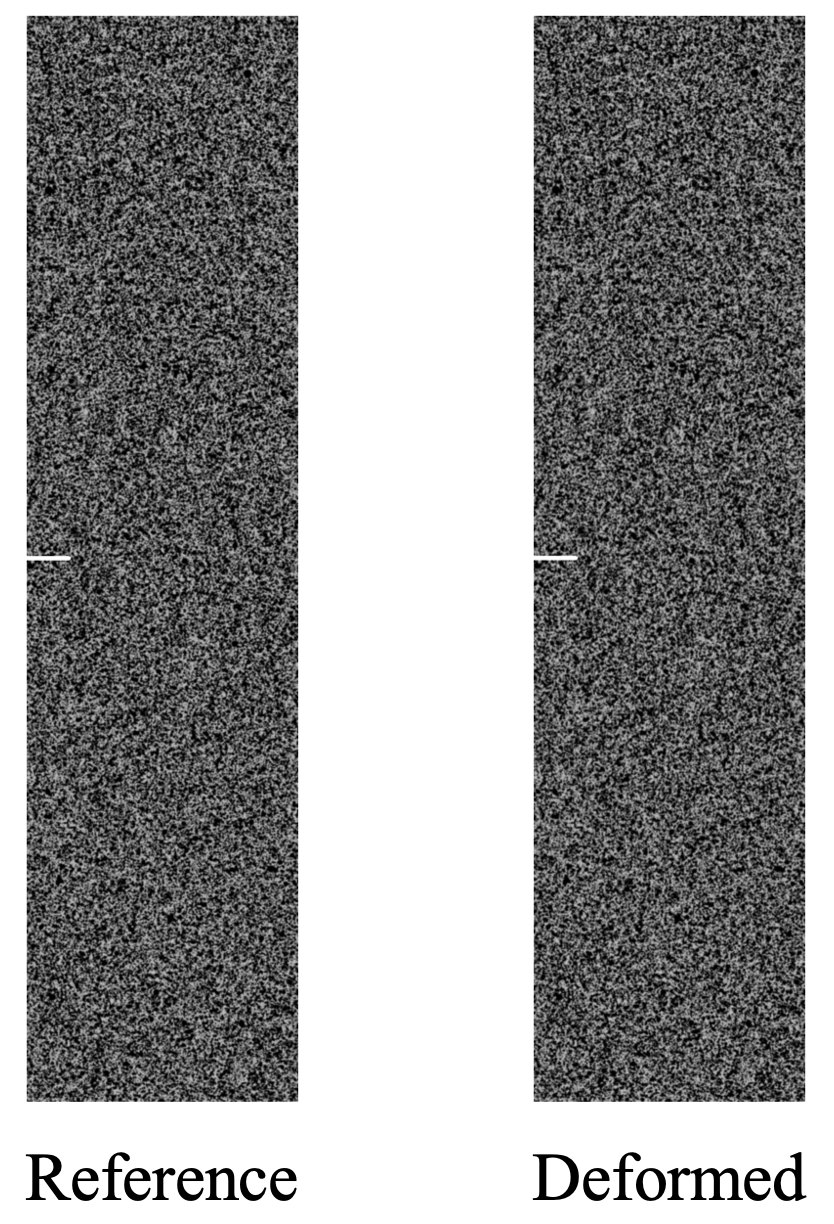}
  \caption{Speckle images}
  \label{fig:chap3_3row_ref}
\end{subfigure}
\hspace{0.07\textwidth}
\begin{subfigure}[b]{0.2\textwidth}
  \centering
  \includegraphics[height=4.5cm,keepaspectratio]{ 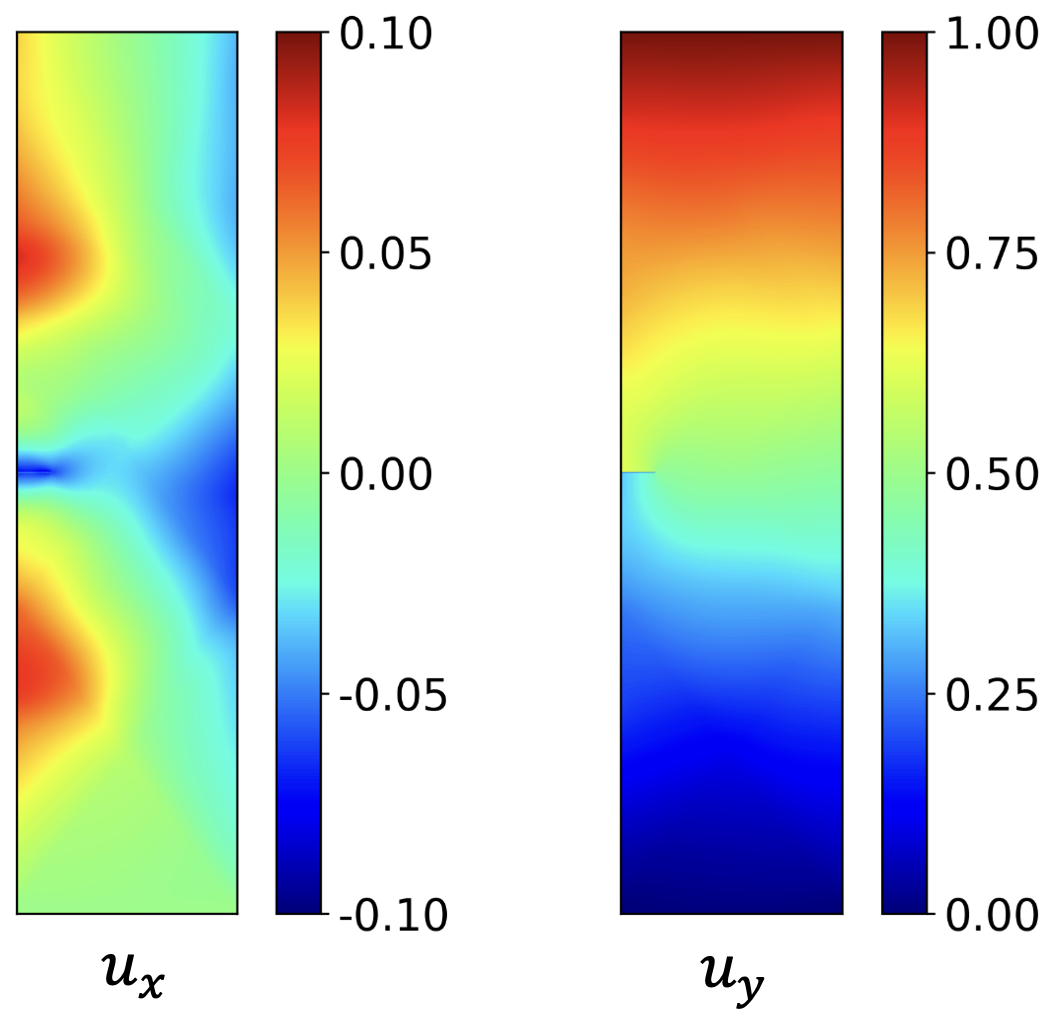}
  \caption{Displacement}
  \label{fig:chap3_3row_ux}
\end{subfigure}
\caption{ {(a) Geometry, (b) reference and deformed images, and (c) ground-truth displacement for the specimen with a small straight crack}}
\label{fig:chap3_3row_all}
\end{figure}

We first tested the proposed PF-DIC using a single edge notched specimen. The problem setup is illustrated in Fig.~\ref{fig:chap3_3row_all}. The displacement is prescribed on the top boundary of the specimen, while the bottom boundary is fixed. The speckle images corresponding to the reference and deformed stages are illustrated in Fig.~\ref{fig:chap3_3row_all}\textup{(\subref{fig:chap3_3row_ref})}. The resulting displacement is illustrated in Fig.~\ref{fig:chap3_3row_all}\textup{(\subref{fig:chap3_3row_ux})}. This displacement is considered as the reference solution for the PF-DIC measurements for displacement. In addition, it is expected that the PF-DIC can automatically detect the small crack in the speckle images, although it seems difficult to visually identify the small crack from the displacement field.

Regarding the PF-DIC, we used an FE mesh of 4-node elements with the element size $h_u = 16$ (pixel) for the displacement field and $h_d=1$ (pixel) for the damage field. This choice is appropriate for the given speckle size, as a smaller $h_u$ would lead to a noisy displacement field. This limit on the spatial resolution is common for standard DIC. For the damage driving force energy, we adopted the definition in Eq.~\eqref{eq:fracture_energy_w_background}. The associated hyperparameters $w_{1,2}$, $a_{1,2}$, and $c_{1,2}$ are given in Table~\ref{tab:parameters}. In general, these parameters should be calibrated against specific specimens. In this example, different choices of  the damage driving force energy, e.g., Eq.~\eqref{eq:fracture_energy}, can also be used, and the calibrated hyperparameters could be different.

\begin{table}[H]
\centering
\setlength{\abovecaptionskip}{2pt}
\caption{Hyperparameters of $\mathit{\Phi}$  in tensile tests.}
\label{tab:parameters}
\small
\begin{tabular}{ccccccc}
\toprule
Parameter & $w_1$ & $a_1$ & $c_1$ & $w_2$ & $a_2$ & $c_2$ \\
\midrule
Value & 0.5 & 0.003 & $1\times10^{-5}$ & 0.1 & 0.6 & 0.05 \\
\bottomrule
\end{tabular}
\end{table}

\begin{figure}[!htbp]
  \centering
  \begin{subfigure}[t]{0.42\textwidth}
    \centering
    \includegraphics[width=\linewidth]{ 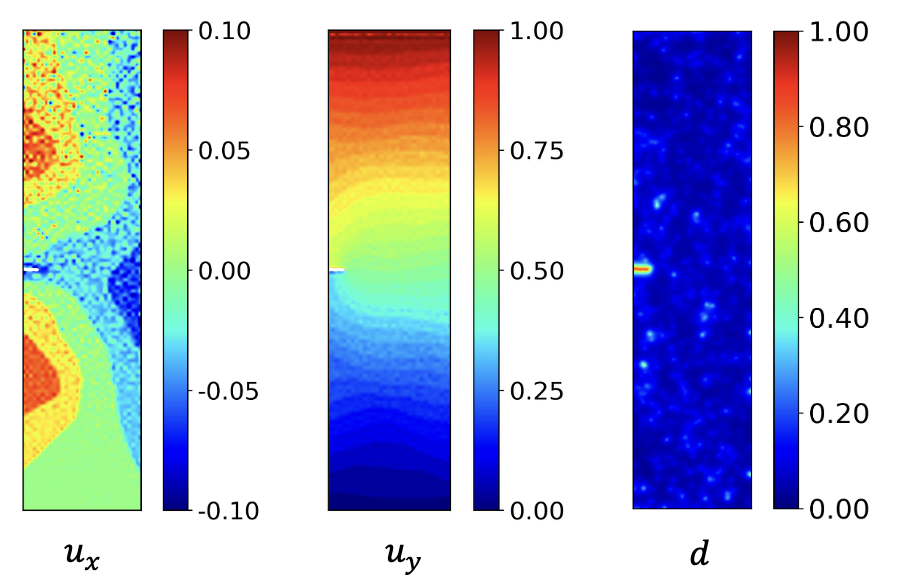}
    \caption{PF-DIC results at iteration 1}
    \label{fig:three_row_a}
  \end{subfigure}
  \hspace{0.08\textwidth}
  \begin{subfigure}[t]{0.42\textwidth}
    \centering
    \includegraphics[width=\linewidth]{ 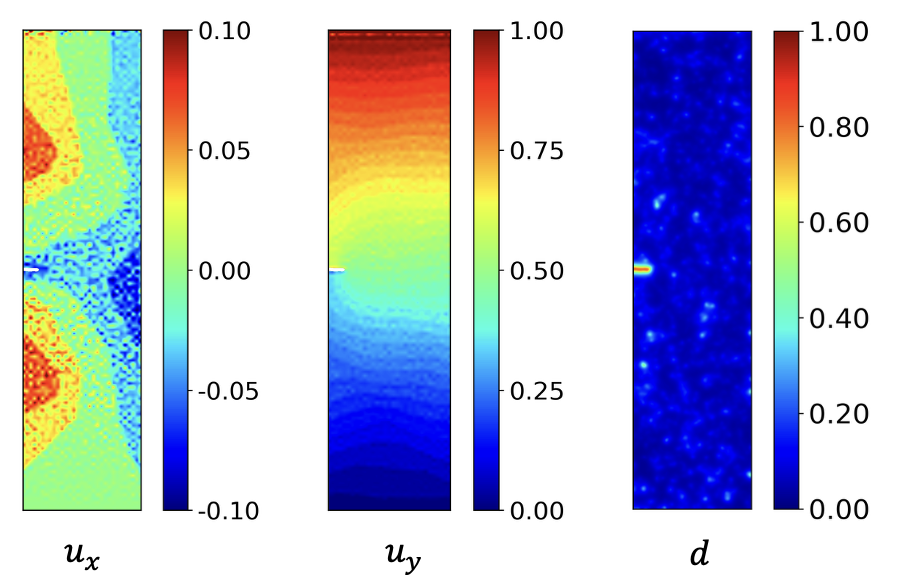}
    \caption{PF-DIC converged results (iteration 3)}
    \label{fig:three_row_b}
  \end{subfigure}\\[0.25cm]
  
  \begin{subfigure}[t]{0.35\textwidth}
    \centering
    \includegraphics[width=\linewidth]{ 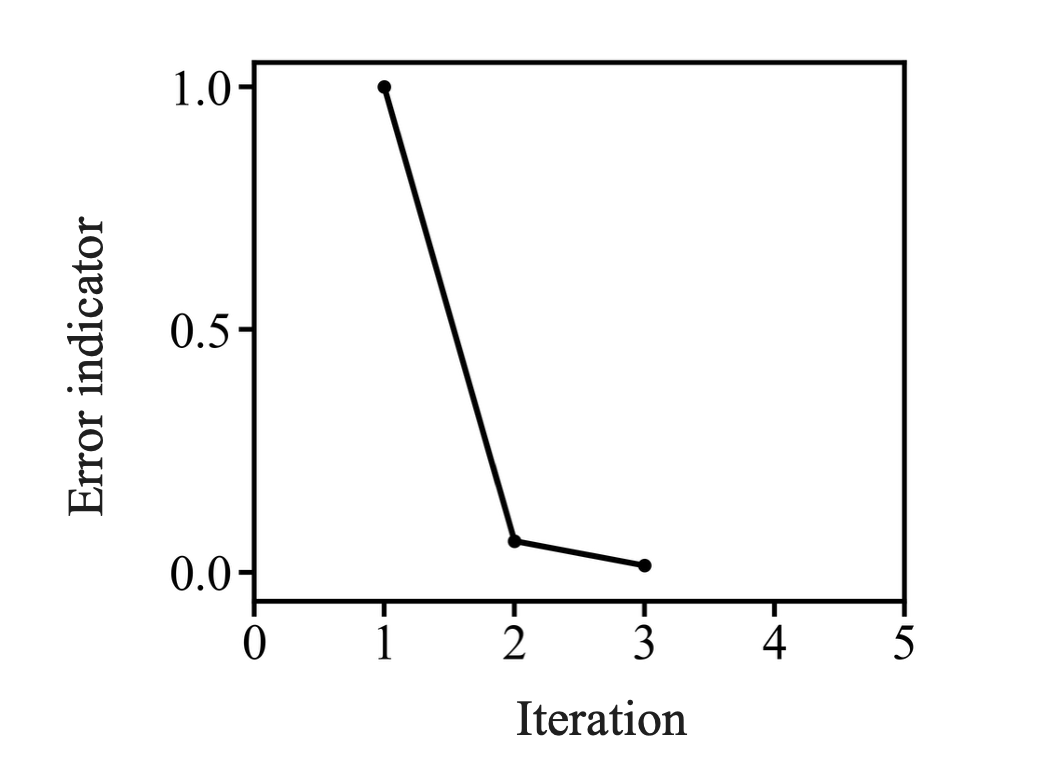}
    \caption{ {Convergence curve for the criterion $\frac{\lVert M(d)\, \delta \boldsymbol{u}\rVert_2}{\lVert M(d)\, \boldsymbol{u}_1\rVert_2 }\leq  0.05$}}
    \label{fig:three_row_c}
  \end{subfigure}

  \caption{ {Displacement and damage measurements from PF-DIC for the specimen with a straight crack. (a) PF-DIC results at the first iteration; (b) PF-DIC converged results at the third iteration; (c) History of the error indicator  $\frac{\lVert M(d)\, \delta \boldsymbol{u}\rVert_2}{\lVert M(d)\, \boldsymbol{u}_1\rVert_2 }$.}}
  \label{fig:chap3_iter_results}
\end{figure}

The results of the PF-DIC are illustrated in Fig.~\ref{fig:chap3_iter_results}. Different from conventional DIC, the PF-DIC provides a full-field measurement of damage in addition to displacement. For illustration purposes, the spurious displacement on the crack was removed from all displacement measurements.  {This can be easily done by removing the pixels where the damage value is high , e.g., $d>0.7$. Fig.~\ref{fig:chap3_iter_results} illustrates the iterative results of the PF-DIC for a convergence criterion $\epsilon=0.05$. A smaller value of this criterion could lead to more iterations, but the current  results converged at iteration 3 seem well acceptable.}  Both displacement and damage results converged quickly, which confirms the feasibility of the staggered solution framework  {and the convergence criterion}.  {In terms of computational cost, the PF-DIC  incurs only a small additional overhead compared with standard DIC. In our implementation, the PF-DIC took 8.8 s for this example, whereas standard DIC took 6.1 s.} We can see that the overall smoothness of displacement measurements improved with the iterative scheme, and the identified displacement agrees well with the reference solution in Fig.~\ref{fig:chap3_3row_all}\textup{(\subref{fig:chap3_3row_ux})}. Regarding the damage measurement, it seems that the damage field converged at the first iteration in this example. As expected, the crack was clearly identified by the  PF-DIC, which confirms the appropriateness of the choice of the damage driving force energy $\mathit{\Phi}$.

\begin{figure}[!htbp]
  \centering
  \includegraphics[height=6.5cm]{ 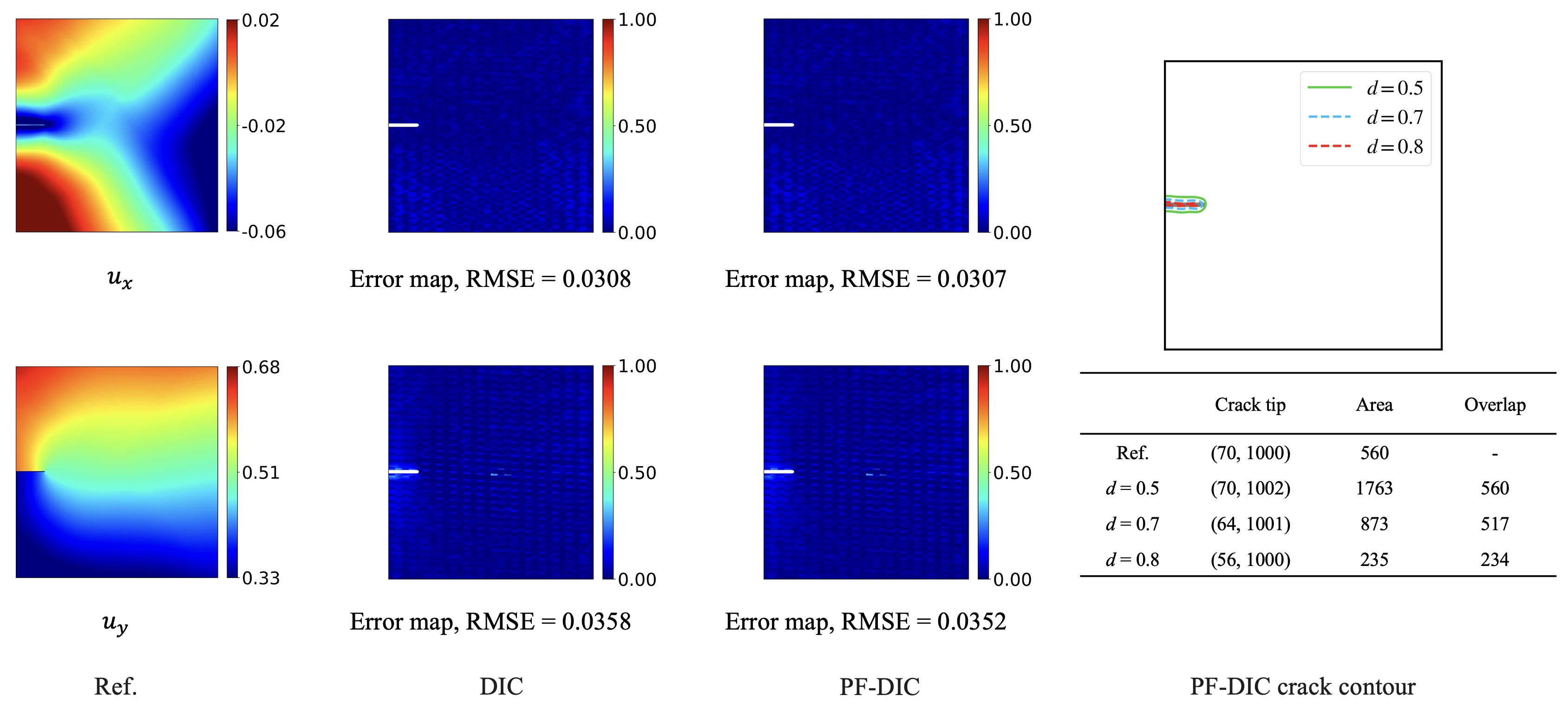}
  \caption{ {Error of the displacement fields around the straight crack and verification of the crack contour identified by the PF-DIC. All units are in pixels. }}
  \label{fig:three_row_reference}
\end{figure}
\FloatBarrier

 {For a quantitative comparison, we calculated the pixel-wise displacement error based on $\frac{|u-u_{\text{ref.}}|}{|u_{\text{ref.},\max}|}$, where $u$ denotes the final  displacement  identified by the DIC or PF-DIC, and $u_{\text{ref.}}$ denotes the ground-truth reference. Particularly, we focused on the displacement around the crack, and removed all the pixels representing the actual crack when calculating the pixel-wise error. These results  are illustrated  in   Fig.~\ref{fig:three_row_reference}. We can see that the DIC and PF-DIC both provided consistent and relatively accurate displacement around the crack. The root mean square error (RMSE) confirmed the overall agreement in displacement measurements. Regarding the damage measurement, we compared the crack contour identified by the PF-DIC with the actual crack. As also shown in Fig.~\ref{fig:three_row_reference}, the PF-DIC measured crack shape matches well the actual one (represented by the gray region) and the crack tip location is well identified, especially when the contour $d=0.5$ is considered. Given the diffusive nature of the representation, the crack contour is wider than the actual crack. This transition region should be controllable by the characteristic length $l_c$. A sensitivity study will be conducted later for this parameter.}


\begin{figure}[!htbp]
  \centering

\begin{subfigure}[b]{0.2\textwidth}
  \centering
  \includegraphics[height=4.5cm,keepaspectratio]{   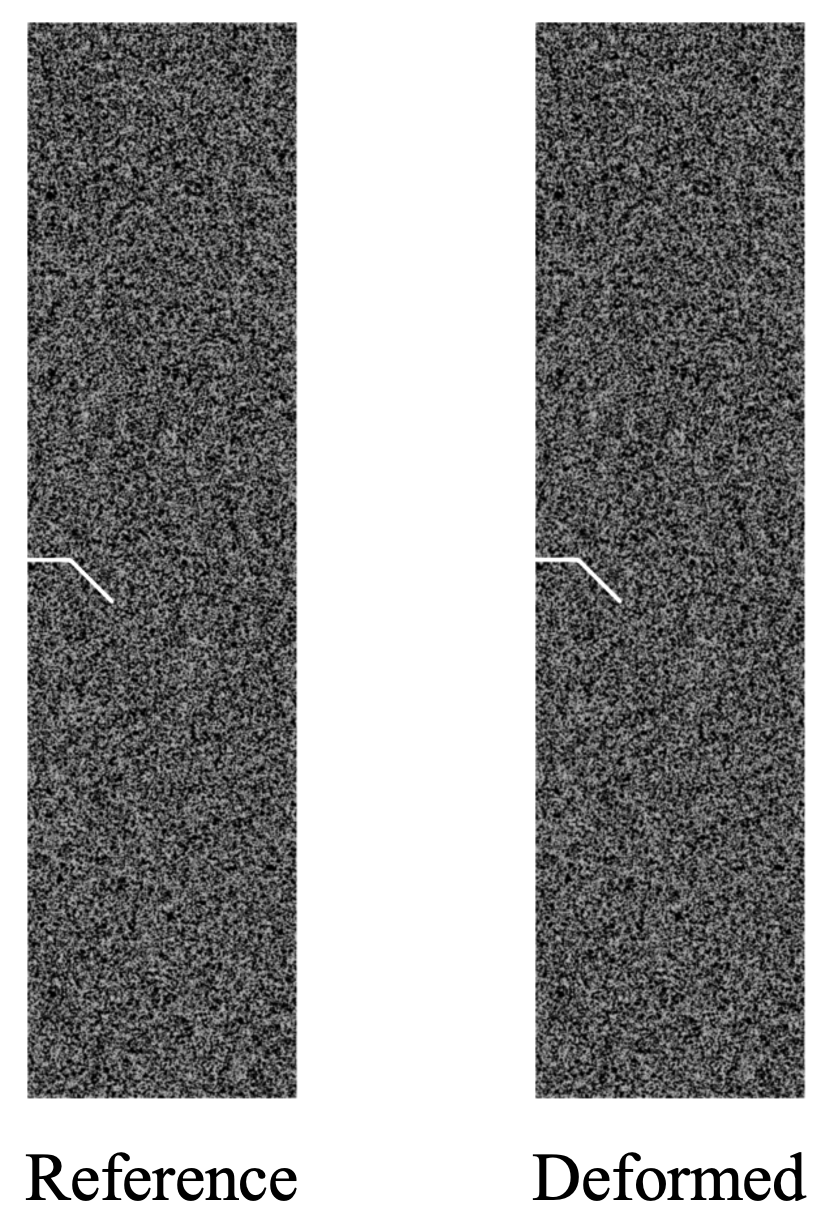}
  \caption{Speckle images}
  \label{fig:deflection_speckle}
\end{subfigure}
\hspace{0.1\textwidth}
\begin{subfigure}[b]{0.2\textwidth}
  \centering
  \makebox[\linewidth][c]{%
    \hspace*{0.5cm}%
    \includegraphics[height=4.5cm,keepaspectratio]{   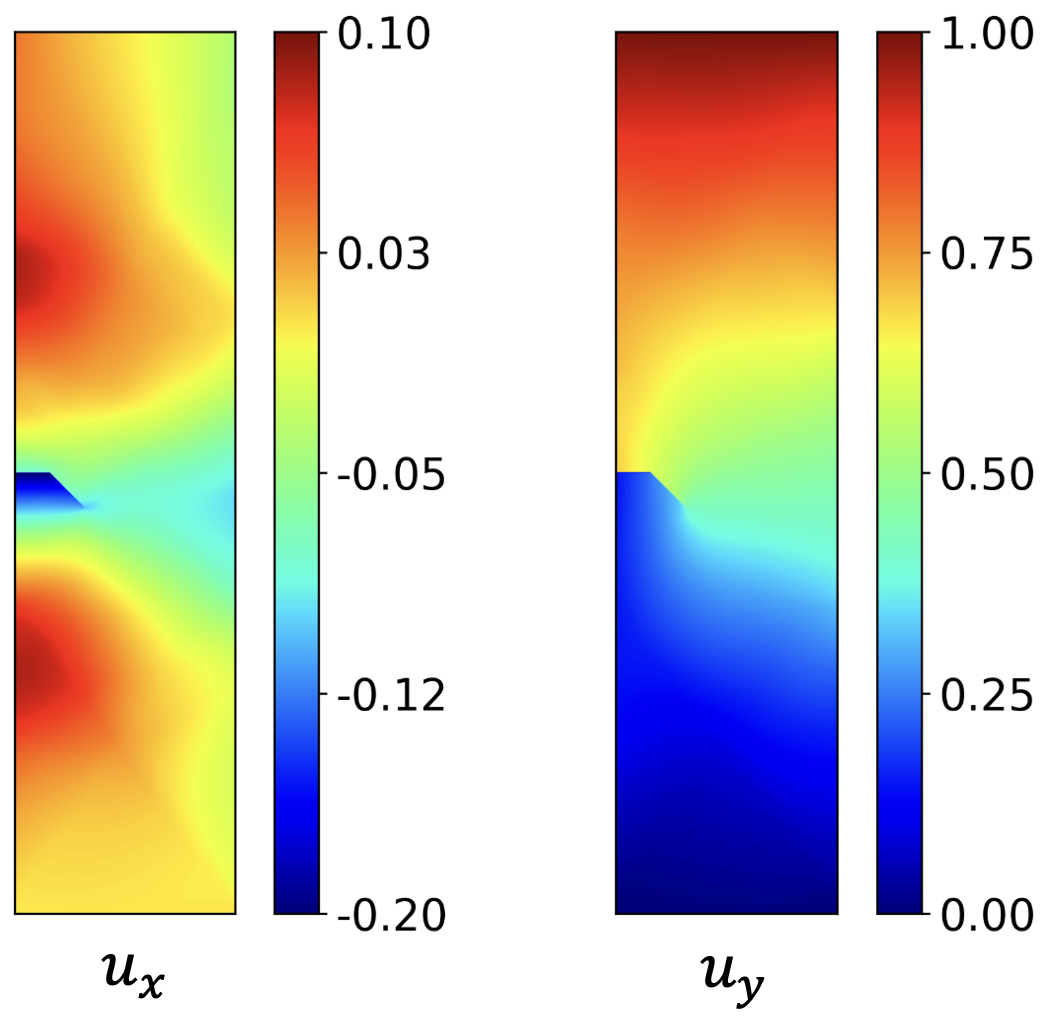}
  }
  \caption{Displacement}
  \label{fig:deflection_ref}
\end{subfigure}

  \caption{ {(a) Reference and deformed images, and (b) ground-truth displacement fields for the case involving crack path deviation.}}
  \label{fig:chap3_deflection}
\end{figure}



\begin{figure}[!htbp]
\centering

\begin{subfigure}[b]{0.35\textwidth}
\centering
\includegraphics[height=4.5cm]{   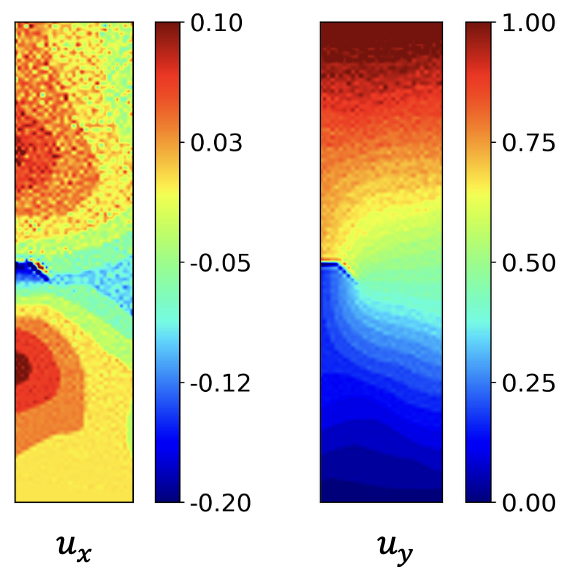}
\caption{DIC results at iteration 1}
\label{fig:deflect_nocrack_1}
\end{subfigure}%
\hspace{0.01\textwidth}%
\begin{subfigure}[b]{0.5\textwidth}
\centering
\includegraphics[height=4.58cm]{   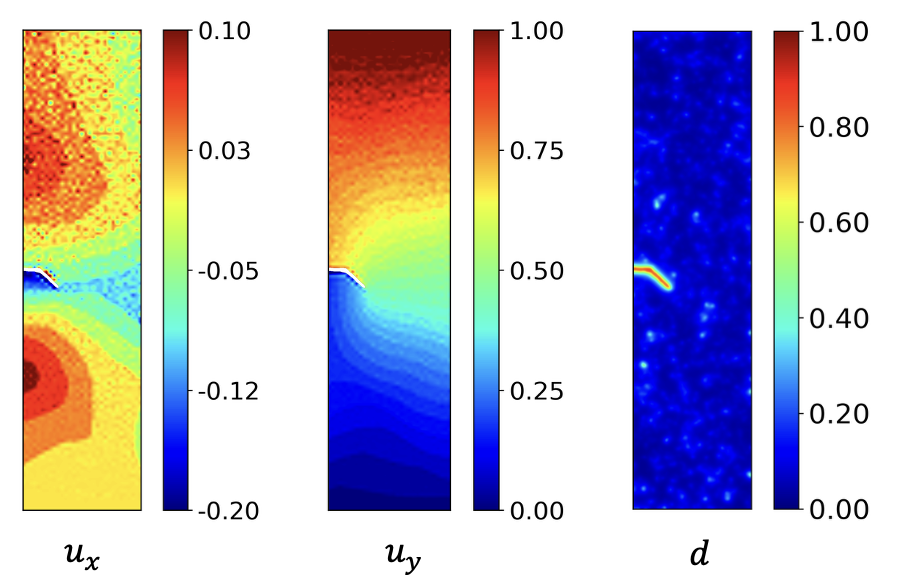}
\caption{PF-DIC results at iteration 1}
\label{fig:deflect_1}
\end{subfigure}

\vspace{0.25cm}



\begin{subfigure}[b]{0.35\textwidth}
\centering
\includegraphics[height=4.5cm]{   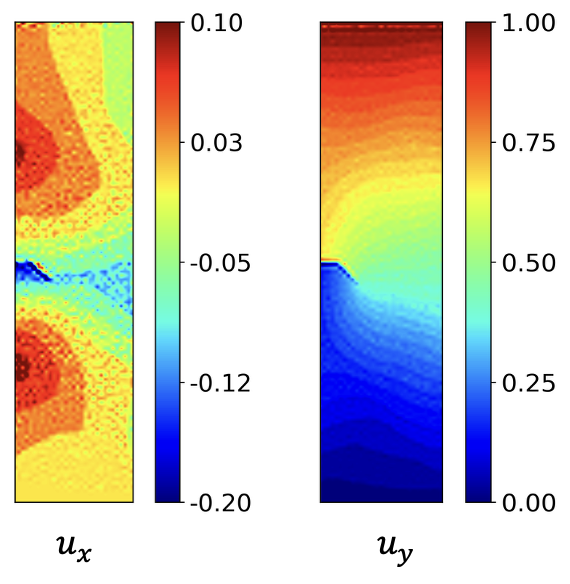}
\caption{DIC converged results (iteration 3)}
\label{fig:deflect_nocrack_3}
\end{subfigure}%
\hspace{0.01\textwidth}%
\begin{subfigure}[b]{0.5\textwidth}
\centering
\includegraphics[height=4.58cm]{   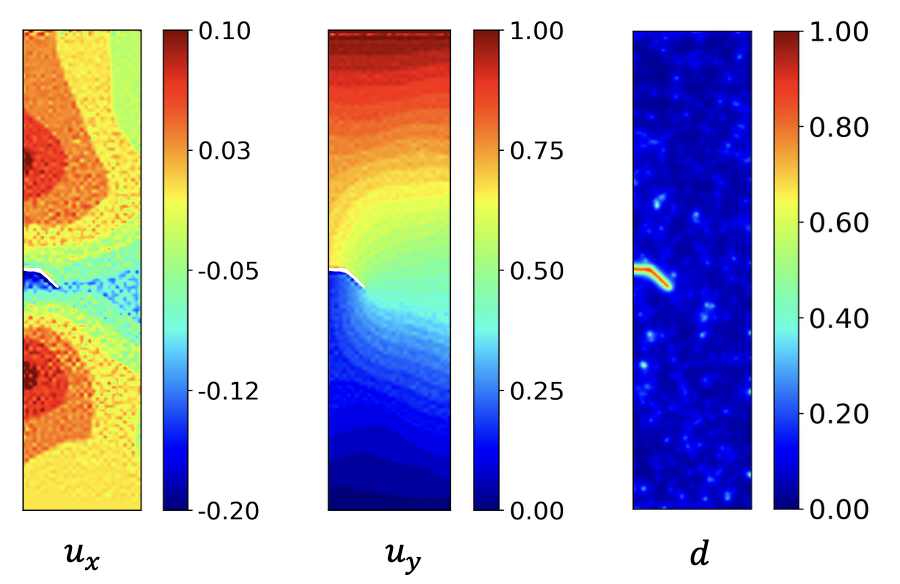}
\caption{PF-DIC converged results (iteration 3)}
\label{fig:deflect_3}
\end{subfigure}

\caption{ {Comparison of standard DIC and the proposed PF-DIC  for the case involving crack path deviation.}}
\label{fig:deflect_iter_results}
\end{figure}

To further confirm the capability of the PF-DIC in terms of measuring more complex cracks, we modified the single edge notched specimen to include a crack path deviation. The newly generated reference displacement and speckle images are illustrated in Fig.~\ref{fig:chap3_deflection}. Using the same hyperparameters in \tablename~\ref{tab:parameters}, we obtained the PF-DIC results in Fig.~\ref{fig:deflect_iter_results}. For comparison purposes, we applied the standard (FE-)DIC to the same speckle images and obtained the results of first three iterations by a modified Gauss-Newton algorithm, as shown in Fig.~\ref{fig:deflect_iter_results}. Again, the PF-DIC converged with few iterations and produced accurate displacement and damage measurements. Particularly, the crack deviation was well identified with the PF-DIC without any manual intervention. We can notice that the overall displacement looks close to that of standard DIC and agrees well with the reference solution, especially for areas far away from the damaged area.


\FloatBarrier

\begin{figure}[!htbp]
\centering

\begin{subfigure}[t]{0.38\textwidth}
\vspace{0pt}
\centering
\parbox[c][2.55cm][c]{\linewidth}{%
    \centering
    \includegraphics[height=2.40cm]{   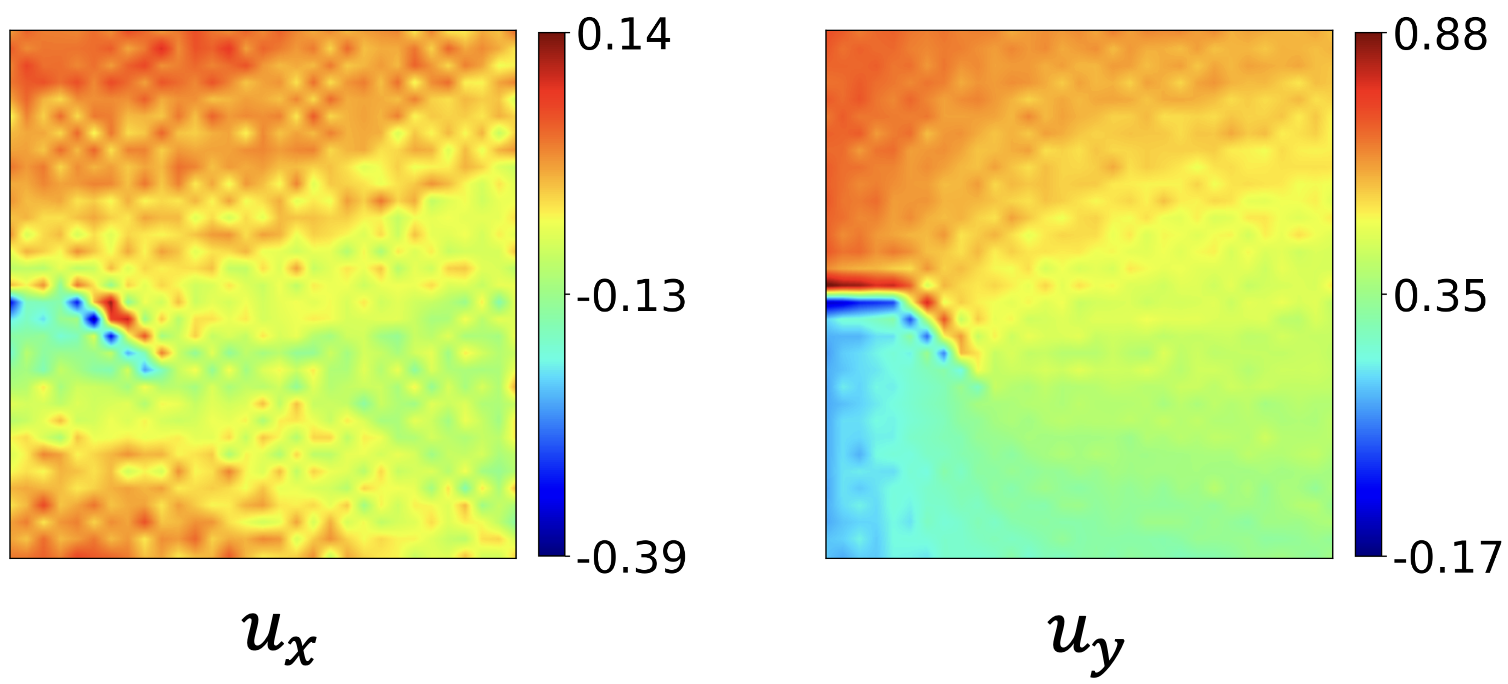}
}
\caption{DIC results near the crack tip at iteration 3}
\label{fig:crack_tip_dic_i3}
\end{subfigure}%
\hspace{0.02\textwidth}%
\begin{subfigure}[t]{0.56\textwidth}
\vspace{0pt}
\centering
\parbox[c][2.55cm][c]{\linewidth}{%
    \centering
    \includegraphics[height=2.53cm]{   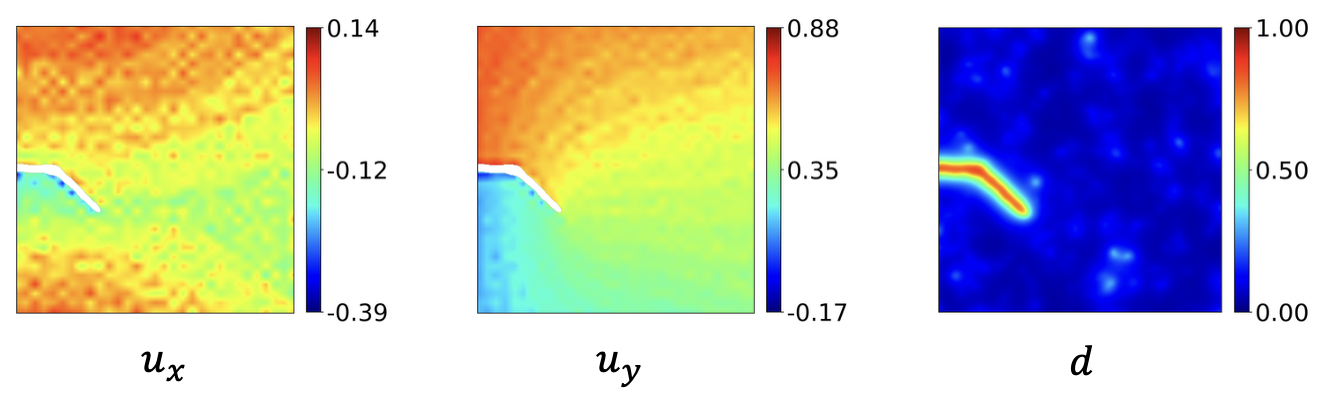}
}
\caption{PF-DIC results near the crack tip at iteration 3}
\label{fig:crack_tip_pfdic_i3}
\end{subfigure}

\caption{ {Comparison of the DIC and PF-DIC results near the deflected crack.}}
\label{fig:deflect_iter_cracktip_results}
\end{figure}

\FloatBarrier

\begin{figure}[!htbp]
  \centering
  \includegraphics[height=6.5cm]{   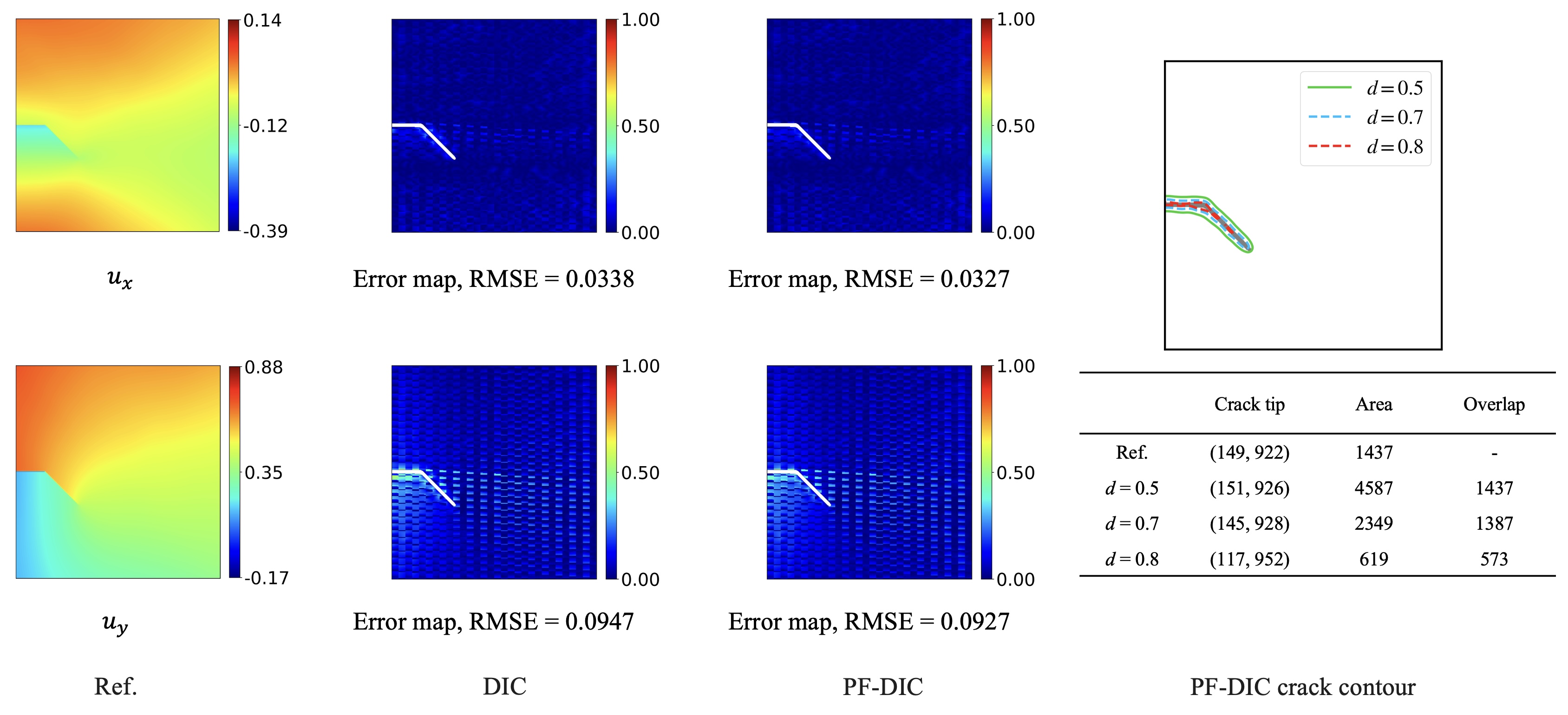}
  \caption{ {Error of the displacement fields around the deflected crack and verification of the crack contour identified by the PF-DIC.}}
  \label{fig:deflection_errormap}
\end{figure}

Fig.~\ref{fig:deflect_iter_cracktip_results} shows a detailed comparison of the displacement field around the crack from the same example. As we can see, the crack created a strong discontinuity in the displacement field across the cracked region, along with  {some} spurious displacement fluctuations. This is why a mask is usually needed to remove the relevant pixels from the speckle images for standard DIC techniques to avoid such non-physical measurements in the cracked regions, based on a prejudgment of users.  {Furthermore, it seems difficult to identify the precise location of the crack tip from the displacement measurements in the current example}.  In the PF-DIC,  {the mask is not needed as the damage field plays a similar role in removing or reducing} the influence of  spurious displacement fluctuations near the crack.
 {To better visualize this influence of crack on displacement measurements, we computed again the pixel-wise diplacement error for both DIC and PF-DIC, as shown in Fig.~\ref{fig:deflection_errormap}. We can see that the crack may cause erroneous displacement near the crack, which seems slightly improved by the PF-DIC. This improvement could be more significant if the error is larger, as shown in the next example. Regarding the damage measurement of PF-DIC, we compared  the crack contours of three different values of $d$ with the reference crack, in terms of crack tip, area, and overlap. Again, the contour corresponding to $d=0.5$ appears to most closely match the reference: both the detected crack topology and the crack tip location  match the pre-defined crack. This confirms the capability of the PF-DIC for automatic damage detection and crack tracking. }

\FloatBarrier

\subsection{Tensile specimen with crack propagation}
\label{sec:crack_propagation}
\begin{figure}[!htbp]
  \centering

\begin{subfigure}[b]{0.2\textwidth}
  \centering
  \includegraphics[height=4.5cm,keepaspectratio]{ 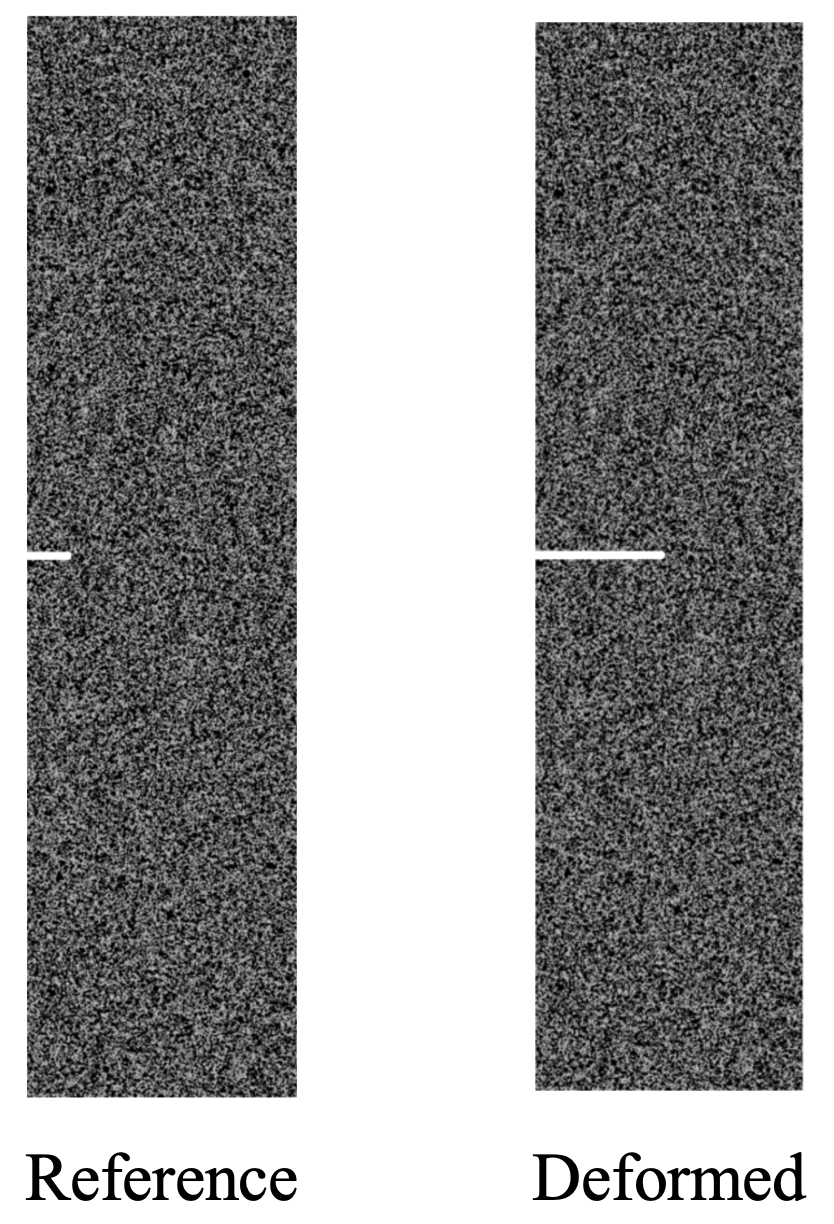}
  \caption{Speckle images}
  \label{fig:extension_speckle}
\end{subfigure}
\hspace{0.1\textwidth}
\begin{subfigure}[b]{0.2\textwidth}
  \centering
  \makebox[\linewidth][c]{%
    \hspace*{0.5cm}%
    \includegraphics[height=4.5cm,keepaspectratio]{ 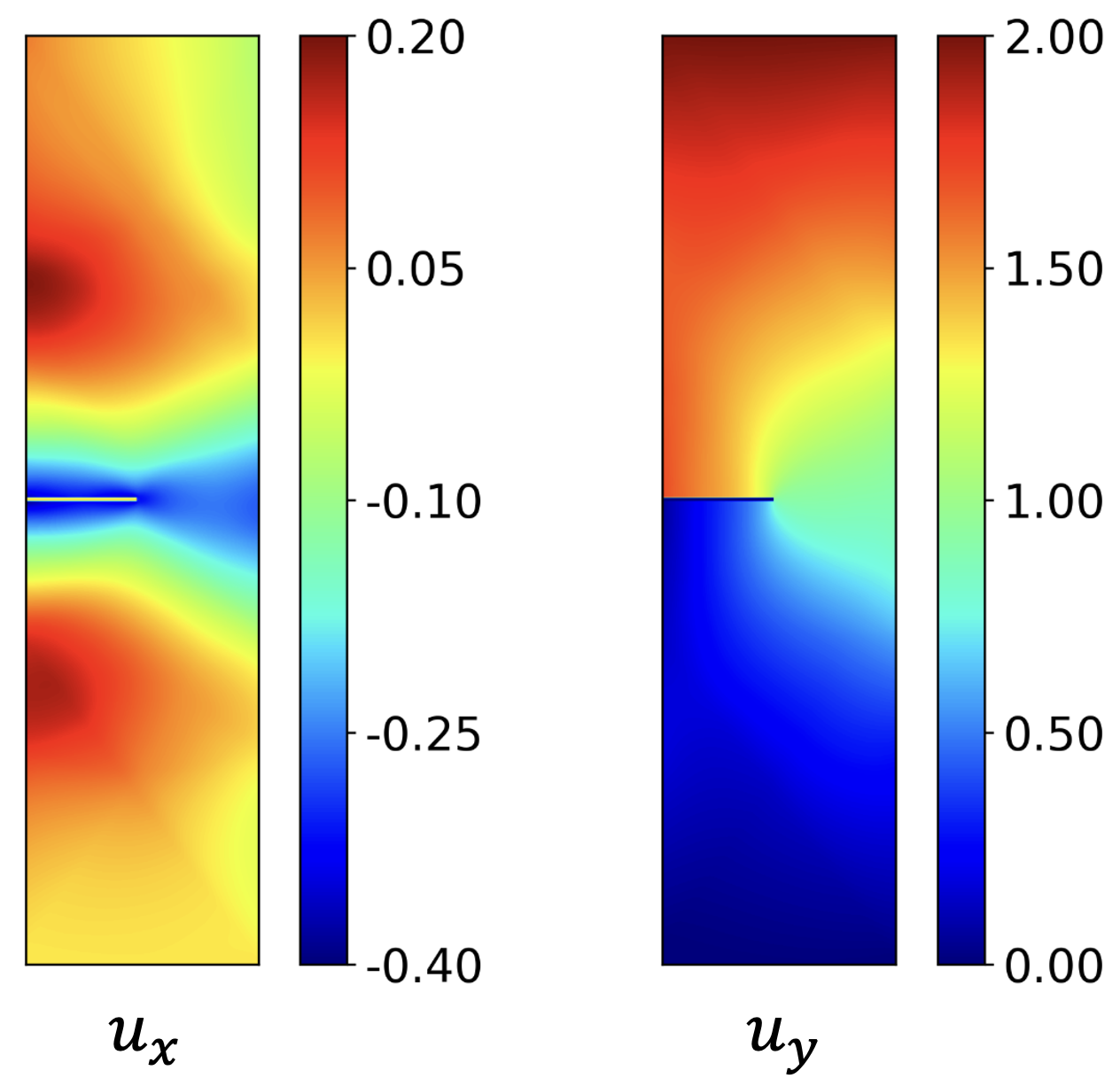}
  }
  \caption{Displacement}
  \label{fig:extension_dis_ref}
\end{subfigure}

  \caption{ {(a) Reference and deformed speckle images, and (b) ground-truth displacement fields for the tensile specimen with crack propagation.}}
  \label{fig:chap3_ref_extension}
\end{figure}

This example aimed to test the capability of the PF-DIC in the case of crack propagation. This could be important for characterizing the cracking behaviors of materials. We used the same single edge notched tensile test as the previous example, but considered an extension of crack in the deformed stage, as shown in Fig.~\ref{fig:chap3_ref_extension}. The relevant displacement between the reference and deformed images was calculated using the sample with the extended crack, which was assumed to represent an equilibrium state after crack propagation. 

\begin{figure}[!htbp]
\centering

\begin{subfigure}[t]{0.35\textwidth}
\centering
\includegraphics[height=4.5cm]{ 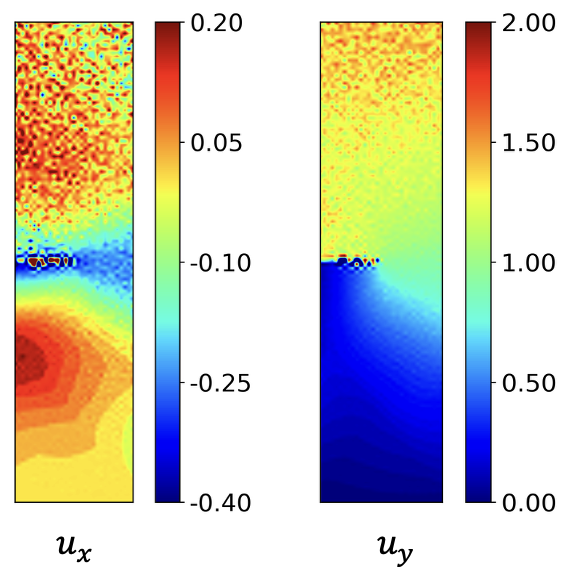}
\caption{DIC results at iteration 1}
\label{fig:chap3_compare_crack_extension_a}
\end{subfigure}
\hspace{0.002\textwidth}
\begin{subfigure}[t]{0.5\textwidth}
\centering
\includegraphics[height=4.58cm]{ 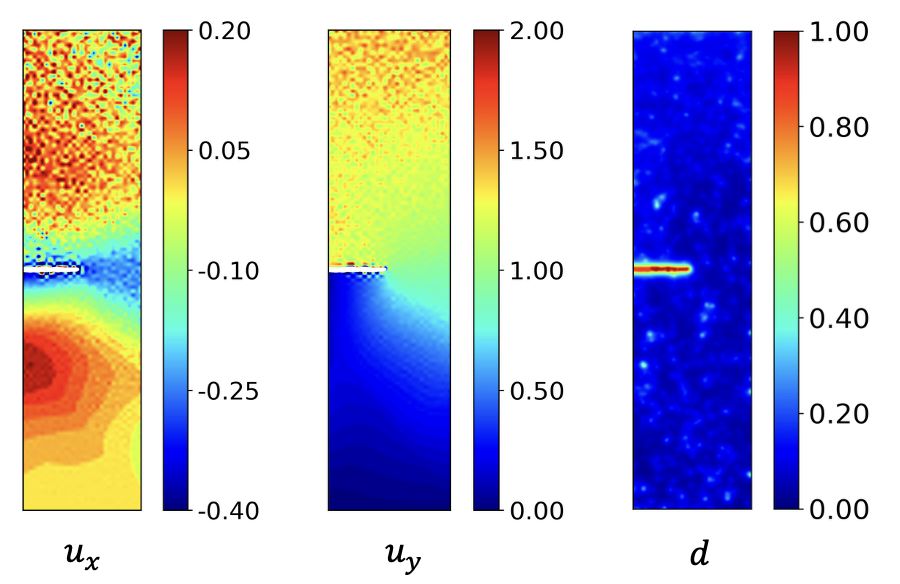}
\caption{PF-DIC results at iteration 1}
\label{fig:chap3_compare_crack_extension_b}
\end{subfigure}

\vspace{0.25cm}

\begin{subfigure}[t]{0.35\textwidth}
\centering
\includegraphics[height=4.5cm]{ 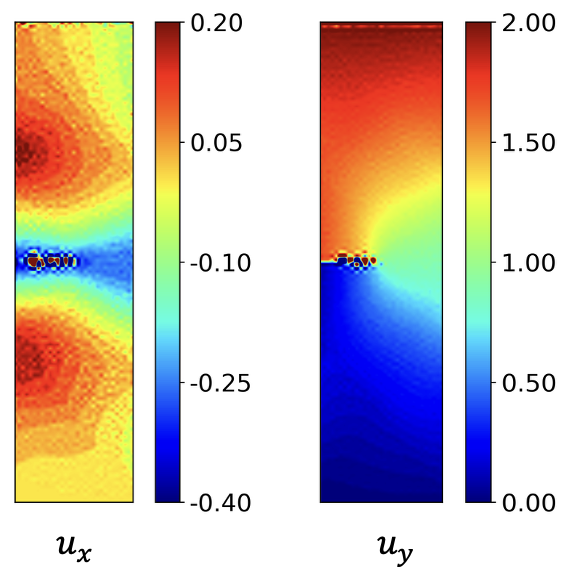}
\caption{DIC converged results (iteration 3)}
\label{fig:chap3_compare_crack_extension_e}
\end{subfigure}
\hspace{0.002\textwidth}
\begin{subfigure}[t]{0.5\textwidth}
\centering
\includegraphics[height=4.58cm]{ 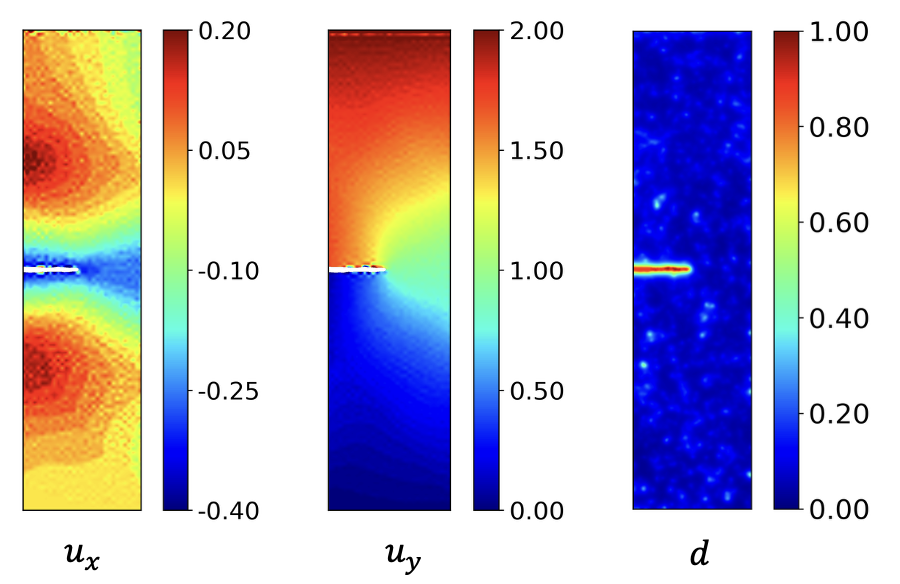}
\caption{PF-DIC converged results (iteration 3)}
\label{fig:chap3_compare_crack_extension_f}
\end{subfigure}

\caption{ {Comparison of standard DIC and the proposed PF-DIC in the case of crack propagation.}}
\label{fig:chap3_compare_crack_extension}
\end{figure}

Using the same parameter settings as previously, including the mesh size and the damage driving force energy and its associated hyperparameters,  we applied the PF-DIC to analysis the displacement and damage fields from the speckle images. As a comparison, we also applied the standard DIC to the images without any predefined mask. The results are illustrated in Fig.~\ref{fig:chap3_compare_crack_extension}. As we can see, the PF-DIC results again converged quickly within a few iterations. More importantly, the extended crack was well identified  by the PF-DIC automatically. The final crack length and the overall displacement agree with the reference in Fig.~\ref{fig:chap3_ref_extension}. Compared with standard DIC results, we can see a significant improvement of displacement smoothness around the crack by the PF-DIC.  We can notice the amplified spurious displacement fluctuations around the crack, due to the propagated crack and the resulting mismatch between the reference and deformed speckle images. These fluctuations cannot be resolved by the iterative solution scheme and may become even worse  if the iterative procedure continues.   {In practice, a mask should be carefully applied and adapted to cracked samples for standard DIC, to avoid the influence of the spurious displacement on cracks.}

\begin{figure}[htbp]
  \centering

  \begin{subfigure}[t]{0.4\textwidth}
    \centering
    \includegraphics[width=\textwidth]{ 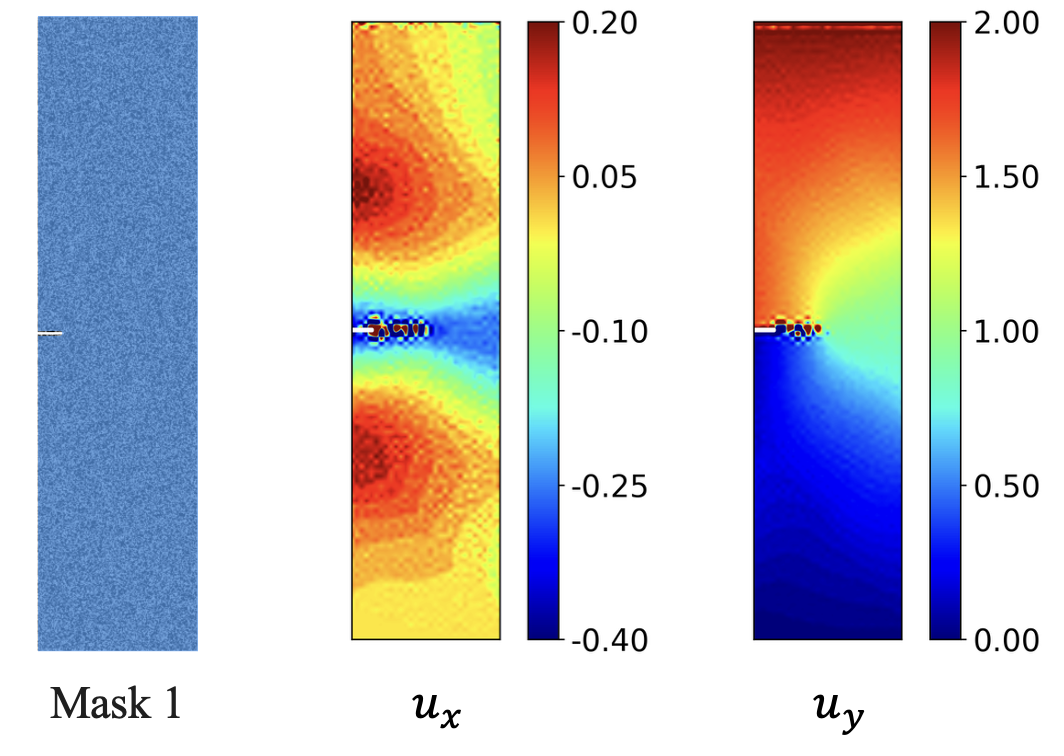}
    \caption{ {Masked DIC-1}}
    \label{fig:mask1}
  \end{subfigure}
  \hspace{0.1\textwidth}
  \begin{subfigure}[t]{0.4\textwidth}
    \centering
    \includegraphics[width=\textwidth]{ 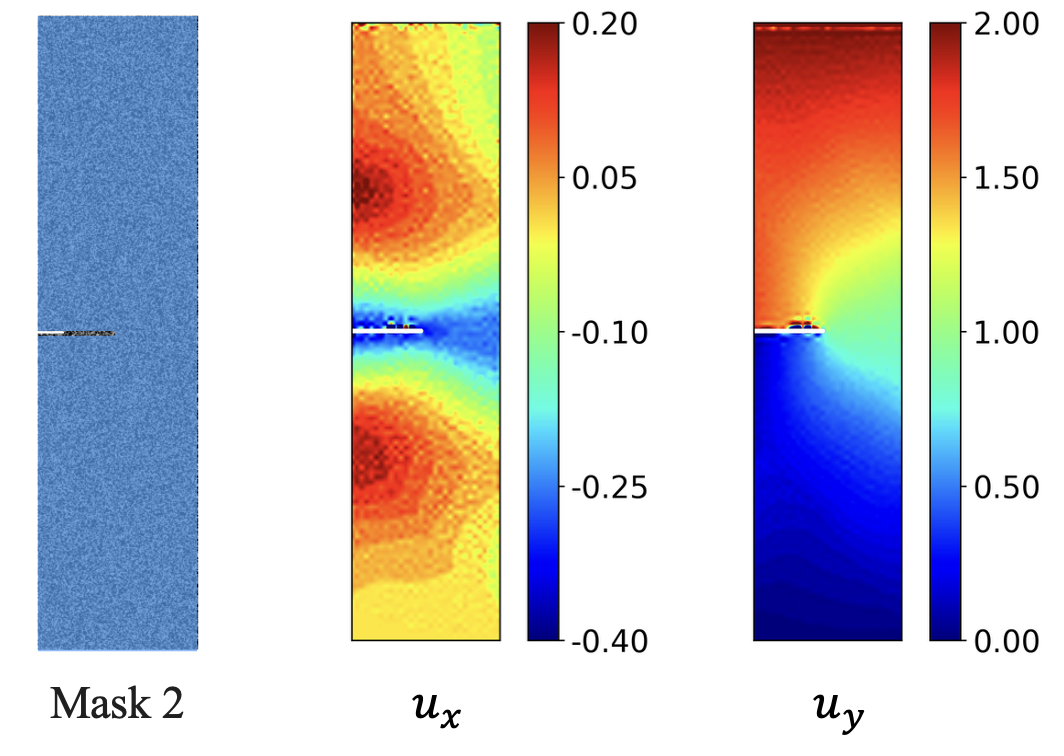}
    \caption{ {Masked DIC-2}}
    \label{fig:mask2}
  \end{subfigure}

  \caption{ {Results of masked DIC based on two pre-defined masks. Mask 1 corresponds to the initial notched specimen by removing the pixels on the crack; Mask 2 corresponds to that of the extended crack.}}
  \label{fig:mask_extraction}
\end{figure}
\FloatBarrier

 {To further evaluate the performance of PF-DIC against masked DIC, we applied two different masks to the standard DIC: one excluding the initial crack and another excluding the final crack.  The results obtained using the two different masks are referred to  as masked DIC-1 and masked DIC-2, respectively, and are illustrated in Fig.~\ref{fig:mask_extraction}. As expected, spurious displacement fluctuations can be reduced by a mask; and the mask 2 seems better than the mask 1 in this example, as it removes more pixels regarding the extended crack. }

\begin{figure}[htbp]
  \centering

  \begin{subfigure}{0.9\textwidth}
    \centering
    \includegraphics[width=\textwidth]{ 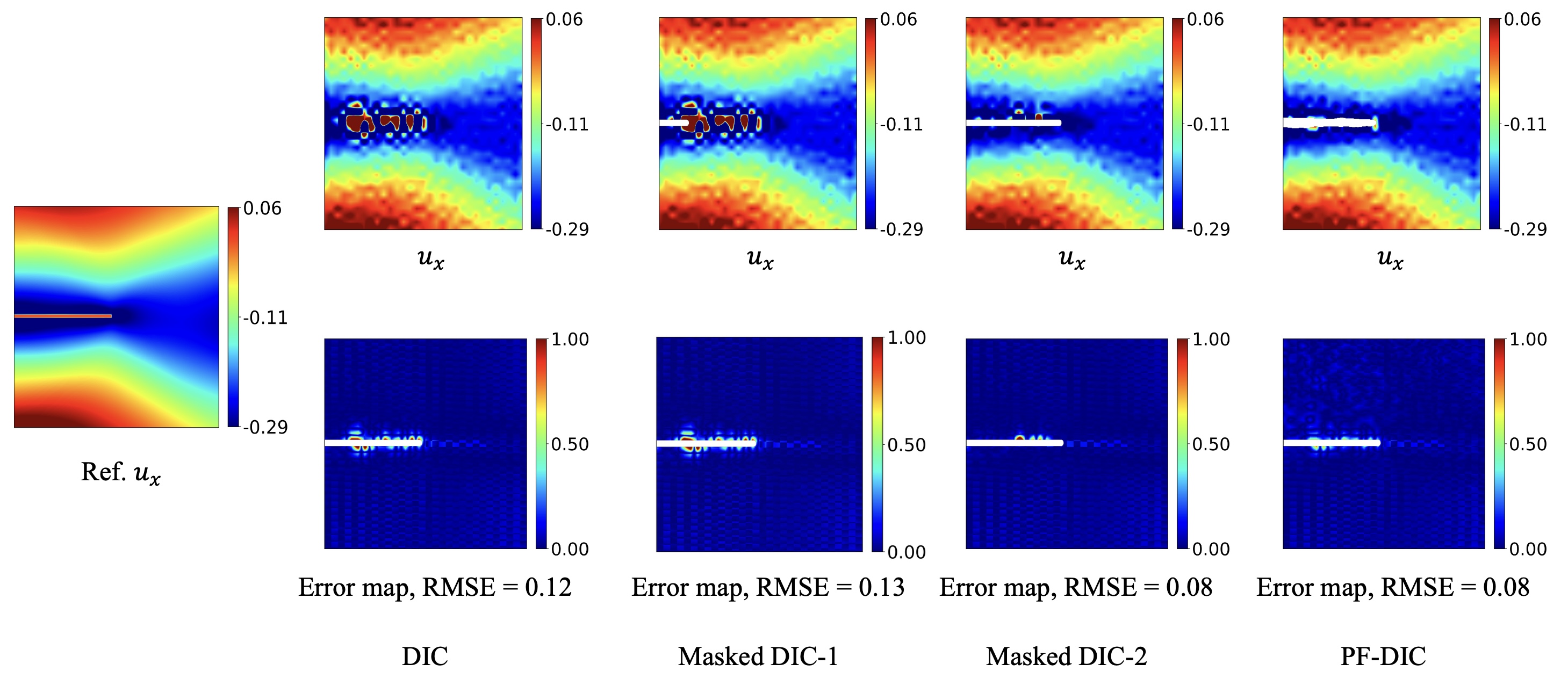}
    \caption{ {Error map $\frac{|u-u_{\text{ref.}}|}{|u_{\text{ref.},\max}|}$ and RMSE for $u_x$}}
    \label{fig:ext_ux}
  \end{subfigure}

  \vspace{0.3cm}

  \begin{subfigure}{0.9\textwidth}
    \centering
    \includegraphics[width=\textwidth]{ 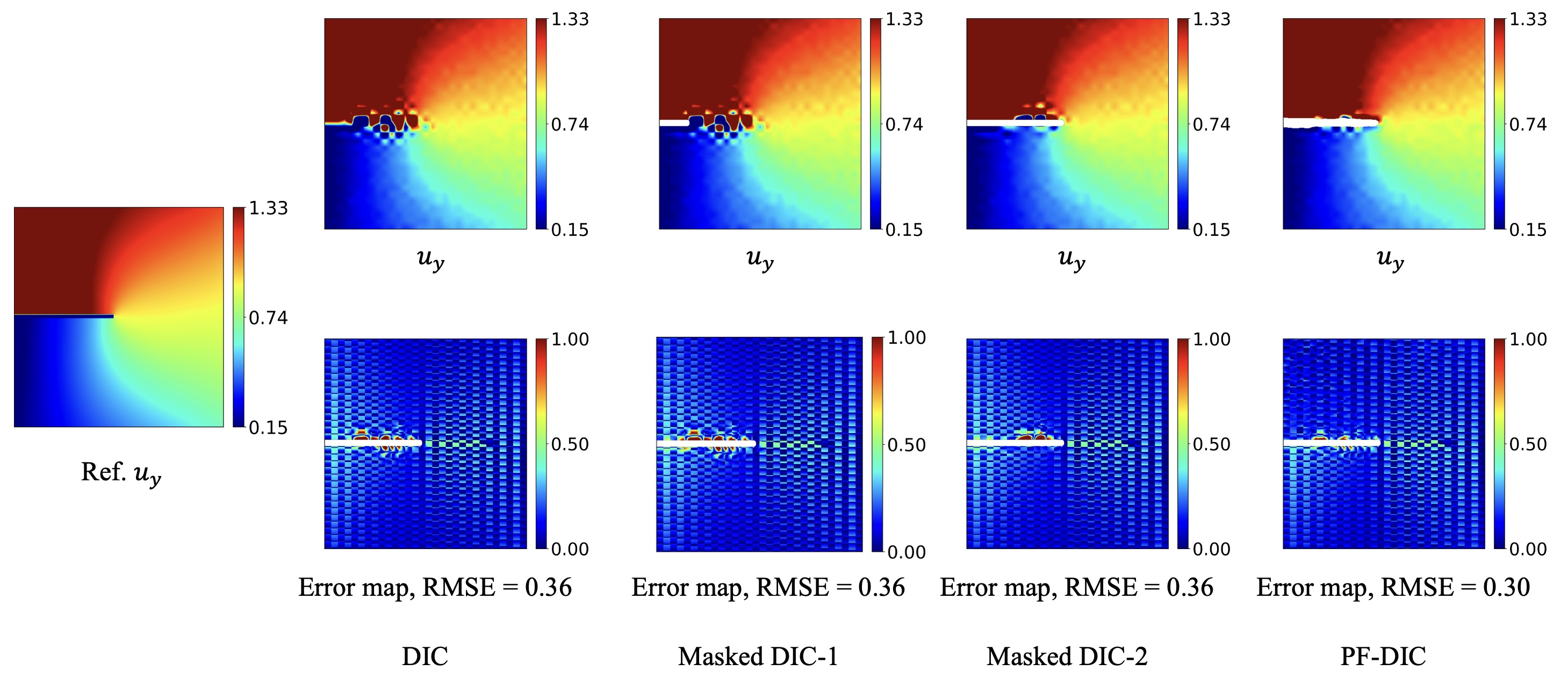}
    \caption{ {Error map $\frac{|u-u_{\text{ref.}}|}{|u_{\text{ref.},\max}|}$ and RMSE for $u_y$}}
    \label{fig:ext_uy}
  \end{subfigure}

  \caption{ {Displacement near the crack  for the (unmasked) DIC, masked DIC, and PF-DIC in the case of crack propagation.}}
  \label{fig:extension_displacement}
\end{figure}
\FloatBarrier

 {A more detailed comparison of the different methods is illustrated in Fig.~\ref{fig:extension_displacement}. In this figure, we compared the original (unmasked) DIC, the two masked DIC, and the (unmasked) PF-DIC for displacement measurements near the crack. Additionally, we excluded all the pixels that represent the actual crack from the error map for better illustrating the discrepancy. The RMSE was then calculated using the remaining part of displacement values. As we can see, an appropriate mask can indeed significantly improve the displacement measurements around the crack, compared to the (unmasked) DIC, especially for $u_x$. The PF-DIC can get similar accuracy for $u_x$ without any mask. Regarding $u_y$, the improvement of masks on standard DIC seems relatively limited, whereas the PF-DIC still improved the displacement near the crack. This improvement may be attributed to the diffusive representation of cracks, which effectively reduces the influence of incomplete speckle information in the vicinity of the crack.} 

 {Applying a larger mask can certainly further suppress the spurious fluctuations in standard DIC. However, this would require users to make a prejudgment regarding the trade-off between the loss of displacement information and the local measurement accuracy near the crack. This again highlights the advantage of PF-DIC, as it eliminates the need for such prejudgments to accommodate the crack-induced displacement discontinuity. In addition, for propagating cracks, it seems difficult to predefine a mask covering all possible crack paths, making masked DIC less practical for real-time displacement monitoring.}

\begin{figure}[!htbp]
  \centering
  \includegraphics[height=6.5cm]{ 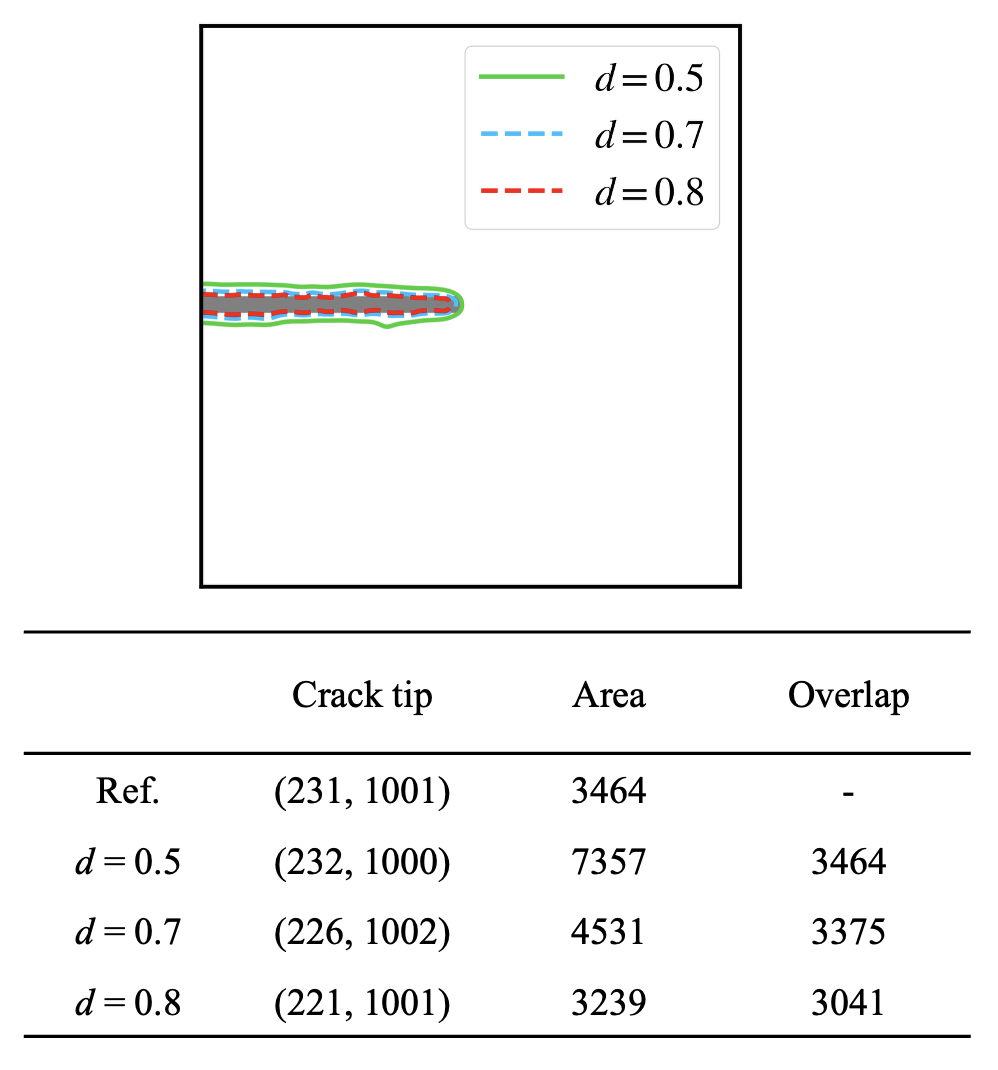}
  \caption{ {The accuracy of the damage measurement by PF-DIC.}}
  \label{fig:extension_contour}
\end{figure}

 {To verify the damage measurement of PF-DIC, we compared the crack contours obtained for different values of $d$ with the  extended crack (reference), as shown in Fig.~\ref{fig:extension_contour}. Consistent with the previous examples, all the contours closely match the reference in terms of crack tip and overleaping area.  Among them, the contour corresponding to $d=5$ seems the most appropriate one for representing the actual crack, despite its relatively larger area. Overall, these results confirm the capability of PF-DIC to reliably track crack propagation.}

\FloatBarrier

\subsection{ {Three-point bending test and selective damage detection}}
\label{sec:bending}
\begin{figure}[!htbp] 
\centering 
\begin{subfigure}[t]{0.45\textwidth}
\hspace{-0.18\textwidth}
\centering
\includegraphics[width=\linewidth]{ 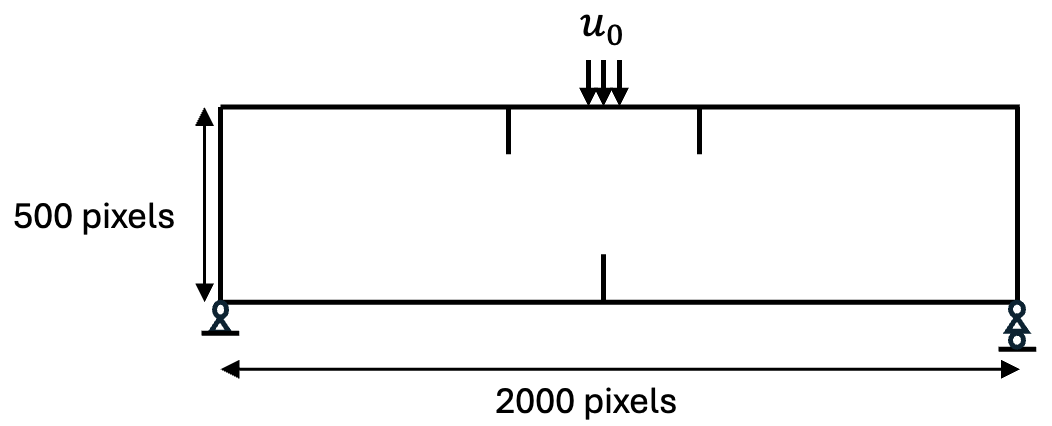}
\caption{Geometry and boundary conditions}
\label{fig:three_point_geometry}
\end{subfigure}

\vspace{0.3cm}
\begin{subfigure}[t]{0.40\textwidth} 
\centering 
\includegraphics[width=\linewidth]{ 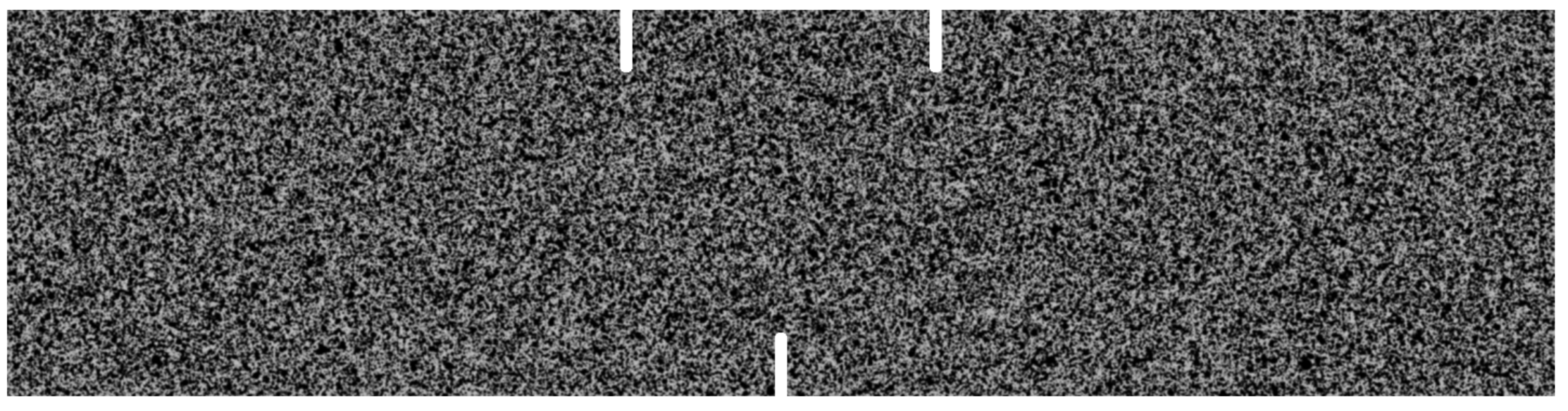} \caption{Reference} 
\label{fig:three_point_f} 
\end{subfigure} 
\vspace{0.2cm} 
\begin{subfigure}[t]{0.4\textwidth} 
\centering 
\includegraphics[width=\linewidth]{ 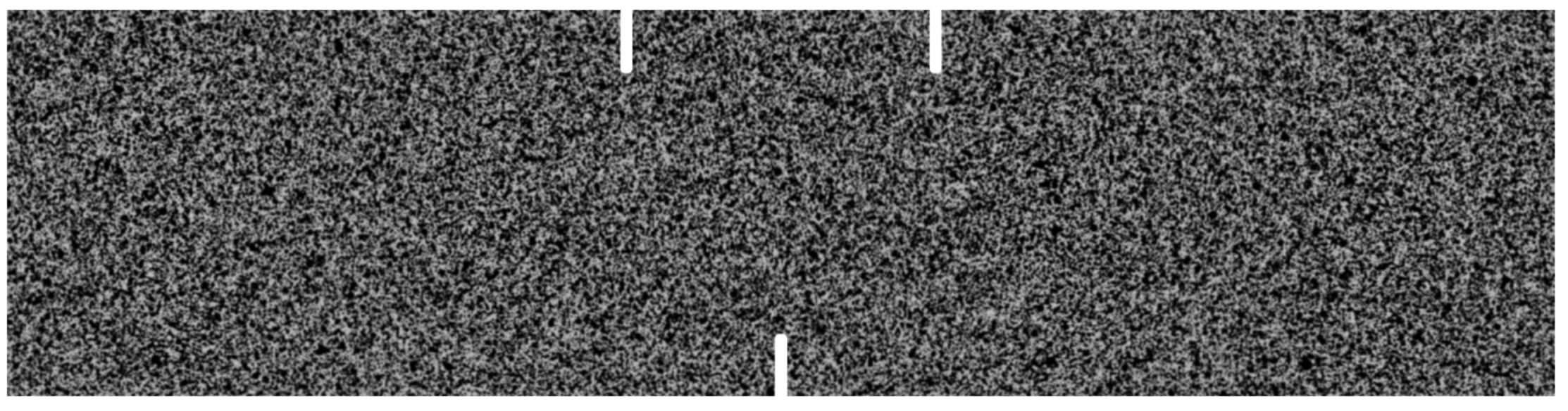} \caption{Deformed} 
\label{fig:three_point_g} 
\end{subfigure} 
\vspace{0.2cm} 
\begin{subfigure}[t]{0.4\textwidth} 
\centering 
\includegraphics[width=\linewidth]{ 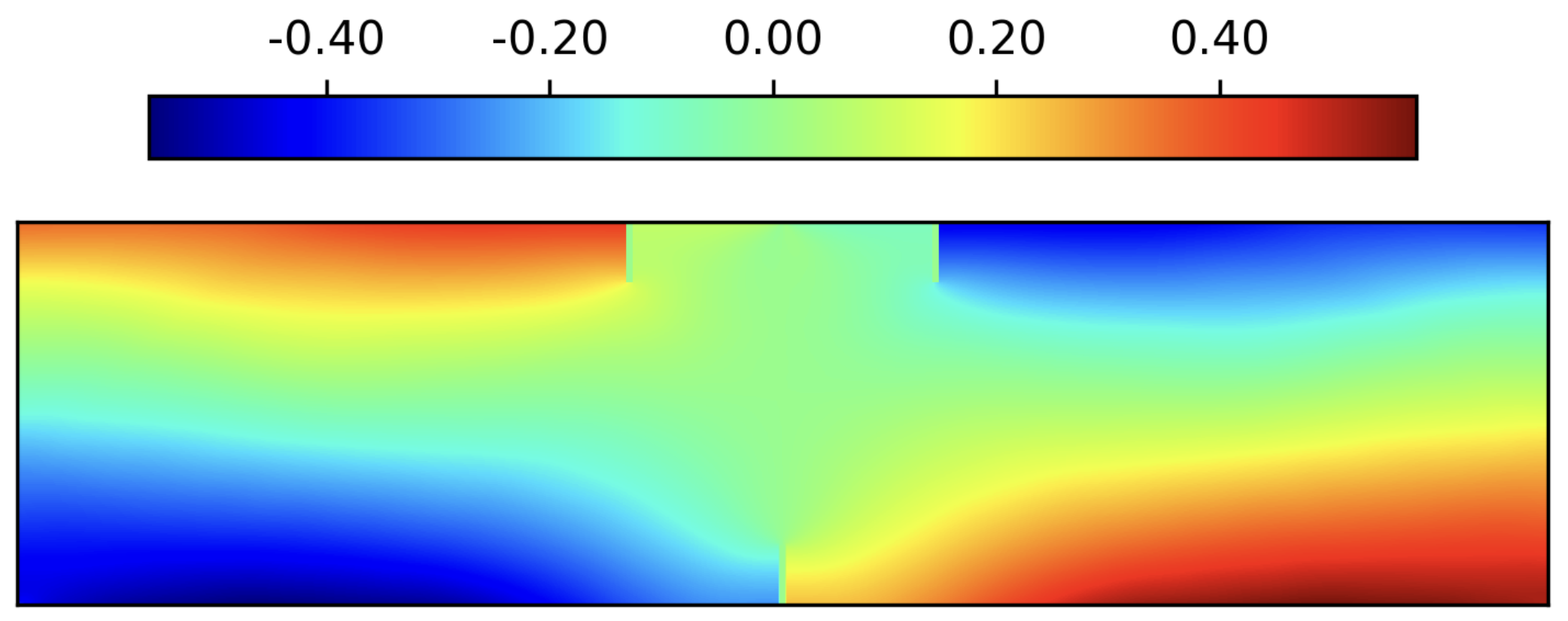} \caption{$u_x$} 
\label{fig:three_point_ux} 
\end{subfigure} 
\vspace{0.2cm} 
\begin{subfigure}[t]{0.4\textwidth} 
\centering 
\includegraphics[width=\linewidth]{ 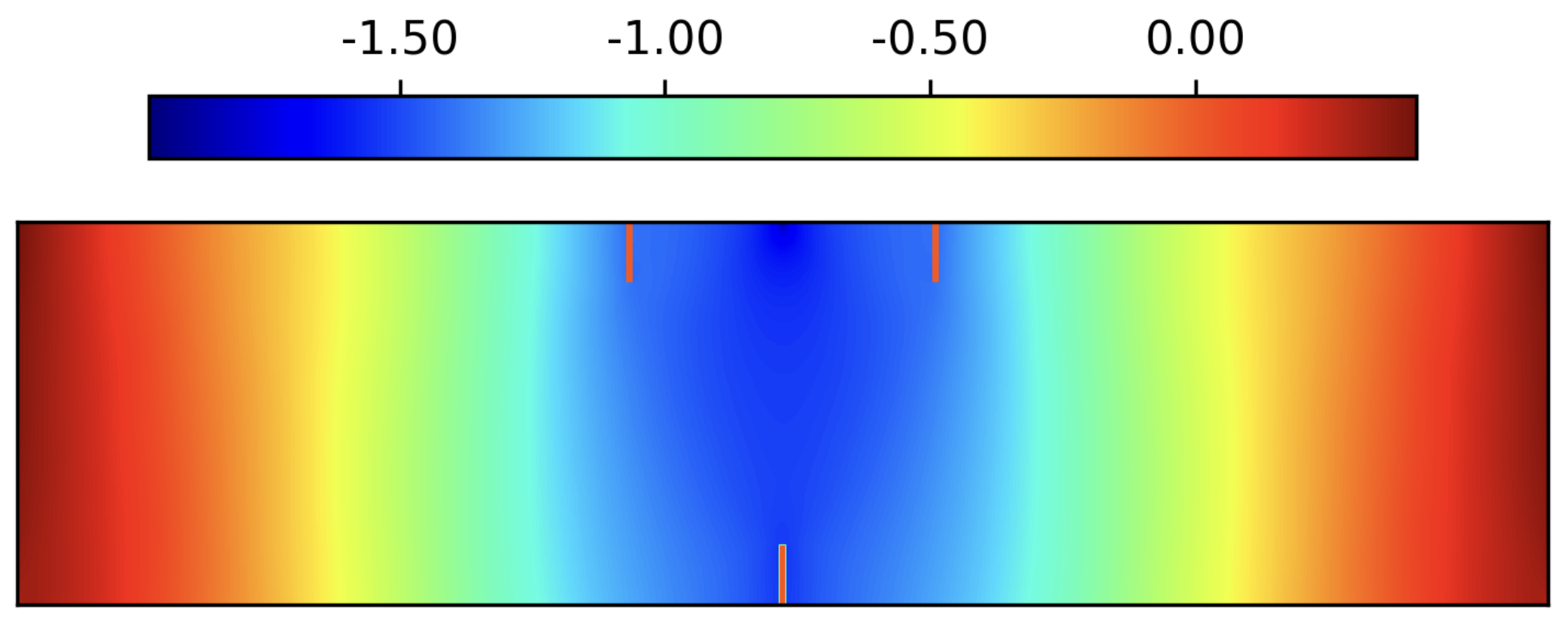} \caption{$u_y$} 
\label{fig:three_point_uy} 
\end{subfigure} 
\caption{Geometry and boundary conditions, reference and deformed speckle images, and  ground-truth displacement fields for the three-point bending test.}
\label{fig:chap3_three_point_bending} 
\end{figure}

To validate the performance of PF-DIC under different loading conditions, we considered a three-point bending test, as shown in Fig.~\ref{fig:chap3_three_point_bending}. A vertical displacement was prescribed at the top boundary to mimic the loading condition in the middle of the top edge. Three preexisting cracks were assumed in the specimen: one located in the lower region and two located in the upper region. The ROI in the speck images contains $2000\times500$ pixels. 

Again, we used an FE mesh of 4-node elements with element size $h_u = 16$ (pixel) for  displacement and $h_d=1$ (pixel) for damage for the PF-DIC. Eq.~\eqref{eq:fracture_energy_w_background} was adopted for the damage driving force energy $\mathit{\Phi}$, with the hyperparameters $w_{1,2}$, $a_{1,2}$, and $c_{1,2}$ given in Table~\ref{tab:3pb_parameters}.


\begin{table}[H]
\centering
\setlength{\abovecaptionskip}{2pt}
\caption{Hyperparameters of $\mathit{\Phi}$  in three-point bending tests.}
\label{tab:3pb_parameters}
\small
\begin{tabular}{ccccccc}
\toprule
Parameter & $w_1$ & $a_1$ & $c_1$ & $w_2$ & $a_2$ & $c_2$ \\
\midrule
Value & 0.5 & $3\times10^{-8}$ & $1\times10^{-9}$ & 0.005 & 0.6 & 0.05 \\
\bottomrule
\end{tabular}
\end{table}


\begin{figure}[htbp]
  \centering

  \begin{subfigure}{\textwidth}
    \centering
    \includegraphics[width=\textwidth]{ 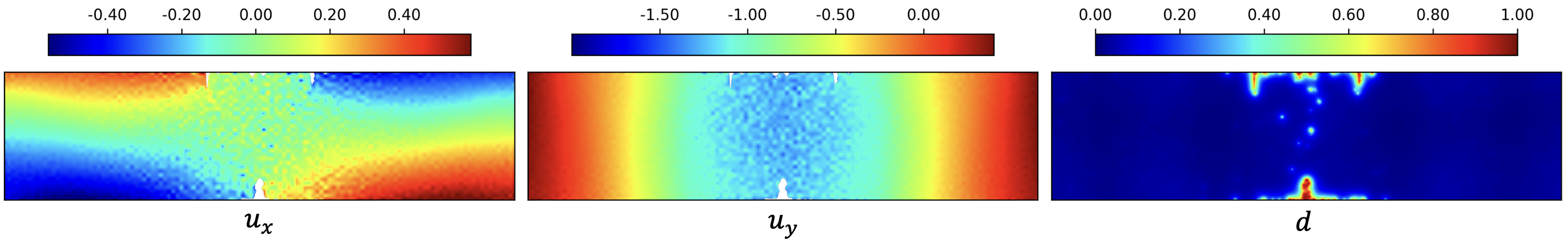}
    \caption{PF-DIC results at iteration 1}
  \end{subfigure}

  \vspace{0.3cm}

  \begin{subfigure}{\textwidth}
    \centering
    \includegraphics[width=\textwidth]{ 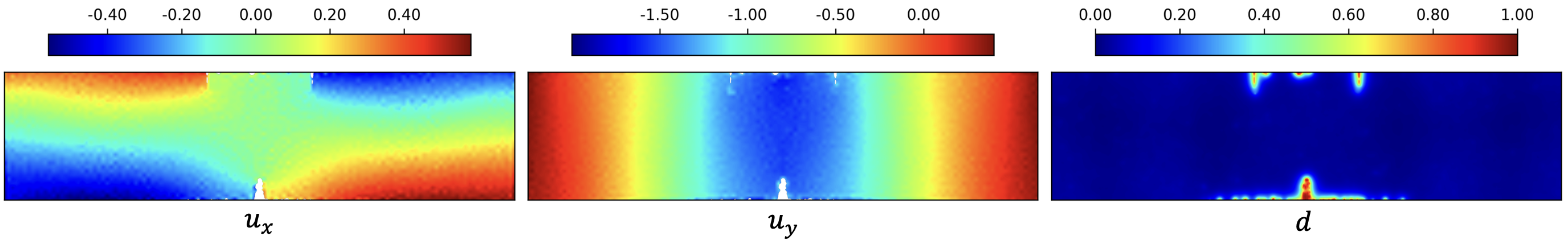}
    \caption{PF-DIC converged results (iteration 3)}
  \end{subfigure}

  \caption{Displacement and damage fields measured by PF-DIC for the three-point bending test.}
  \label{fig:chap3_3pb_original_results}

\end{figure}

The displacement and damage fields obtained by PF-DIC for this three-point bending test are illustrated in Fig.~\ref{fig:chap3_3pb_original_results}. We can observe similar iterative behavior in the solutions as before: the iterative procedure converged with a few iterations  and improved both the smoothness and accuracy of the overall displacement field, and the final converged result matched the ground-truth displacement. Regarding the damage measurement, the PF-DIC was able to detect all three cracks located at different positions within the specimen. The spurious fluctuations in the damage field present in the first iteration were effectively removed over successive iterations. This is due to the improvement of displacement measurement over iterations and the appropriate choice of the damage driving force energy. Therefore, these results confirmed again the unique capability of PF-DIC for effectively measuring displacement and damage fields in a unified framework. 

In the next, we used the same bending specimen to test another important capability of the PF-DIC, i.e., selective detection of cracks. As mentioned earlier, the PF-DIC can selectively detect cracks and damage that are considered critical to users by modifying the damage driving force energy $\mathit{\Phi}$. In this three-point bending test, we can expect that the upper portion of the specimen experiences compressive strain, whereas the lower portion is subjected to tensile strain. Given this context, the crack in the lower portion may be identified as the critical crack with the greatest potential for propagation. This identification or diagnostic process can be done automatically by the PF-DIC, by considering the following definition of damage driving force energy  
\begin{equation}
\mathit{\Phi}=\mathcal{A}_1
\!\left(
\boldsymbol{\varepsilon}^+, r,f
\right)+
\mathcal{A}_2\!\left(f\right), \quad \text{with}\quad \boldsymbol{\varepsilon}^+\coloneqq \langle \varepsilon_{xx}\rangle^+,
\label{eq:fracture_energy_select_modif}
\end{equation}
where $\langle\cdot\rangle^+ = \frac{1}{2}(\cdot+|\cdot|)$, and $\varepsilon_{xx}$ denotes the normal strain component in $x$-direction. This definition allows to consider only the Mode I fracture of the crack and identify the most dangerous crack among the existing ones.

\begin{figure}[htbp]
  \centering

  \begin{subfigure}{\textwidth}
    \centering
    \includegraphics[width=\textwidth]{ 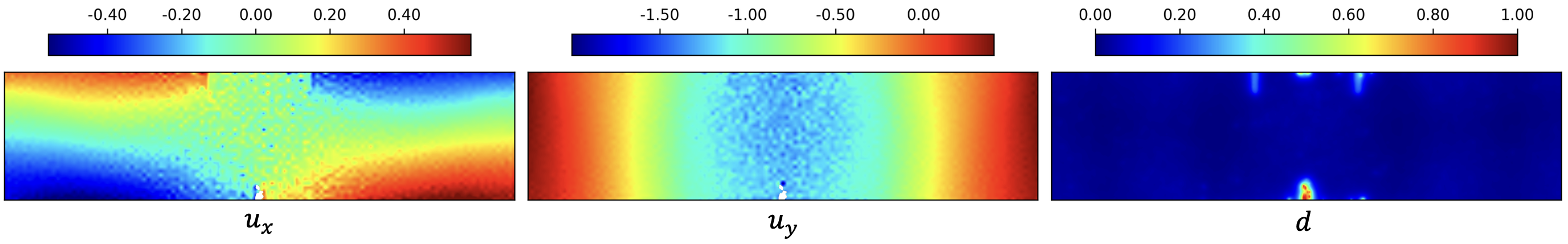}
    \caption{PF-DIC results at iteration 1 with the modified $\mathit{\Phi}$}
  \end{subfigure}

  \vspace{0.3cm}

  \begin{subfigure}{\textwidth}
    \centering
    \includegraphics[width=\textwidth]{ 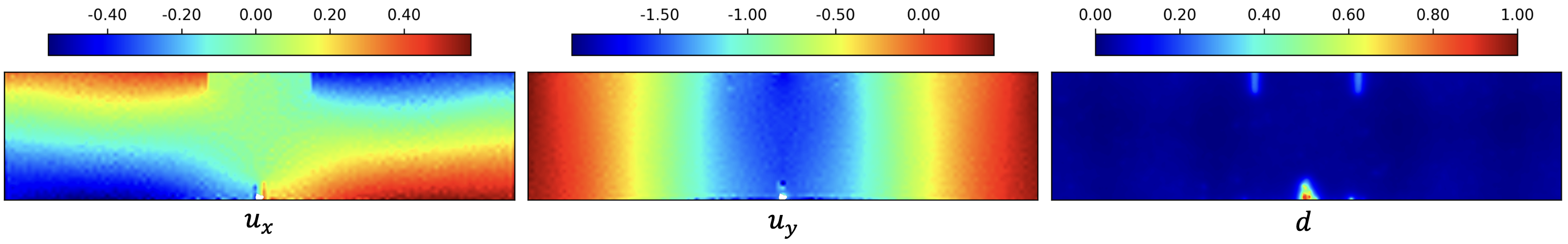}
    \caption{PF-DIC converged results (iteration 3) with the modified $\mathit{\Phi}$}
  \end{subfigure}

  \caption{Displacement and damage fields measured by PF-DIC using the modified damage driving force energy.}

  \label{fig:chap3_3pb_modified_results}

\end{figure}

This modified definition of $\mathit{\Phi}$ can use the same hyperparameters as in Table~\ref{tab:3pb_parameters}. The corresponding results of PF-DIC are illustrated in Fig.~\ref{fig:chap3_3pb_modified_results}. In this case, the PF-DIC identified only the crack in the lower portion of the specimen, clearly distinguishing it as the critical crack, while the overall measurement of displacement remained similar to that of the previous case. This demonstrated  the diagnostic capability of PF-DIC for identifying critical cracks and damage in specimens. In more general cases, the stress state and crack morphology can become highly complex, making the identification of critical cracks among multiple cracks particularly challenging. In this context, the diagnostic capability of PF-DIC can become especially valuable.


\subsection{Experimental images of fractured samples}
\begin{figure}[htbp]
  \centering
  \includegraphics[width=0.5\textwidth]{ 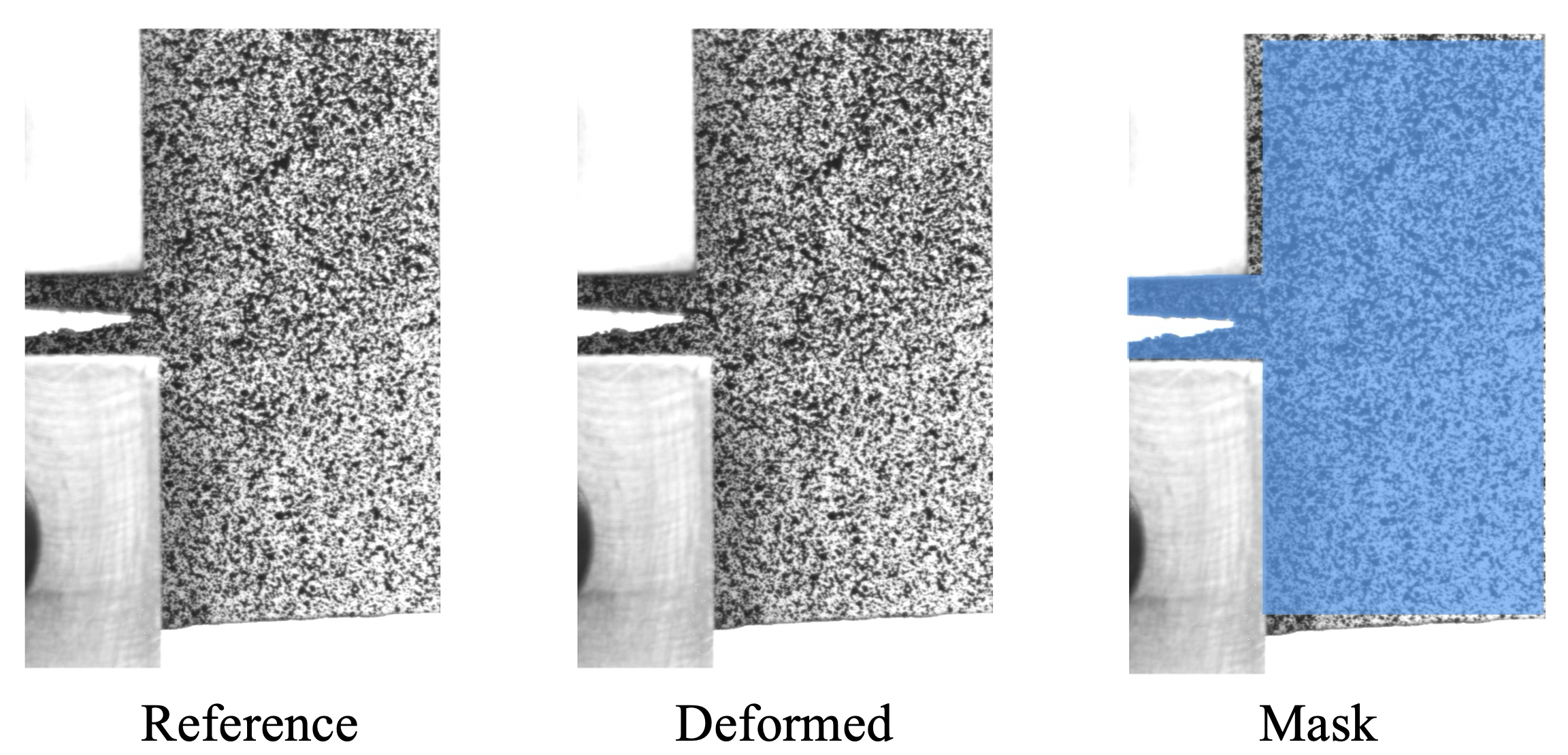}
  \caption{ {Reference and deformed images for $t_{465-480}$ and a potential mask for standard DIC.}}
  \label{fig:exp_ref}
\end{figure}

\begin{figure}[htbp]
  \centering

  \begin{subfigure}{0.45\textwidth}
    \centering
    \includegraphics[width=\textwidth]{ 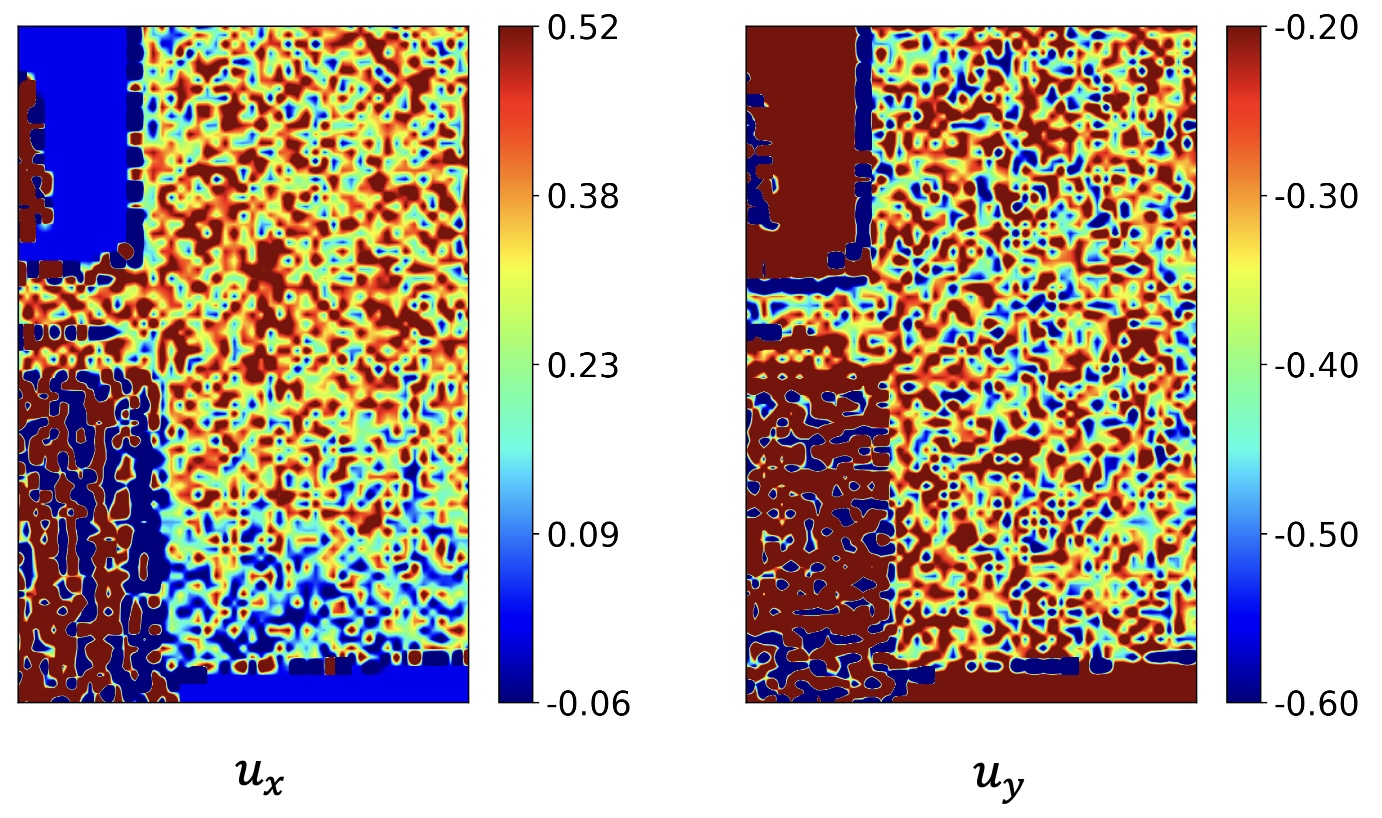}
    \caption{ {Displacement measurements from unmasked DIC}}
    \label{fig:0465_u_dic}
  \end{subfigure}

  \vspace{0.3cm}

  \begin{subfigure}{0.45\textwidth}
    \centering
    \includegraphics[width=\textwidth]{ 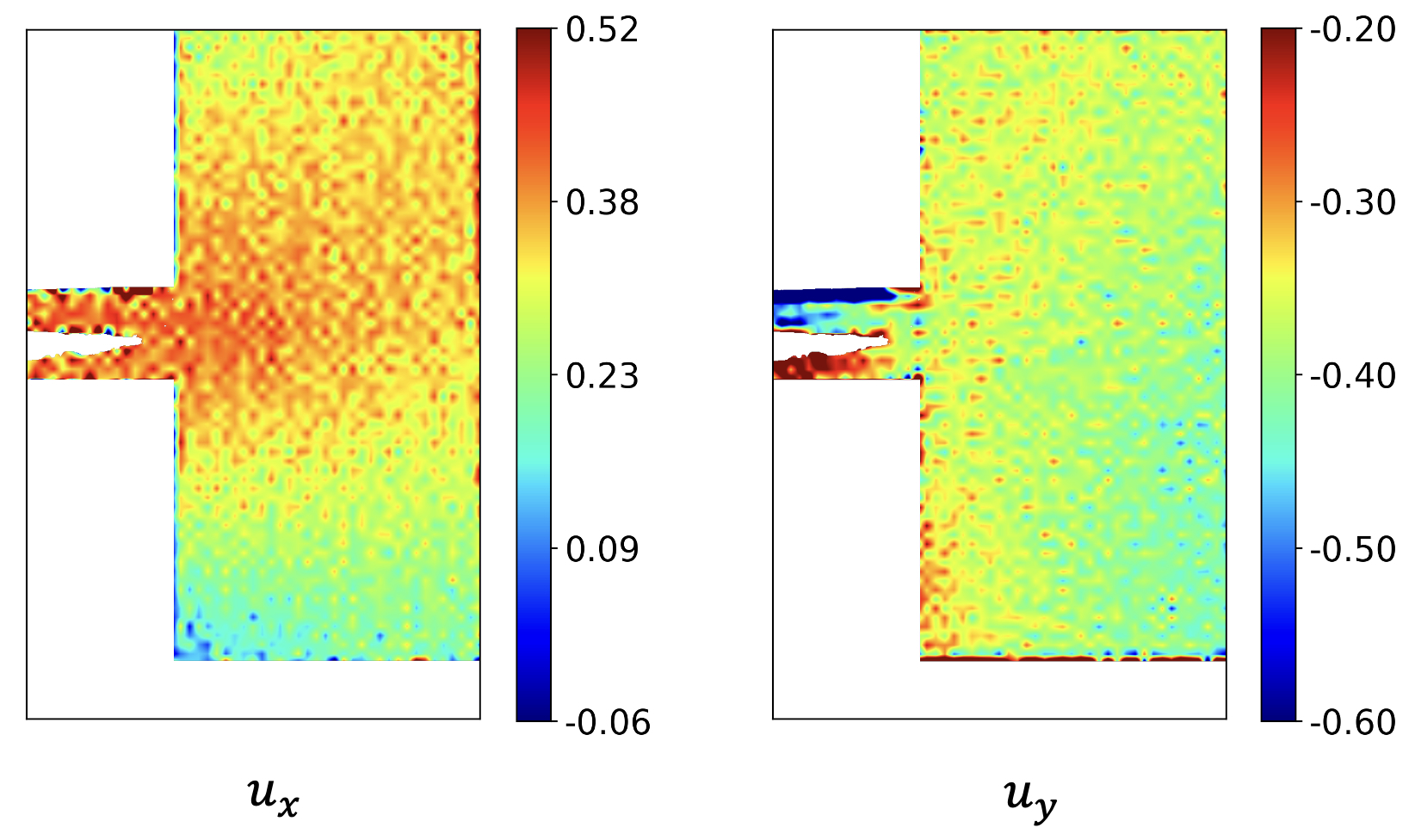}
    \caption{ {Displacement measurements from masked DIC}}
    \label{fig:0465_u_maskdic}
  \end{subfigure}

  \vspace{0.3cm}

  \begin{subfigure}{0.7\textwidth}
    \centering
    \includegraphics[width=\textwidth]{ 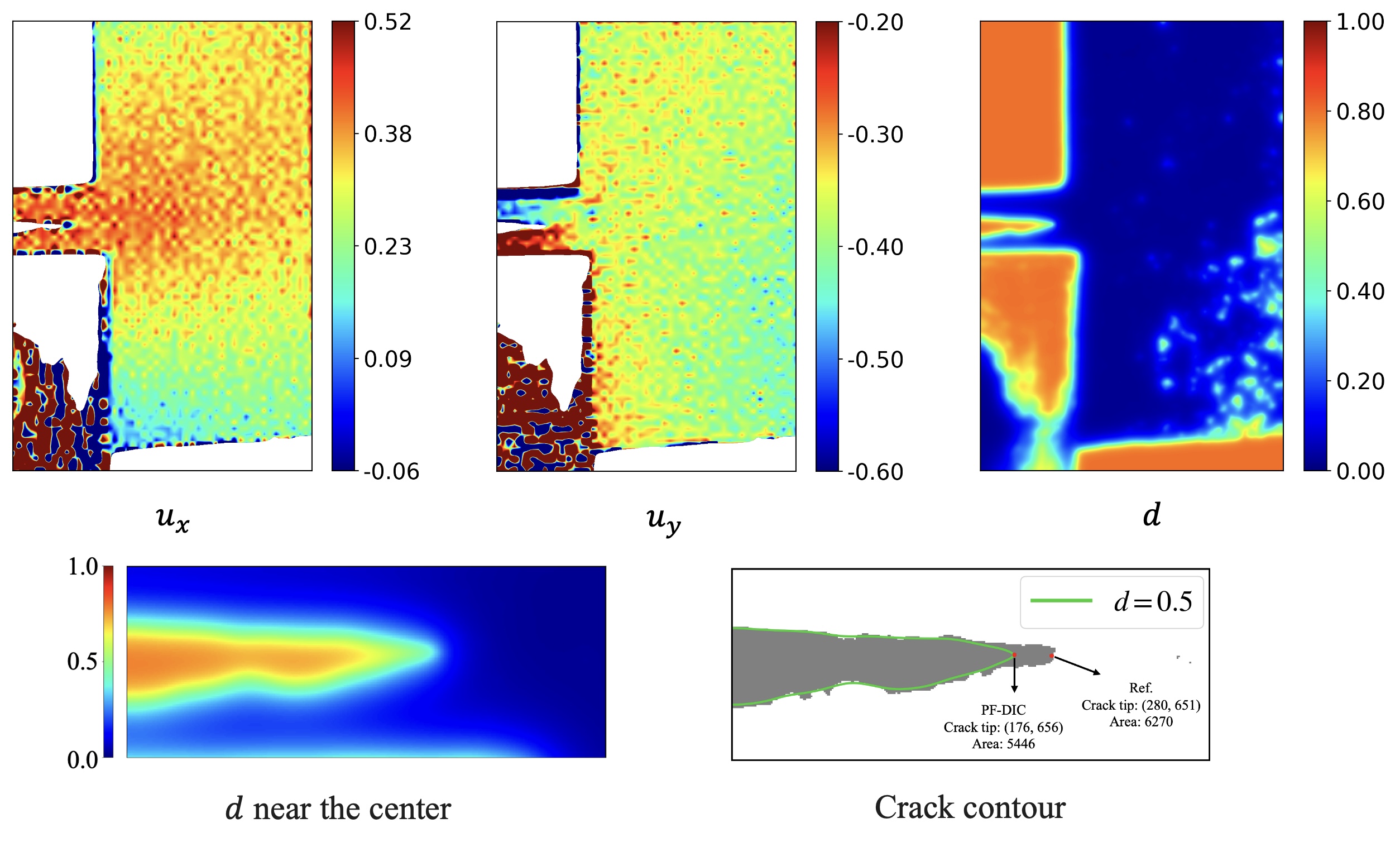}
    \caption{ {Displacement and damage measurements from PF-DIC}}
    \label{fig:0465_u_pfdic}
  \end{subfigure}

  \caption{ {Comparison of (a) standard unmasked DIC, (b) masked DIC, and (c) PF-DIC measurements for $t_{465-480}$.}}
  \label{fig:0465_u_comparison}
\end{figure}

Finally, we tested the PF-DIC using experimental speckle images. We relied on a set of experimental images of fractured samples that represent the fracture of a heterogeneous medium.  Details of the relevant experimental setup and speckle pattens can be found in the literature (e.g., \cite{yang2019augmented}). The image data are available on the GitHub website \cite{yang2020}, comprising 38 speckle images of $800\times1228$ pixels recorded at different time instances. Here, we keep the original ID of each image \cite{yang2020} and use the notation $t_{A-B}$ to denote the time interval between $t_{A}$ and $t_{B}$. Therefore, for a given time interval $t_{A-B}$, the image taken at $t_{A}$ becomes the reference image and that of $t_{B}$ becomes the deformed image.  {Fig.~\ref{fig:exp_ref}} shows an example of the reference and deformed images of the fractured sample taken from the time interval  $t_{465-480}$.  {For the analysis, we considered a ROI of $ 800 \times 1200$ pixels in the middle of the image frame, and applied the standard (unmasked) DIC and PF-DIC with the same $h_u$ and $h_d$ as previously. We used again Eq.~\eqref{eq:fracture_energy_w_background} for the damage driving force energy $\mathit{\Phi}$, whose hyperparameters are given in Table~\ref{tab:experimental_parameters}.  For comparison purposes, we also applied a mask excluding the damaged area and background pixels, as shown in Fig.~\ref{fig:exp_ref}, to obtain masked DIC results.} 
 {The results of different methods for $t_{465-480}$ are illustrated in Fig.~\ref{fig:0465_u_comparison}}.

\FloatBarrier

\begin{table}[H]
\centering
\setlength{\abovecaptionskip}{2pt}
\caption{Hyperparameters of $\mathit{\Phi}$ for the case with experimental images.}
\label{tab:experimental_parameters}
\small
\begin{tabular}{ccccccc}
\toprule
Parameter & $w_1$ & $a_1$ & $c_1$ & $w_2$ & $a_2$ & $c_2$ \\
\midrule
Value & 0.05 & $3\times10^{-5}$ & $1\times10^{-5}$ & 0.04 & 0.85 & 0.05 \\
\bottomrule
\end{tabular}
\end{table}






As we can see from  {Fig.~\ref{fig:0465_u_comparison}}, the standard  {(unmasked)} DIC was difficult to converge and resulted in an over-iterated solution. This is due to the spurious displacement in the background region that is supposed to be removed by a pre-defined mask. By contrast, the proposed  {PF-DIC} does not rely on the mask and can converge quickly and robustly. In the current example, the PF-DIC solution converged at iteration 4 and provided a reasonable displacement field comparable to the original work of \cite{yang2019augmented}.  {This result also agrees with that of masked DIC. Both $u_x$ and $u_y$ show an excellent agreement, although the PF-DIC does not use any pre-defined mask. This again highlights the importance of a mask for handling cracks and removing unnecessary backgrounds in standard DIC.}
 
 {Regarding the damage measurement, we can see that the cracking area was successfully identified by the PF-DIC. In terms of crack tip location and crack area, the PF-DIC results can be considered  $90\%$ accurate for a crack contour corresponding to $d=0.5$. This is well acceptable for PF-DIC, given the diffusive nature of the PF crack representation. Additionally, due to the choice of $\mathit{\Phi}$, the damage field in this example does not distinguish the crack and the actual background, as expected. The white color in the PF-DIC displacement results indicates the area where the value of $d$ is relatively high. Ideally, high damage values should occur in regions where the available image information is insufficient for reliable displacement measurement. This consideration was incorporated into the selection of the hyperparameters reported in Table~\ref{tab:experimental_parameters}. A more systematic calibration procedure for these parameters will be presented later.}

\begin{figure}[htbp]
  \centering
\hspace*{-0.1\textwidth}
  \includegraphics[width=0.7\textwidth]{ 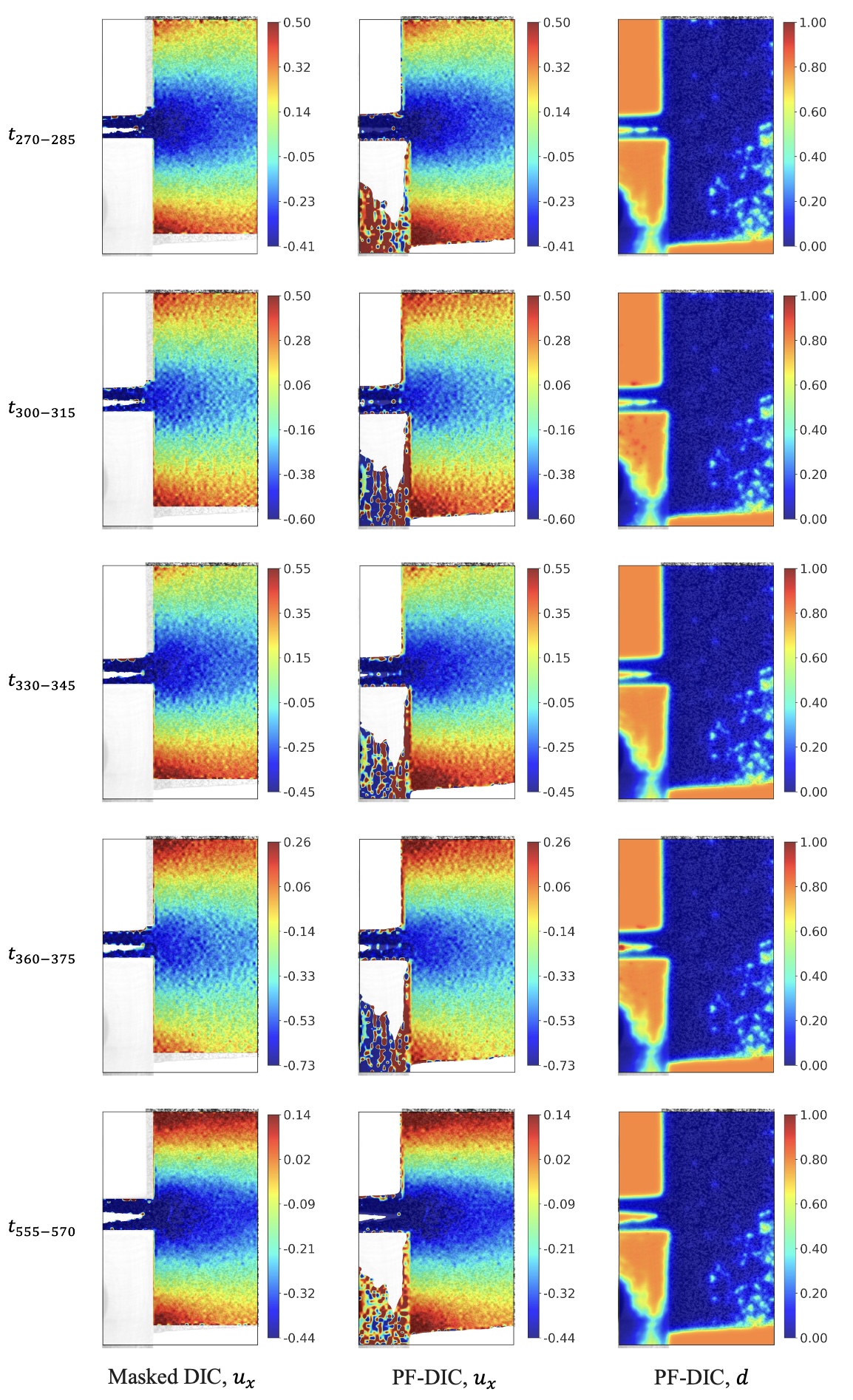}
  \caption{ {Masked DIC and PF-DIC results for unseen experimental images.}}
  \label{fig:exp_u}
\end{figure}

To further evaluate the performance of the PF-DIC with  {unseen} experimental images, we tested it on  {five} additional pairs of images captured at different time instances:  {$t_{270-285}$, $t_{300-315}$, $t_{330-345}$, $t_{360-375}$, and $t_{555-570}$}.  {Here, the hyperparameters in Table~\ref{tab:experimental_parameters} stayed unchanged and were calibrated only once using the previous images of $t_{465-480}$. {Fig.~\ref{fig:exp_u}} summarizes the results for the {five} image pairs using the masked DIC and (unmasked) PF-DIC. For the displacement comparisons, we focused on $u_x$. As we can see, the PF-DIC results again exhibit good agreement with those obtained with masked DIC, across a wide range of displacement that is different from the one used for calibration. This agreement confirmed the robustness of the PF-DIC for displacement measurements.}

\begin{figure}[htbp]
  \centering
    \hspace*{-0.1\textwidth}
  \includegraphics[width=0.7\textwidth]{ 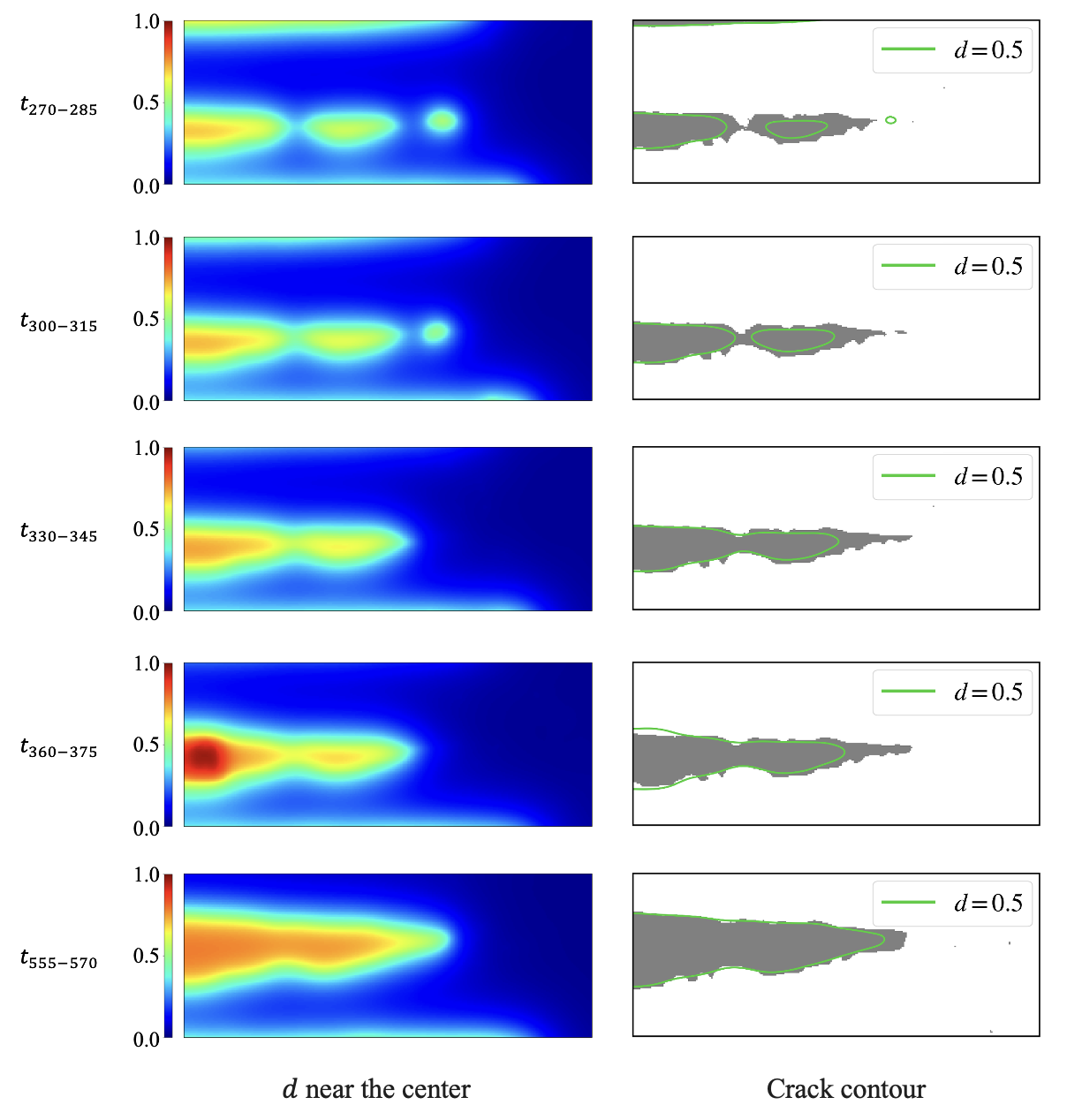}
  \caption{ {Verification of the crack contour in the PF-DIC.}}
  \label{fig:all_exp_contour}
\end{figure}

 {To verify the damage measurements obtained by the PF-DIC, we extracted the crack from the speckle images at the different time and compared it with the crack contour represented by $d=0.5$. As shown in Fig.~\ref{fig:all_exp_contour}, the evolution of the damage field $d$ can effectively represent the evolution of the actual crack during this experiment: the initial small cracks expanded and merged into a big one; and the overall advancement of the crack is well captured. The crack contour at different time also looks similar to the reference. For quantitative comparisons, we calculated the damaged area enclosed by the contour and its overlap with the reference region. As shown in Table~\ref{tab:crack_overlap_metrics}, the overall damage measurement achieves an accuracy of more than $80\%$ in most cases. Since this is for a contour around the midline of the diffusive crack boundary, the $80\%$  accuracy is considered well acceptable in our opinion.}

\FloatBarrier

\FloatBarrier

\FloatBarrier

\begin{table}[H]
\centering
\setlength{\abovecaptionskip}{2pt}
\captionsetup{font=small}

\caption{ {Crack-area and overlap comparison between the reference and PF-DIC results.}}
\label{tab:crack_overlap_metrics}

\small
\setlength{\tabcolsep}{14pt}
\renewcommand{\arraystretch}{1.25}

\begin{tabular}{ccccc}
\toprule
Image group & Method & Area & Overlap & Ratio to ref. \\
\midrule

\multirow{2}{*}{$t_{270\text{--}285}$}
& Ref.
& 3209
& --
& -- \\

& PF-DIC
& 2505
& 2324
& 0.72 \\[0.4em]

\multirow{2}{*}{$t_{300\text{--}315}$}
& Ref.
& 3640
& --
& -- \\

& PF-DIC
& 3172
& 3031
& 0.83 \\[0.4em]

\multirow{2}{*}{$t_{330\text{--}345}$}
& Ref.
& 4159
& --
& -- \\

& PF-DIC
& 3827
& 3610
& 0.87 \\[0.4em]

\multirow{2}{*}{$t_{360\text{--}375}$}
& Ref.
& 4591
& --
& -- \\

& PF-DIC
& 4661
& 4103
& 0.89 \\[0.4em]

\multirow{2}{*}{$t_{555\text{--}570}$}
& Ref.
& 7008
& --
& -- \\

& PF-DIC
& 7011
& 6679
& 0.95 \\

\bottomrule
\end{tabular}
\end{table}

 {In summary, the PF-DIC successfully identified both the displacement and damage fields. The unmasked PF-DIC displacement measurements achieved accuracy comparable to that of standard masked DIC, while the additional damage field effectively captured the evolution of the experimental cracks. These results validate the capability of PF-DIC to identify displacement and damage fields in previously unseen cracked specimens, provided that the method is calibrated against a baseline setup. }

\section{ {Parameter sensitivity analysis}}
 {This section presents a sensitivity analysis of the various hyperparameters in the PF-DIC, followed by a general discussion on the calibration procedure. We used the synthetic images in Section \ref{sec:crack_propagation} to illustrate the effects of different parameters, as the ground-truth results are given. The resulting observations and conclusions should be generally applicable to both synthetic and experimental images.}

\subsection{ {Influence of $l_c$}}
\begin{figure}[htbp]
  \centering
  \includegraphics[width=0.5\textwidth]{ 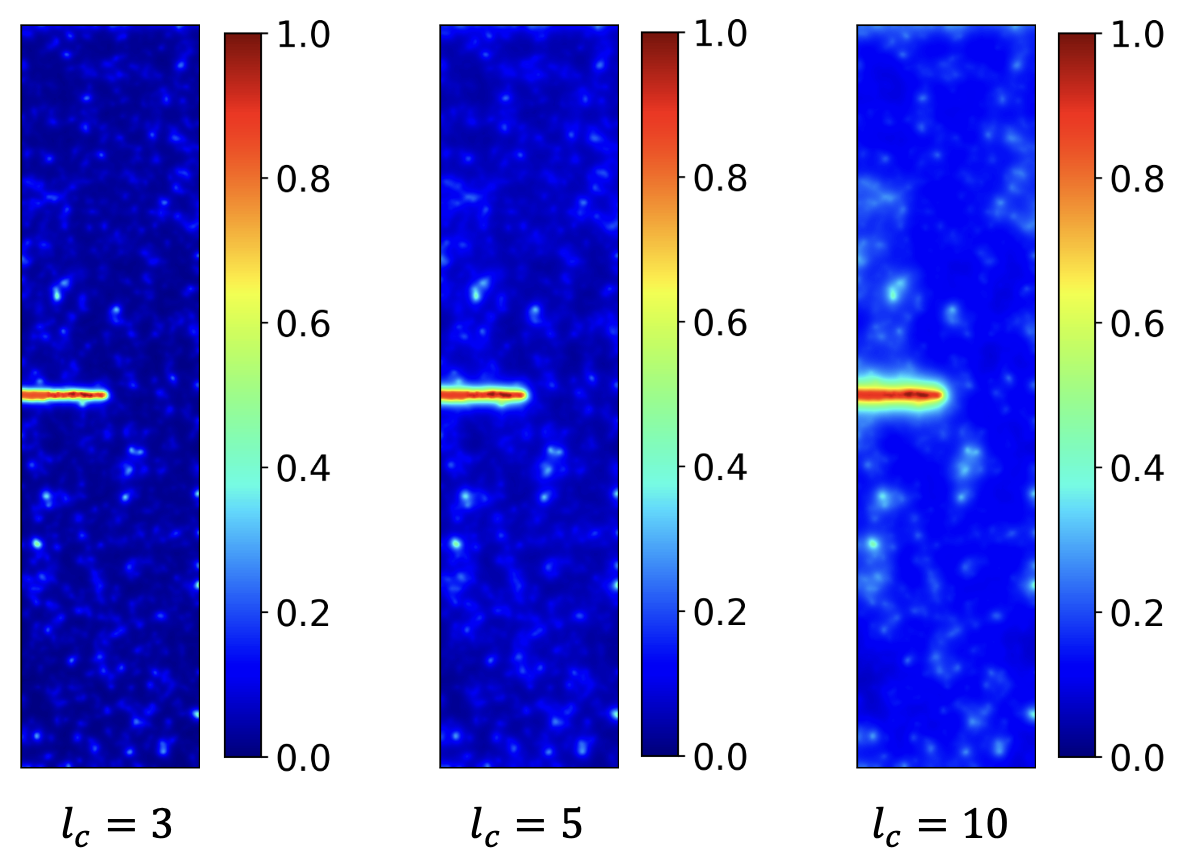}
  \caption{ {Damage fields obtained by PF-DIC using different $l_c$.}}
  \label{fig:lc_crack}
\end{figure}

 {Although the characteristic length $l_c=5$ was fixed in the examples, we first wanted to illustrate the effect of this parameter on the PF-DIC measurements. As mentioned earlier, this parameter is expected to control the size of the transition zone in the diffusive representation of cracks. To illustrate this effect, we conducted the PF-DIC for three different values of $l_c$, while the other settings stayed unchanged. As shown in Fig.~\ref{fig:lc_crack}, a large $l_c$ could increase the thickness of the crack significantly while reducing the contrast between damaged and undamaged regions. This is expected, as $l_c$ controls the contribution of the gradient of damage. In PF models, a general rule is that $l_c$ should be at least 3-5 times less than the mesh size \cite{Miehe2010,molnar20172d,molnar2020open}. $l_c=5$ is a reasonable choice, given the spatial discretization of the damage field ($h_d=1$) and the observations in Fig.~\ref{fig:lc_crack}. In general, we expect that $h_d$ is small enough to capture the length scale of the physical crack, and $h_d=1$ seems the smallest size we may use for a given speckle image.}

\FloatBarrier

\begin{figure}[htbp]
  \centering

  \begin{subfigure}[t]{0.25\textwidth}
    \centering
    \includegraphics[width=\linewidth]{ 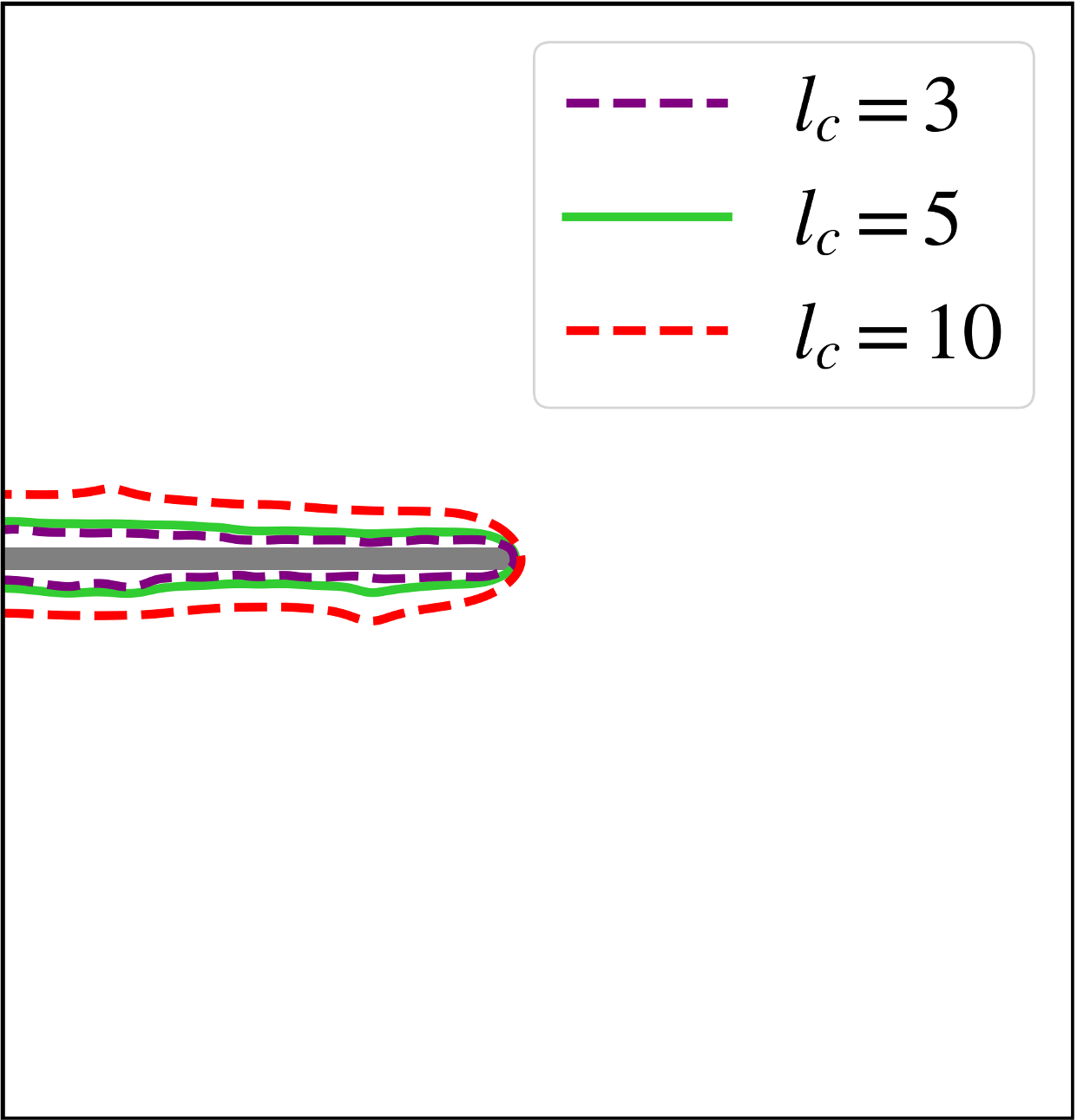}
    \caption{ {Crack thickness = 5 pixels}}
    \label{5contour}
  \end{subfigure}
  \hspace{0.02\textwidth}
  \begin{subfigure}[t]{0.25\textwidth}
    \centering
    \includegraphics[width=\linewidth]{ 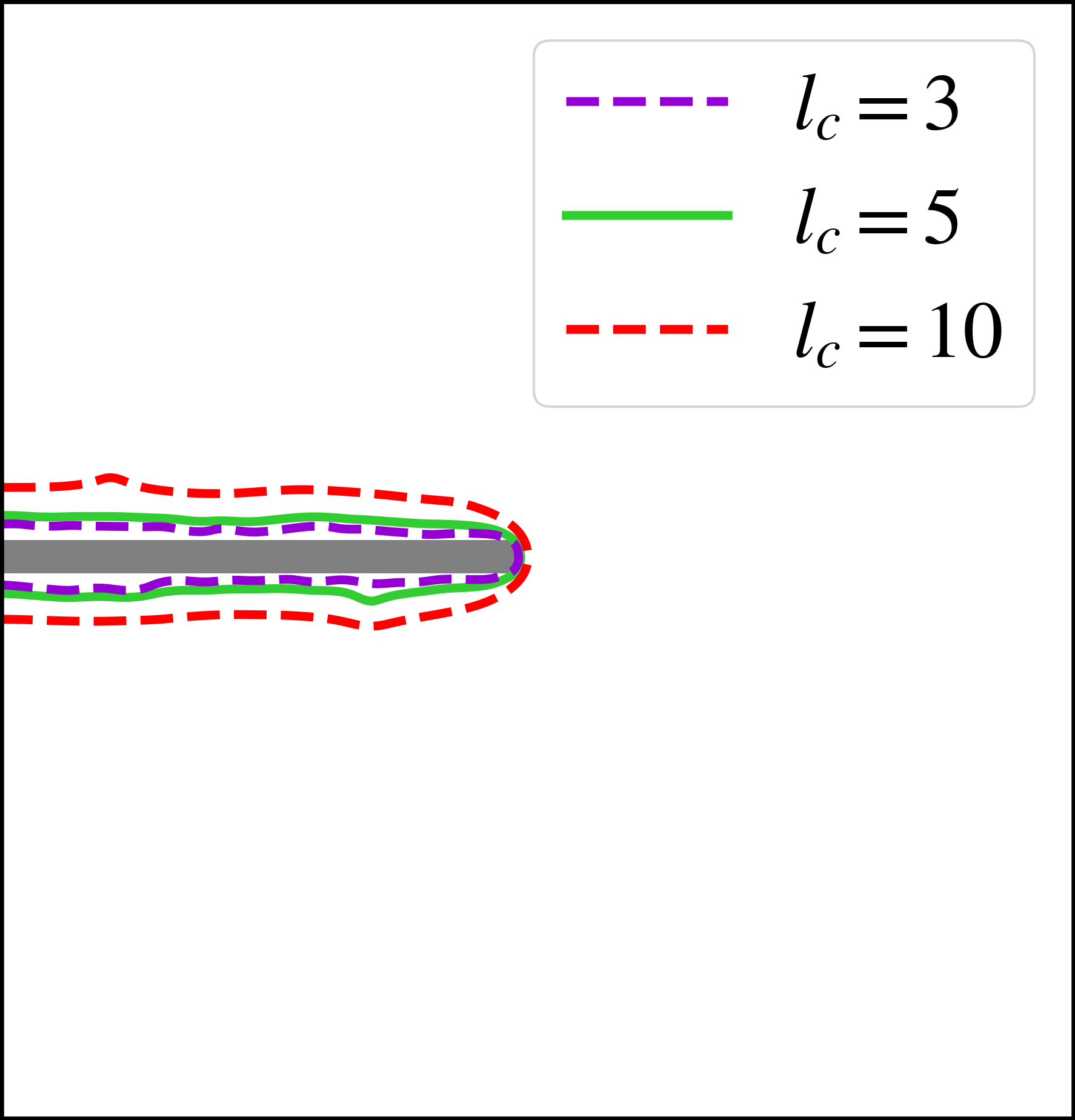}
    \caption{ {Crack thickness = 10 pixels}}
    \label{10contour}
  \end{subfigure}
  \hspace{0.02\textwidth}
  \begin{subfigure}[t]{0.25\textwidth}
    \centering
    \includegraphics[width=\linewidth]{ 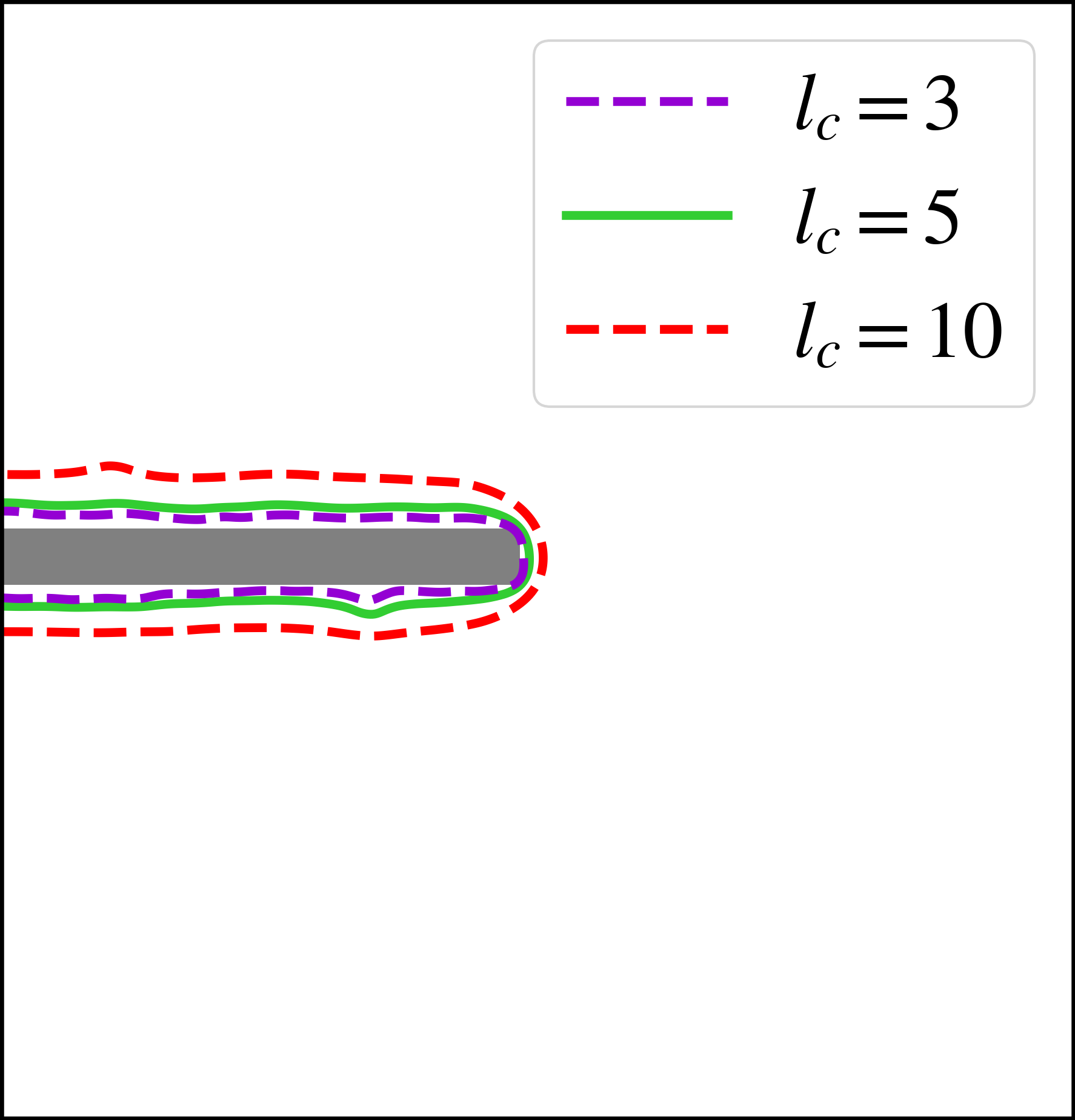}
    \caption{ {Crack thickness = 20 pixels}}
    \label{20contour}
  \end{subfigure}

  \caption{ {Damage contours for cracks of different sizes.}}
  \label{fig:different_lc}
\end{figure}

 {To verify if this choice depends on the size of the actual crack, we generated results for cracks of different thickness. We expected that a reasonable $l_c$ should always give the correct crack contour and tip location. Fig.~\ref{fig:different_lc} illustrates the contours of different $l_c$ for $d=0.5$ in the case of different crack thickness. As we can see, the contour given by $l_c=10$ is very sensitive to the actual crack size, and it overestimated the advancement of crack tip. By contrast, the results for $l_c=3,5$  are less sensitive to the actual crack size, with both the crack contour and crack tip location accurately identified. This observation supports our choice of $l_c$, and it can be generally fixed to a given value when the mesh resolution of $d$ is fixed.}

\FloatBarrier

\subsection{ {Influence of $w_1$, $a_1$, $c_1$}}
 {Next, we studied the hyperparameters used in the damage driving force energy $\mathit{\Phi}$, under the definition of Eq.~\eqref{eq:fracture_energy_w_background}. To isolate the effects of $w_1$, $a_1$, $c_1$, we set $w_2=0$ in this study. In this case, the damage driving force energy only contains  one function: $\mathit{\Phi}=\mathcal{A}_1
\!\left(
\boldsymbol{\varepsilon}, r,f
\right)=w_1\tanh (\frac{f\,r^2\,
\boldsymbol{\varepsilon}:\boldsymbol{\varepsilon}-a_1}{c_1})$.}  

\begin{figure}[htbp]
  \centering

  \begin{subfigure}[t]{0.6\textwidth}
    \centering
    \includegraphics[width=\linewidth]{ 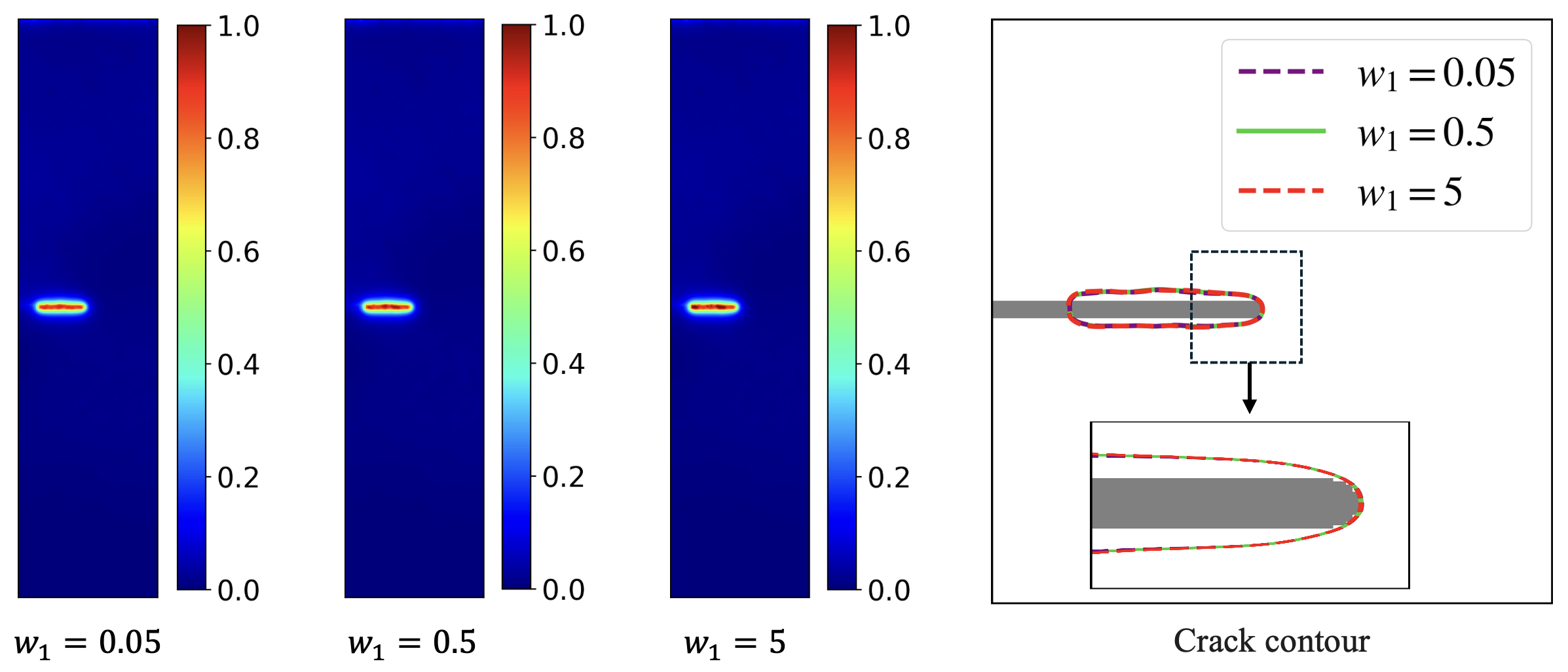}
    \caption{ {Damage fields obtained using different values of $w_1$}}
    \label{fig:w1}
  \end{subfigure}

  \vspace{0.3cm}

  \begin{subfigure}[t]{0.6\textwidth}
    \centering
    \includegraphics[width=\linewidth]{ 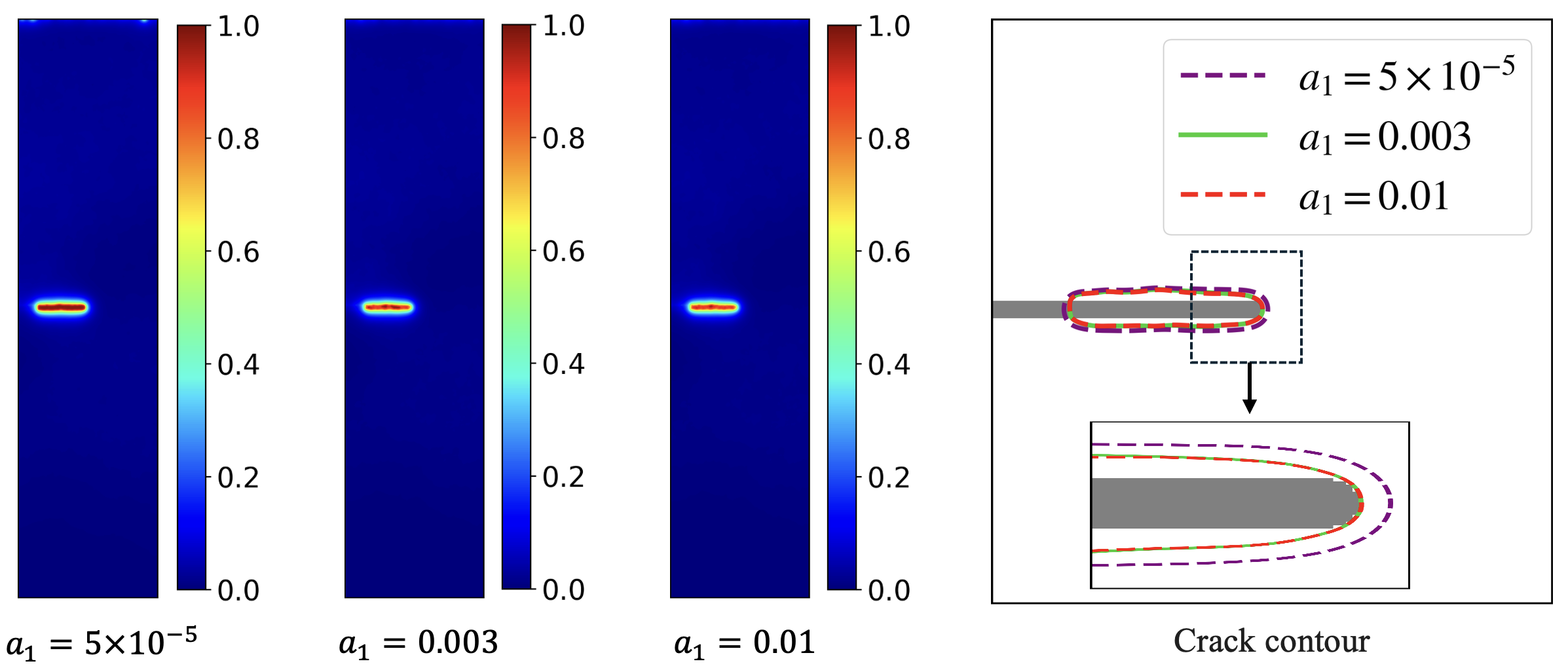}
    \caption{ {Damage fields obtained using different values of $a_1$}}
    \label{fig:a1}
  \end{subfigure}

  \vspace{0.3cm}

  \begin{subfigure}[t]{0.6\textwidth}
    \centering
    \includegraphics[width=\linewidth]{ 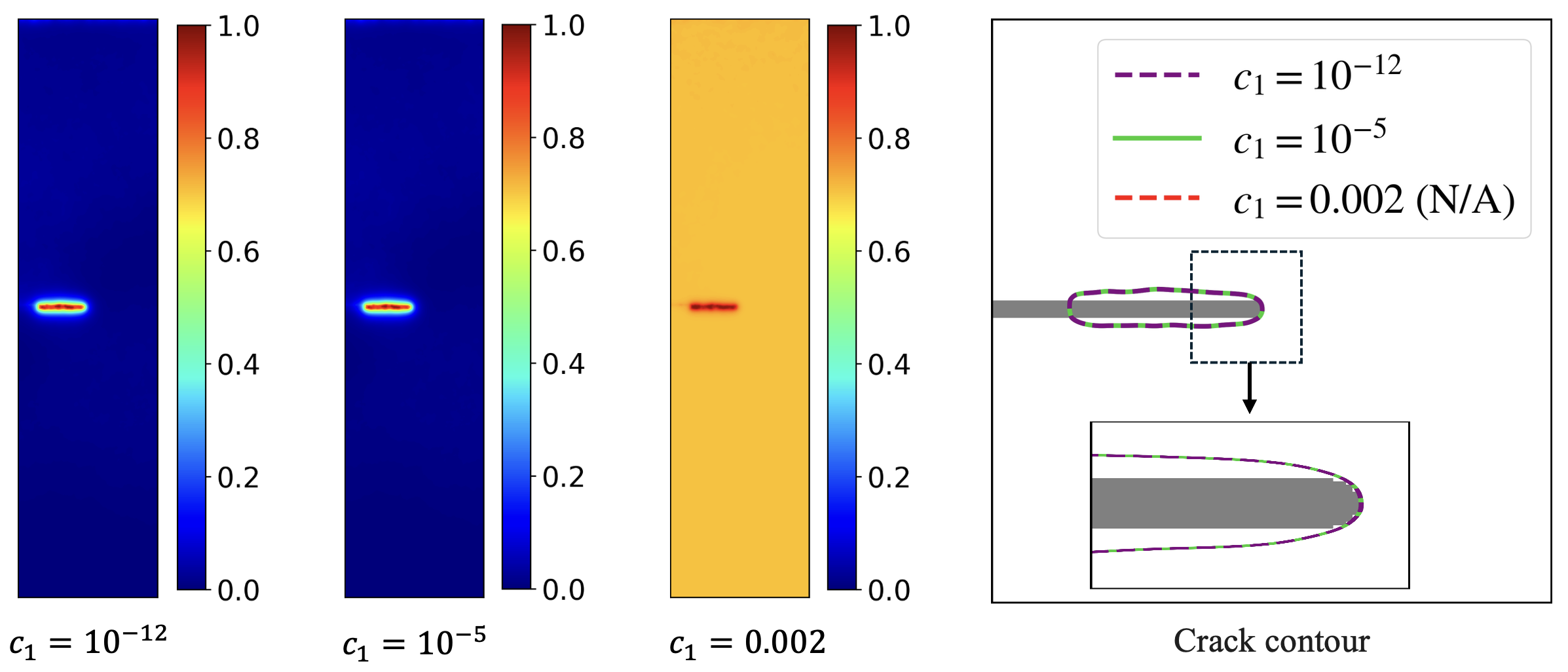}
    \caption{ {Damage fields obtained using different values of  $c_1$}}
    \label{fig:c1}
  \end{subfigure}

  \caption{ {Parametric study for $w_1$, $a_1$, and $c_1$.}}
  \label{fig:w1_a1_c1}
\end{figure}

 {We first studied the parameter $w_1$ by changing it to three different values: 0.05, 0.5, and 5, while fixing $a_1=0.003$ and $c_1=10^{-5}$. The corresponding damage results are illustrated in Fig.~\ref{fig:w1_a1_c1}(a). As we can see, the three values have very minor effects on the final crack contour, meaning the flexible choice of this parameter. In addition, we can see that the current choice of $\mathit{\Phi}$ can only identify the extended part of the crack. The initial crack in $f$ was not captured. This is expected, as the second term in Eq.~\eqref{eq:fracture_energy_w_background} is responsible for identifying existing cracks and it is set to 0 in the current study.}

 {The second target parameter is $a_1$. We can expect that this parameter can significantly influence the damage measurement, as it serves as a threshold in the damage function. The value of this parameter should be generally between 0 and $\max(f\,r^2\,
\boldsymbol{\varepsilon}:\boldsymbol{\varepsilon})$. To illustrate its effect, we conducted the analysis with three different values: $5\times10^{-5}$, 0.003, and 0.01, while fixing $w_1=0.5$ and $c_1=10^{-5}$. The results can be found in Fig.~\ref{fig:w1_a1_c1}(b). As we expected, the final damage field was significantly affected by this parameter. We can observe that a small value would overestimate the damage. As this parameter is gradually increased, the identified crack area decreased and the contour may converge toward the actual crack. Again, this is consistent with our expectation.}

 {Regarding the parameter $c_1$, it controls the overall window size in the hyperbolic tangent function. The value of this parameter can be as small as possible, and it should not affect the final damage result, if $a_1$ is well calibrated. Again, we tested three different values for $c_1$: $10^{-12}$, $10^{-5}$, and ${0.002}$, as shown in Fig.~\ref{fig:w1_a1_c1}(c). We can see that the results for $c_{1}=10^{-12}$ and $c_{1}=10^{-5}$ do not exhibit any noticeable difference, and they both successfully captured the extended part of the crack. For a large value of $c_1$, the damage field was not well identified. This makes the choice of this parameter straightforward, since a small value is sufficient.}

\FloatBarrier

\subsection{ {Influence of $w_2$, $a_2$, $c_2$}}
\begin{figure}[htbp]
  \centering

  \begin{subfigure}[t]{0.6\textwidth}
    \centering
    \includegraphics[width=\linewidth]{ 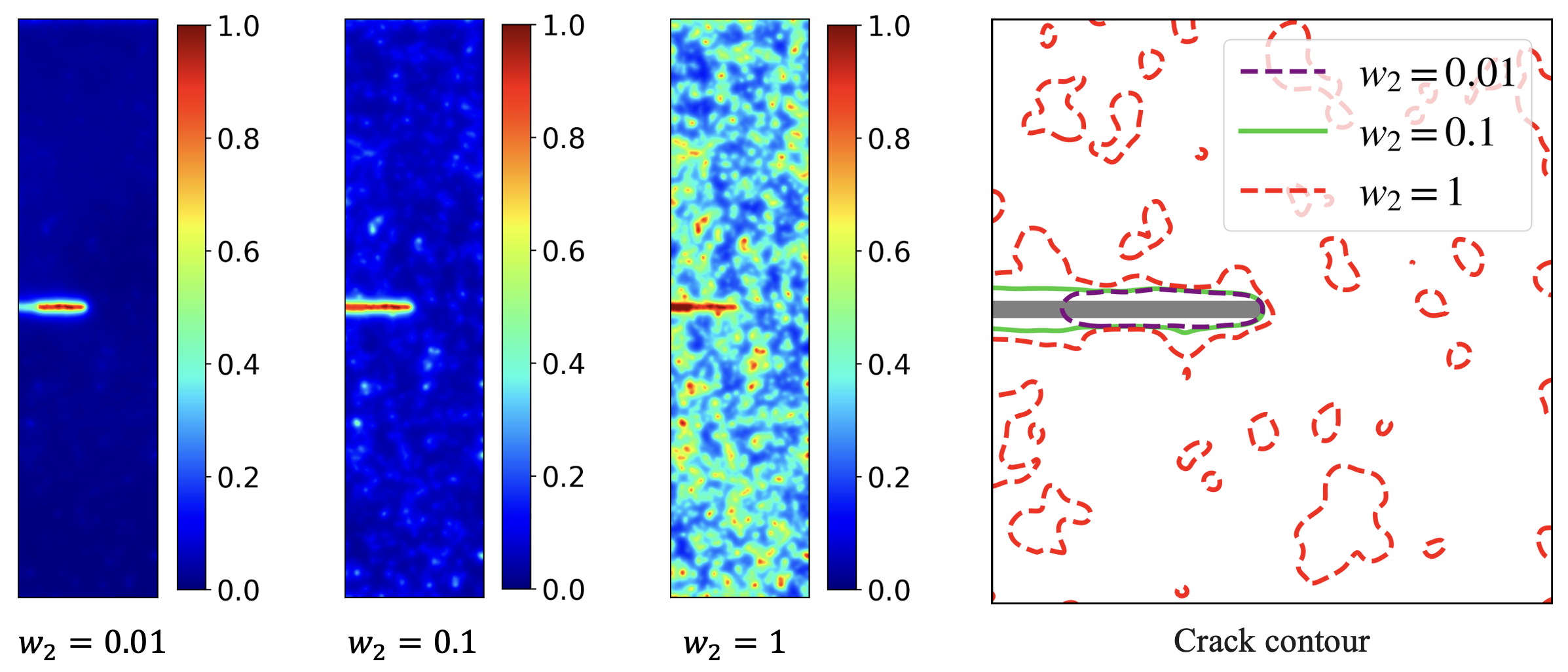}
    \caption{ {Damage fields obtained using different values of $w_2$}}
    \label{fig:w2}
  \end{subfigure}

  \vspace{0.3cm}

  \begin{subfigure}[t]{0.6\textwidth}
    \centering
    \includegraphics[width=\linewidth]{ 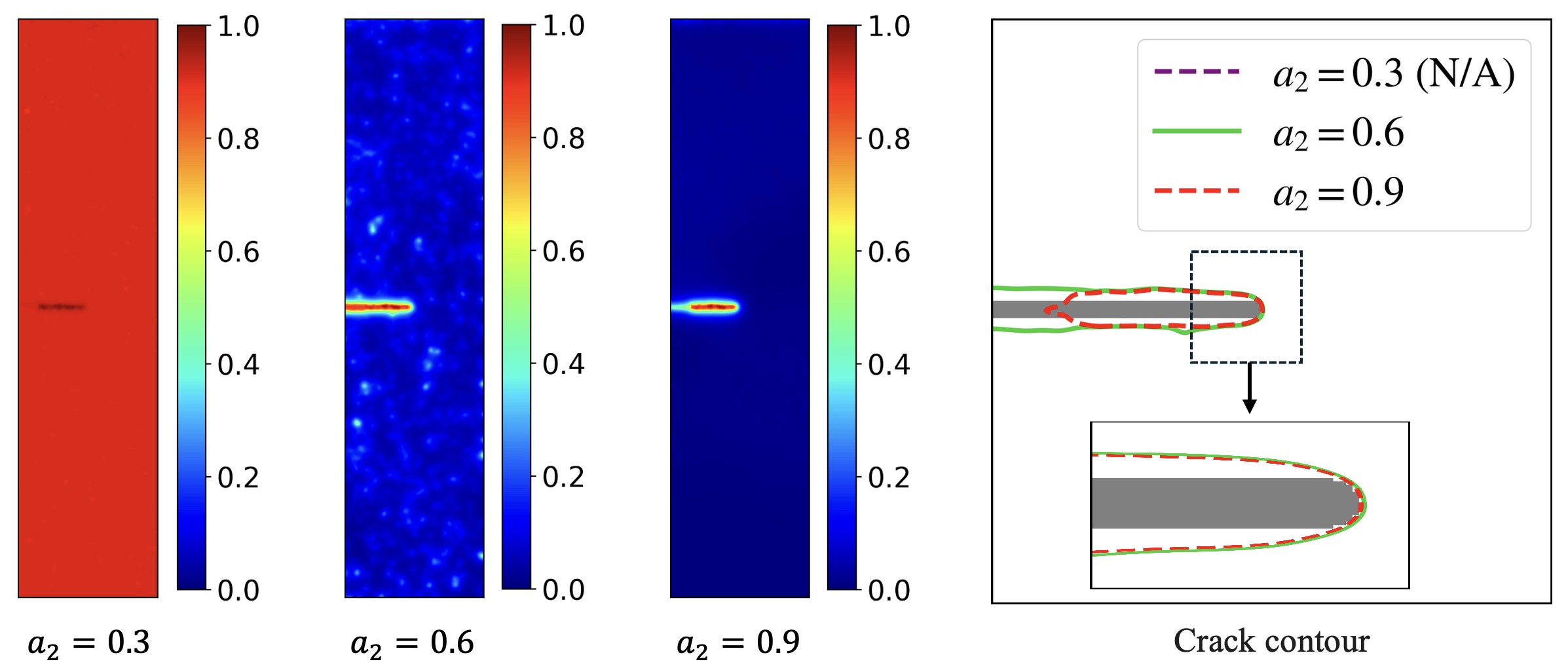}
    \caption{ {Damage fields obtained using different values of $a_2$}}
    \label{fig:a2}
  \end{subfigure}

  \vspace{0.3cm}

  \begin{subfigure}[t]{0.6\textwidth}
    \centering
    \includegraphics[width=\linewidth]{ 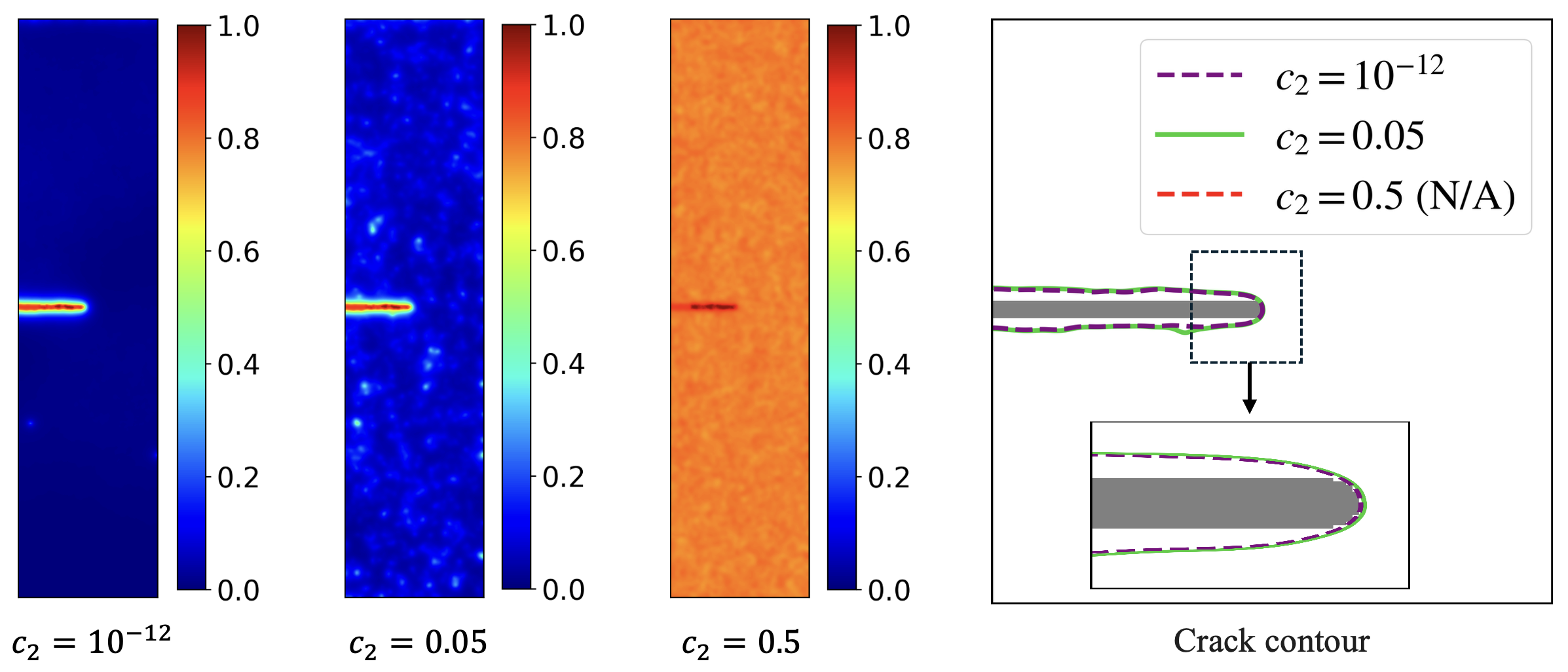}
    \caption{ {Damage fields obtained using different values of $c_2$}}
    \label{fig:c2}
  \end{subfigure}

  \caption{ {Parametric study for $w_2$, $a_2$, and $c_2$.}}
  \label{fig:w2_a2_c2}
\end{figure}

 {Similarly to the above, we studied the hyperparameters in the second term of Eq.~\eqref{eq:fracture_energy_w_background}. Here, we fixed the values of the previous parameters with $w_1=0.5$, $a_1=0.003$, and $c_1=10^{-5}$.  Fig.~\ref{fig:w2_a2_c2} illustrates the results for different values of the target parameters: $w_2$, $a_2$, $c_2$. As we can see, $w_2$ can magnify the overall damage value in the ROI. If an appropriate value is selected, this parameter only magnifies the region containing the initial crack. This can be seen by comparing the first two figures in Fig.~\ref{fig:w2_a2_c2}(a). $w_2=0.1$ is a good choice in this case. Regarding the  parameter $a_2$, we can expect it to have an effect similar to that of $a_1$: a small value would lead to overestimated damage. As shown in Fig.~\ref{fig:w2_a2_c2}(b), a small value $a_2=0.3$ did not provide any meaningful damage measurement. As the value increased, we started seeing the expected damage measurement. A large value could underestimate the initial crack. Therefore, this is also a critical value that needs to be carefully selected, typically within the range from 0 to $\max(f)$. Finally, for the parameter $c_2$, it plays a role similar to that of $c_1$, and a small value should be sufficient. As confirmed by Fig.~\ref{fig:w2_a2_c2}(c), both 0.05 and $10^{-12}$ would be acceptable in the current study.}

\FloatBarrier

\subsection{ {Discussion on the calibration of parameters}}
 {Based on the previous sensitivity analysis, we can propose a general procedure to calibrate the hyperparameters. Assuming a pair (or a set) of baseline speckle images of a cracked specimen has already been selected, we may use the following strategy to adjust the parameters in the damage driving force energy, based on the definition of Eq.~\eqref{eq:fracture_energy_w_background}.}
 {\begin{enumerate}
    \item Give some initial values to the six parameters: $w_1$, $a_1$, $c_1$, $w_2$, $a_2$, $c_2$. A simple choice can be $w_1=1$, $a_1=\max(f\,r^2\,
\boldsymbol{\varepsilon}:\boldsymbol{\varepsilon})/2$, $c_1=10^{-12}$,  $w_2=0$, $a_2=\max(f)/2$, $c_2=10^{-12}$, where $\max(f\,r^2\,
\boldsymbol{\varepsilon}:\boldsymbol{\varepsilon})$ can be known after the first iteration for displacement measurements in the PF-DIC. 
\item  Run the PF-DIC and check the damage measurement and potentially extract the crack contour for verification. If the damage measurement does not match the target crack, we can first adjust the parameter $a_1\in[0,\ \max(f\,r^2\,
\boldsymbol{\varepsilon}:\boldsymbol{\varepsilon})]$ according to whether the current values result in overestimated or underestimated the damage. If needed, adjust $w_1$  based on the contrast between damaged and undamaged regions.
\item  Gradually increase $w_2$, and adjust the value of $a_2\in[0,\ \max(f)]$ by verifying the resulting damage measurement in the PF-DIC.
\end{enumerate}}
 {In the above procedure, $c_1$ and $c_2$ are simply set to a small value of $10^{-12}$, so that they do not require further calibration. This procedure should be general and can be adapted  according to the chosen damage driving force energy and the judgment of users on how and when the damage should be detected. In general, the hyperparameters need to be recalibrated for each new experimental setup.}

\section{Conclusions}
This work proposed a PF-DIC framework that combined the optical flow conservation with the PF approach for integrated displacement and damage measurements. The PF-DIC has several advantages over standard DIC when dealing with cracked samples. First, unlike the standard DIC, the PF-DIC eliminates the need for a pre-defined mask, which can be difficult to define for complex crack morphologies. Second, the  {PF-DIC can effectively capture crack propagation while achieving displacement measurement accuracy comparable to that of masked DIC, with notable improvements observed in some cases}. Third, the PF-DIC enables automatic damage detection and diagnostic capabilities that can selectively identify critical cracks in defective materials and structures.  

The numerical examples have demonstrated the advantages of PF-DIC using both synthetic and experimental speckle images of fractured samples. These examples are representative in material characterization, including tensile and bending loading conditions and heterogeneous displacement.  {The PF-DIC outperformed standard unmasked DIC in all examples and produced displacement results comparable to those of masked DIC, thereby validating the displacement measurement capability of PF-DIC
}.

 {The scope of this work is to demonstrate the feasibility of the proposed PF-DIC formulation. The results obtained from both synthetic and experimental images are encouraging and suggest its potential applicability to various fractured specimens. In the future, we plan to apply the PF-DIC for automatic extraction of fracture properties of materials, by including additional physical parameters, such as energy release rates, in the PF-DIC formulation.}

\section*{Acknowledgments}

The authors would like to acknowledge the support of University of Maryland Baltimore County through the startup fund. 

\section*{ {Data availability}}
 {Data will be made available on request.}

\bibliographystyle{elsarticle-num} 
\bibliography{main}

\end{document}